**Highlights:**

1. Novel fully coupled aero–hydro–mooring model of six hybrid wind–wave configurations
2. WEC reduced platform response in surge and pitch, enhancing overall stability
3. Hybridization reduced WEC LCOE by up to 83%, while FWT had minimal change
4. Power fluctuations reduced by ~50%, enhancing stability and lowering storage needs
5. Synergy-based framework reveals mutualism benefits in hybrid offshore renewables

# Advancing Offshore Renewable Energy: Techno-Economic and Dynamic Performance of Hybrid Wind-Wave Systems

Alaa Ahmed[a,*], Maha N. Haji[a,b]

[a]Cornell University
[b]University of Michigan,

ABSTRACT

Offshore wind and wave energy offer high energy density and availability. While offshore wind has matured significantly, wave energy remains costly and under development. Integrating both technologies into a hybrid system can enhance power generation, stabilize output, and reduce costs. This study explores the benefits of combining an offshore floating wind turbine with the two-body heaving point absorber wave energy converter, Reference Model 3 (RM3). Six configurations are analyzed: RM3 integrated with the National Renewable Energy Laboratory 5 MW and the International Energy Agency 15 MW wind turbines, each tested on both spar and semi-submersible platforms. The analysis examines dynamic response, mooring loads, and power production under varying environmental conditions, considering the influence of the wave energy converter's float motion and an optional reaction plate. Results indicate that the reaction plate improves damping for the spar platform, enhancing wave energy absorption and power output. A comparative analysis indicates that integrating the wave energy converter reduces its levelized cost of energy by 15-83%, while leaving the wind turbine's levelized cost of energy unaffected. Hybridization significantly reduces power fluctuations by 50%, reduces the levelized cost of energy with the 5 MW wind turbine, and slightly increases it with the 15 MW wind turbine. The results highlight a mutualistic relationship between the wave energy converter and the offshore wind turbine, where the former benefits substantially while the latter experiences slight improvements or negligible effects. Additional findings quantify hydrodynamic interactions, mooring performance, and economic feasibility. This research provides insights into optimizing hybrid offshore renewable systems, demonstrating their potential to lower costs and support sustainable energy solutions.

## 1. Introduction

The development of domestic renewable energy systems is essential for strengthening U.S. energy security and building resilient coastal infrastructure. Among various renewable energy technologies, wind and solar energy have been widely adopted for power generation Zhu, Chen, Wang, Wang, Fang, Du and Gao (2024); Khan, Yu and Waseem (2022). However, these sources are often intermittent and subject to seasonal variations, which can limit their reliability and efficiency Dabar, Awaleh, Waberi, Ghiasirad, Adan, Ahmed, Nasser, Juangsa, Guirreh, Abdillahi et al. (2024); Garcia and Oliva (2023). While wind energy offers significant potential, its variability and dependence on weather conditions make it less stable and sustainable for continuous energy supply Ibáñez-Rioja, Järvinen, Puranen, Kosonen, Ruuskanen, Hynynen, Ahola and Kauranen (2023); Tang, Rehme and Cerin (2022). In contrast, wave energy, derived from ocean waves and harvested by wave energy converters (WECs), presents a more consistent and reliable energy source Zou, Robertson, Roach, Mundon, Rosenberg and Penalba (2024), as it is available year-round and is less susceptible to seasonal fluctuations, making it particularly suitable to complement wind energy for powering energy-intensive processes. Exploring the integration of a floating wind turbine (FWT) and a WEC into a single system, named a hybrid wind-wave system, has been the interest of many researchers. Harnessing wind and wave energy through hybrid wind-wave systems could offer a sustainable and cost-effective alternative to conventional energy sources, addressing rising electricity demand Gao, Khan, Sergiienko, Ertugrul, Hemer, Negnevitsky and Ding (2022); Li, Luo, Zhou, Kumar, Guo, Law and Cao (2022). Additionally, this integration represents a groundbreaking advancement in renewable energy technology. These systems enhance energy capture, improve efficiency, and lower costs by integrating complementary renewable resources Kluger, Slocum and Sapsis (2017); Kluger, Haji and Slocum (2023); Moghimi, Derakhshan and Motawej (2020). Beyond energy and economic benefits, hybrid wind-wave systems exhibit remarkable scalability and adaptability, making them well-suited for a variety of geographical locations and environmental conditions Moghimi et al. (2020).

*Corresponding author, Email address: alaa.ahmed@cornell.edu
ORCID(s):

Hybridization allows for infrastructure cost-sharing and creates pathways to enhance the economic feasibility, ultimately contributing to a reduction in the levelized cost of energy (LCOE) He, Zhao, Liu, He and Luan (2024). For instance, Kluger et al. (2023) Kluger et al. (2023) found that the LCOE for the WEC decreases by 43% when integrated with offshore wind farms. Another key advantage of hybrid wind-wave energy systems lies in their capacity to address the variability challenges faced by standalone wind or wave systems, thereby providing a steadier and more dependable power output Gao, Bechlenberg, Jayawardhana, Ertugrul, Vakis and Ding (2024). According to Gao et al. (2024) Gao et al. (2024), these hybrid configurations can lower power storage needs by 20% and energy storage requirements by 35% when compared to wind-only systems, emphasizing their efficiency in both cost and resource utilization. Kluger et al. (2023) Kluger et al. (2023) also found that integrating wave energy with offshore wind farms leads to a 15% improvement in output energy smoothness, a 6% reduction in energy curtailment, and a 2% increase in overall power efficiency relative to standalone wind setups. Collectively, these performance improvements contribute to lower operational expenses and enhance the economic competitiveness of renewable energy technologies Kluger et al. (2023); Feng, Lin, Gu, Li, Chen, Liu and Sun (2023). However, the aforementioned studies analyzed the wind and wave systems separately, without accounting for their dynamic interactions. This gap underscores the need for integrated modeling approaches that capture full system dynamic interactions, enabling more accurate techno-economic assessments.

Recent studies have considered the dynamics of the hybrid system and confirmed that such a system not only enhances power absorption but also improves the dynamic behavior and structural stability of the overall platform. For example, the integration of Salter's Duck WECs with semi-submersible platforms was shown to reduce pitch up to 14.6% and heave responses by 12.34%, leading to more stable wind turbine (WT) operation and further smoothing of power output He et al. (2024). Moreover, Xu et al. (2024) Xu, Ji, Xu and Chen (2024) investigated the dynamic response and power performance of a combined semi-submersible FWT and WEC array using a 3D potential flow theory and numerical modeling. The results show that integrating WECs improves system stability and increases total power output by 2.5%–6.5%, with optimal performance achieved through tailored power takeoff (PTO) damping and conical-bottom WEC float designs. Despite increasing interest, these previous studies rely on certain steady wind and regular wave conditions, limiting their relevance to real-world dynamics. Neshat et al. (2024) Neshat, Sergiienko, Nezhad, da Silva, Amini, Marsooli, Garcia and Mirjalili (2024) proposed a novel adaptive chaotic multi-objective grey wolf optimizer (DMOGWA) to enhance the performance of hybrid wave-wind energy systems by optimizing the trade-off between maximizing WEC power output and minimizing WT nacelle acceleration. The method demonstrated superior convergence speed and solution quality, achieving a 138.5% average power boost and a 41% reduction in acceleration compared to other optimization techniques. Furthermore, the synergy between wind and wave structures can enhance platform stability, reducing the need for costly stabilizing mechanisms. Feng et al. (2023) Feng et al. (2023) successfully demonstrated that hybrid platforms maintain a tilting angle below 10°, while Meng et al. (2023) Meng, Sergiienko, Ding, Zhou, Da Silva, Cazzolato and Li (2023) proposed a control model capable of generating radiated waves equal in magnitude and opposite in direction to the excitation waves, effectively canceling their destabilization effects. Furthermore, research remains underdeveloped with many limitations in experimental modeling. Some studies oversimplify WT representation as a fixed mass above water, neglecting hydrodynamic-aerodynamic interactions Zhao, Stansby, Liao and Li (2024); Wei, Cao, Cao, Han, Shi and Ji (2025), while others exclude wind excitation in wave tank tests, failing to account for aerodynamic effects Zhu, Shi, Zheng, Michele and Cao (2023). Nataliia et al. Sergiienko, Xue, da Silva, Ding and Cazzolato (2025) reviewed 178 cases in the study of wind–wave systems. Their analysis revealed that many studies rely on small-scale or outdated WTs, and approximately 40% neglect turbine dynamics. Furthermore, most studies do not consider realistic site conditions, and techno-economic assessments remain uncertain or have not been performed. These limitations necessitate a comprehensive framework that captures fully coupled aero-hydro-mooring interactions to realistically evaluate the performance of hybrid wind-wave energy systems in terms of dynamic response, power output, and economic viability.

A key question, then, is how WECs and FWTs influence one another when combined. This interaction can be evaluated through the lens of system synergy, where the performance of the integrated system differs from the sum of its parts. This approach emphasizes whether the combined system delivers benefits beyond those of standalone configurations. The relationship between the FWT and WEC components in a hybrid system can take various forms depending on their performance interactions Schupp, Bocci, Depellegrin, Kafas, Kyriazi, Lukic, Schultz-Zehden, Krause, Onyango and Buck (2019). Mutualism occurs when both devices gain, commensalism when one benefits while the other remains unaffected, and parasitism when one is slightly compromised Garcia-Sanz (2024). The degree to which the hybrid system is advantageous depends on whether the combined system outperforms the standalone configurations. Determining the degree of synergy—whether positive, neutral, or negative—provides a practical framework for assessing the viability of hybrid systems and justifying their role as offshore renewable investments.

This study introduces a novel integrated approach to evaluating hybrid wind-wave systems, combining a newly developed fully coupled aero-hydro-mooring dynamic modeling framework with a synergy-based assessment of component interactions. The modeling framework couples the open-source software tools WEC-Sim (Wave Energy Converter SIMulator) Ruehl, Keester, Forbush, Grasberger, Husain, Leon, Ogden and Shabara (2024), MOST (MATLAB for Offshore floating wind turbine Simulation Tool) Cottura, Caradonna, Ghigo, Novo, Bracco and Mattiazzo (2021), and MoorDyn Hall (2020). The coupled framework is used to evaluate the feasibility of hybrid wind-wave energy systems across multiple scales. Specifically, it explores configurations that combine the National Renewable Energy Laboratory (NREL) 5 MW reference WT Jonkman, Butterfield, Musial and Scott (2009) and the International Energy Agency (IEA) 15 MW WT Gaertner, Rinker, Sethuraman, Zahle, Anderson, Barter, Abbas, Meng, Bortolotti, Skrzypinski, Scott, Feil, Bredmose, Dykes, Shields, Allen and Viselli (2020) with the Reference Model 3 (RM3) two-body point absorber WEC Neary, Lawson, Previsic, Copping, Hallett, Labonte, Rieks and Murray (2014). Both spar-type Jonkman (2010) and semi-submersible Robertson, Jonkman, Masciola, Song, Goupee, Coulling and Luan (2014) floating platforms are examined to evaluate system performance and feasibility.

The analysis numerically examines the dynamic response of the WEC, WT, and mooring platform across six hybrid configurations and compares them against their corresponding standalone systems under regular and irregular wave conditions with turbulent wind field. In addition, the annual energy production (AEP) of each configuration is evaluated while fully accounting for coupled aero-hydro-mooring dynamic interaction. The analysis also includes a techno-economic assessment and evaluation of hybridization effects on selected metrics, LCOE, and power coefficient of variation ($P_{CV}$), highlighting the advantages of hybrid wind-wave systems. Finally, the study investigates the degree of synergy between the FWT and the WEC, identifying whether integration yields mutualistic, commensal, or parasitic relationships.

Overall, the paper emphasizes the dynamic response, power generation, $P_{CV}$, and LCOE of hybrid wind-wave systems across different configurations, comparing them to standalone systems at one location: offshore of Humboldt Bay, California. A second site, offshore of Kodiak, Alaska, is used to assess standalone systems and illustrate the role of site selection in hybrid wind–wave performance. All numerical simulations are conducted using the integrated framework (WEC-Sim-MOST-MoorDyn). To the best of the authors' knowledge, this study is the first to couple hydrodynamics, aerodynamics, and mooring dynamics in the analysis of integrated hybrid wind-wave systems and while also applying a synergy-based framework to assess their economic viability.

The paper is organized as follows: Section 2 present the wave and wind energy resources. Section 3 describes different configurations of the hybrid wind-wave systems studied. System synergy is explained in Section 4. The methodology and simulation procedures are detailed in Section 5, followed by a presentation and discussion of the results in Section 6. Finally, conclusions are drawn in Section 7.

## 2. Resource Assessment

### 2.1. Wave energy resource

The potential WEC deployment site considered in this work to analyze the performance of all systems is located offshore of Humboldt Bay, California (referred to as the Humboldt Site). The Humboldt Site Dallman and Neary (2014) is located at 40.8418° N, 124.2477° W with a reported depth up to 456 m. The joint probability distributions (JPD), shown in Figure 1, define the significant wave height ($H_{m0}$) and energy period ($T_e$) for the year 2017, based on measurements from buoy station 46022 National Data Buoy Center (2024b), owned and maintained by the National Data Buoy Center (NDBC). The most commonly occurring sea state has a significant wave height of 1.75 m and an energy period of 6.5 s and occurs for 768 hr annually.

### 2.2. Wind energy resource

The wind resource assessment for the Humboldt Site is conducted using TurbSim Jonkman et al. (2014), a stochastic, full-field turbulence simulator primarily designed for use with InflowWind-based simulation tools. TurbSim generates synthetic wind flow fields that replicate realistic atmospheric turbulence conditions, enabling accurate modeling of inflow conditions for WT simulations. To account for variations in wind shear and turbulence intensity at different heights, the wind field is adjusted according to the hub height of each WT, ensuring that the inflow conditions are representative of the wind experienced at the hub height. This approach allows for a consistent evaluation of turbine performance across different configurations. The wind speed data was obtained from the measurements from buoy station 46022 National Data Buoy Center (2024b) for the year 2017. The wind rose for the Humboldt Site, presented in Figure 2, illustrates the predominant wind directions and their frequency distribution, providing critical insights into the wind characteristics at the target location. As shown, a dominant wind speed range between 5 m/s and 10 m/s is observed at the site.

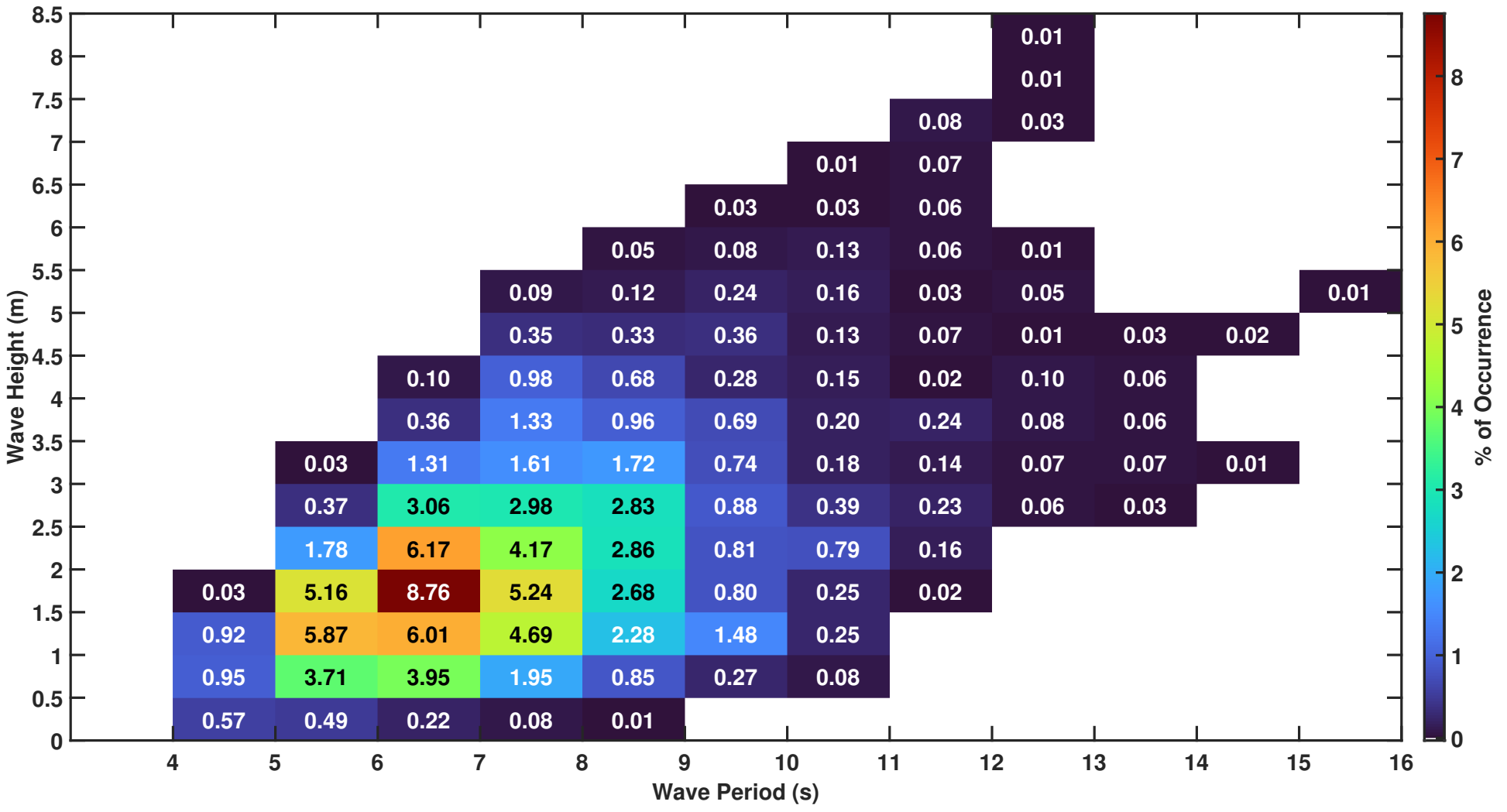

**Figure 1:** JPD of sea states at the Humboldt Site.

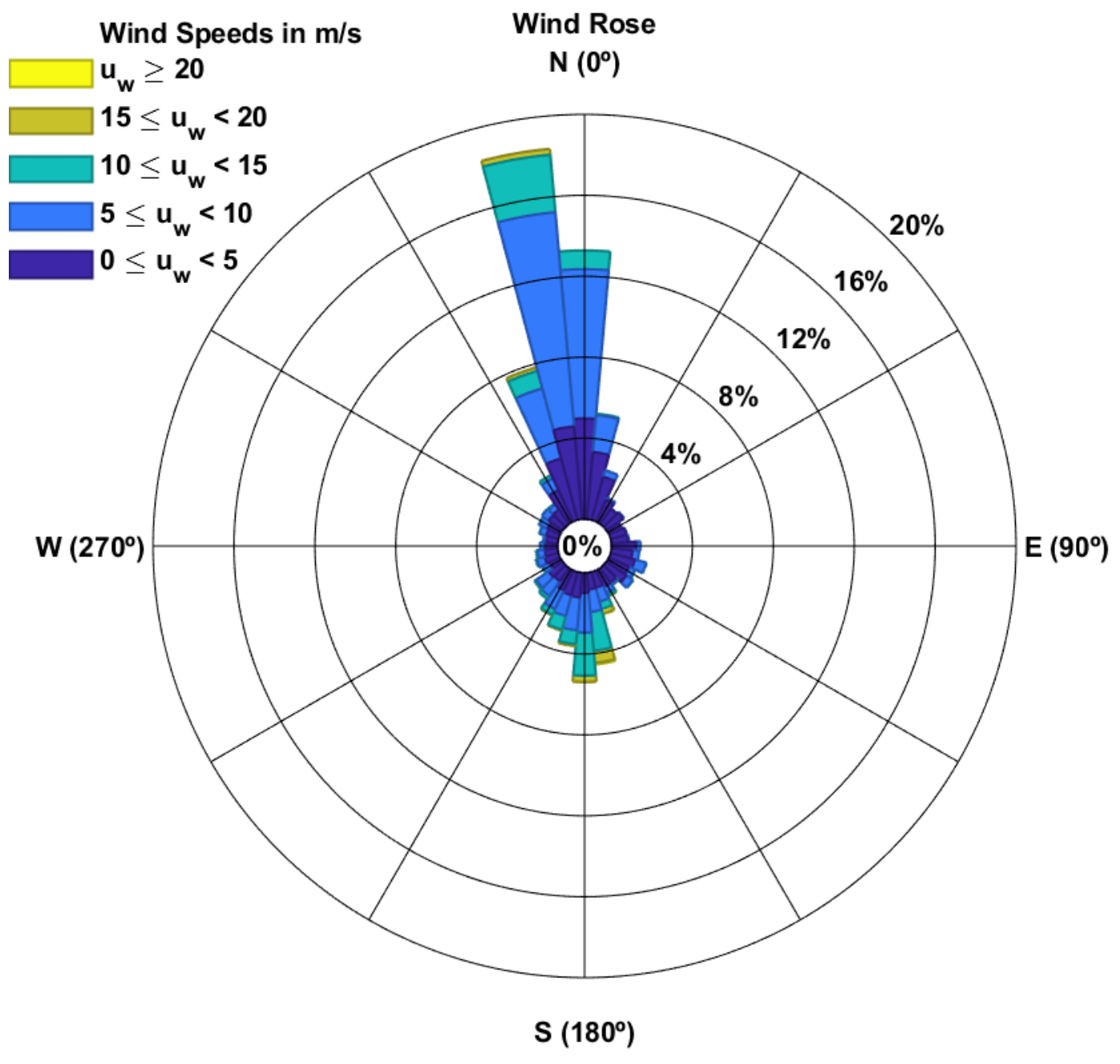

**Figure 2:** Wind rose at the Humboldt Site.

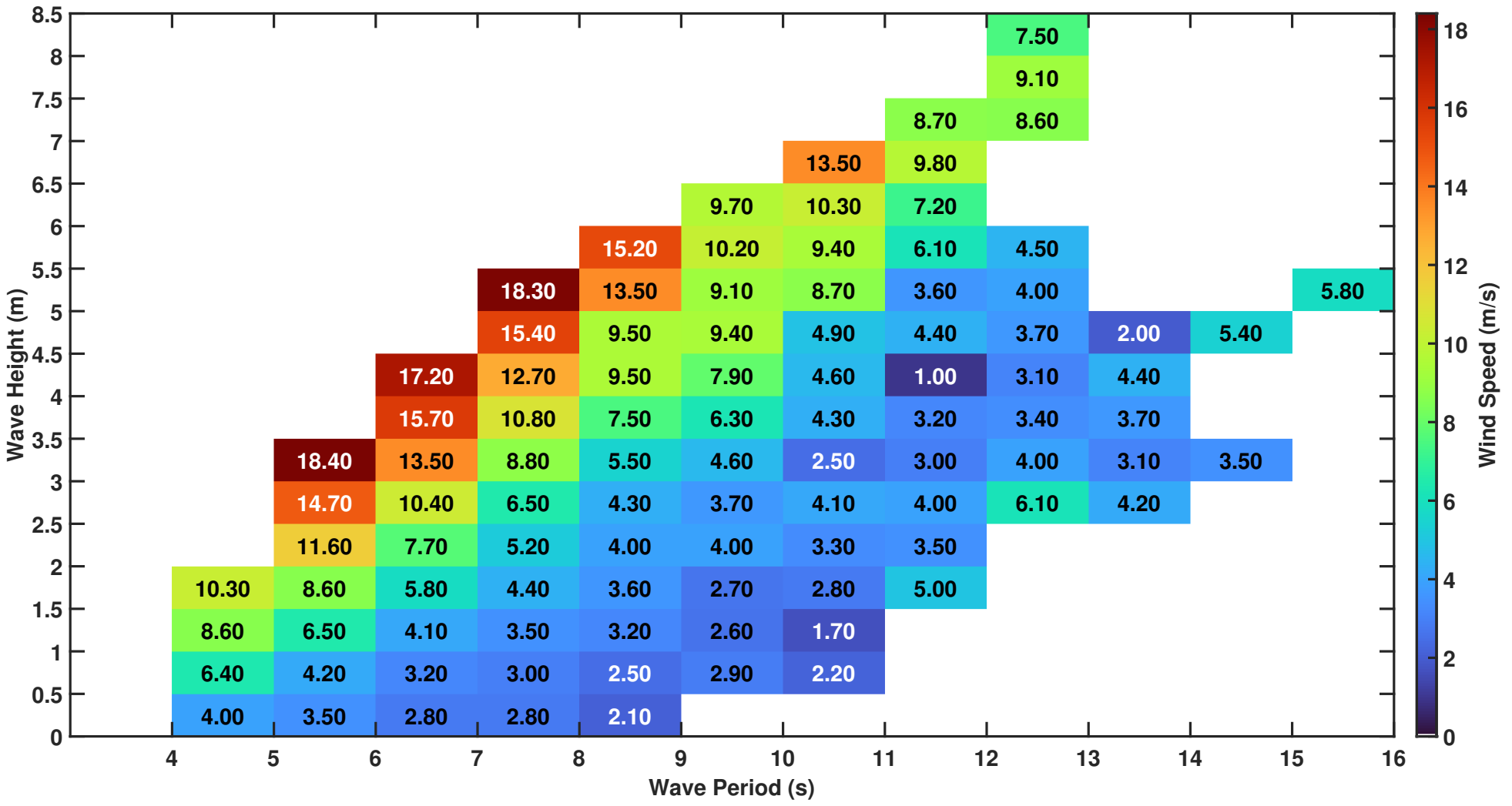

**Figure 3:** Averaged wind speed and wave conditions distribution at the Humboldt Site.

### 2.3. Combined wind-wave resource

To model a hybrid system excited by both waves and wind, a combined dataset of wave and wind conditions at the same location is required. Following a similar approach to that presented by Lou et al. (2024) Lou and Robertson (2024), a corresponding wind speed was assigned to each wave condition at the Humboldt Site. In Figure 3, each cell represents a specific wave condition, with the number within each cell indicating the averaged wind speed associated with that condition. The average wind speed was calculated to be 6.5 m/s, which matches with the wind speed range from the wind rose and the average wind speed reported by NDBC for the same location.

To investigate the impact of the location resources on the power generation by WTs, a second location with higher average wind speeds was selected for both standalone WTs. The second selected location is the Kodiak site at the western Gulf of Alaska (56.296° N 148.027° W), where the average wind speed is approximately 10 m/s. The wave height, energy period, and wind speed calculations are based on measurements from buoy station 46001 National Data Buoy Center (2024a). Figure 4 shows the averaged wind speed associated with each wave condition occurring at the Kodiak site.

## 3. Case studies

To examine the dynamic interaction between WECs and WTs in a hybrid wind-wave system, numerical simulations are performed under varying wave and wind conditions across multiple different configurations. The RM3 WEC Neary et al. (2014) (Figure 5) is selected for its simplicity and versatility, featuring a float that moves vertically along a spar with a reaction plate at its base. This study incorporates the NREL 5MW Jonkman et al. (2009) and IEA 15MW Gaertner et al. (2020) WTs to represent two distinct scales and capacities, broadening the scope of the research. The properties for both WTs are shown in Table 1.

To evaluate integration across different water depths and environmental conditions, two mooring foundations—spar Jonkman (2010) (Figure 6) and semi-submersible Robertson et al. (2014) (Figure 7)—are employed. In this hybrid wind-wave system, the WT tower also serves as the spar in the RM3 design, as shown in Figure 8, and the significance of adding a reaction plate at the base is explored.

The dimensions of the float in RM3 were adjusted to integrate with the 5MW and 15MW FWTs while aligning the float's natural frequency with the predominant wave frequency at the target site (6.5s) in order to maximize the power generated by the WECs. Table 2 shows the dimensions of the modified floats used in our study; float 1 is integrated with the 5MW WT, and float 2 is integrated with the 15MW WT. The dimensions of the spar and the semi-submersible platforms remained unchanged from the OC3 and OC4 platforms Jonkman (2010); Robertson et al. (2014) when assembled with the 5MW WT. However, the platforms required modification to integrate with the 15MW WT. Due to the increased weight of the 15MW WT, additional buoyancy was necessary to ensure system stability. The dimensions are shown in Table 3 and Table 4 for the spar platform and the semi-submersible platform, respectively. In our case studies, we also consider the incorporation of a reaction plate to the base at the base of the spar platform. Its radius is set to 25 m when integrated with the 5MW WT and 30 m for the 15 MW WT, with a thickness of 0.5 m for both. The radii of the reaction plate are larger than the one on the original RM3 design due to the presence of WTs on top of the spar.

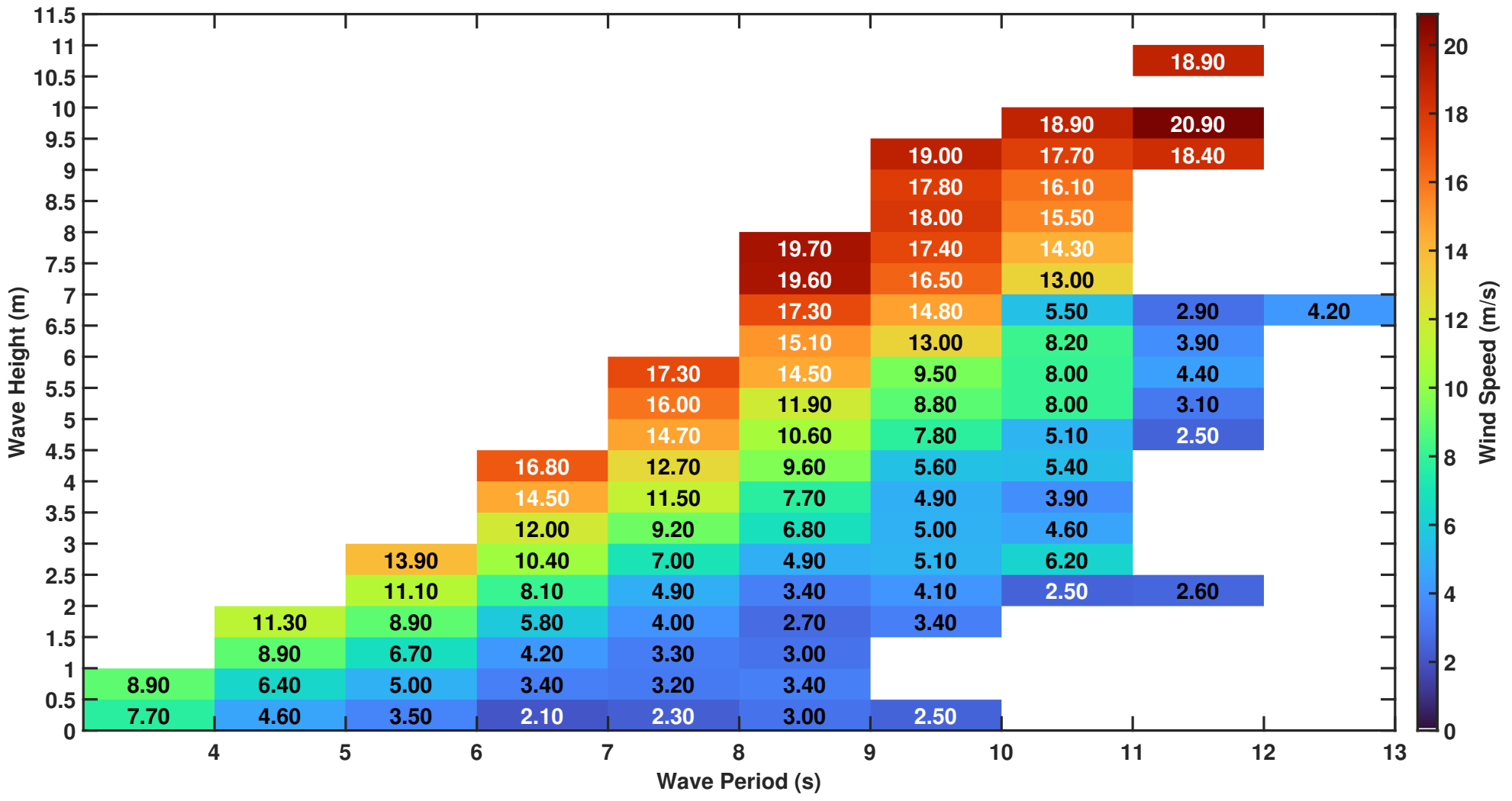

**Figure 4:** Averaged wind speed and wave conditions distribution at the Kodiak Site.

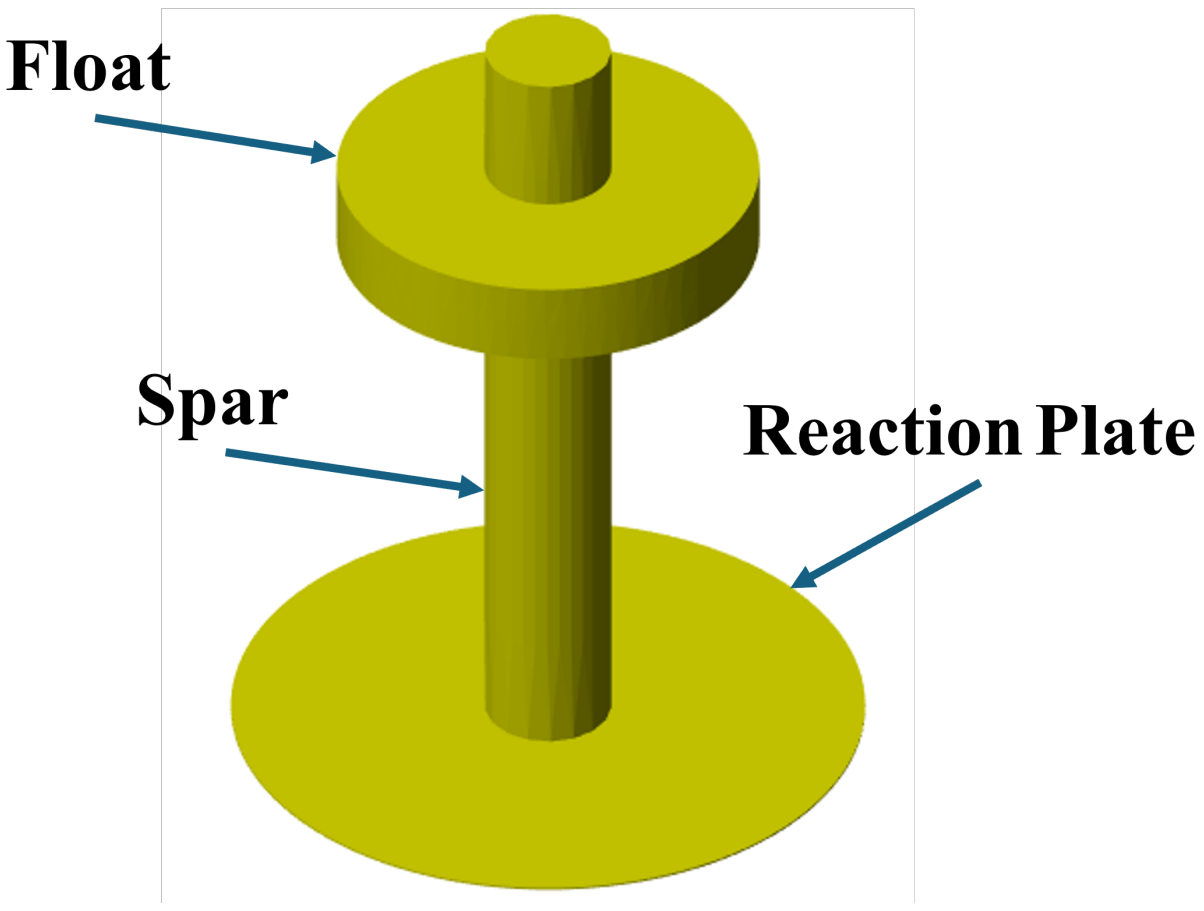

**Figure 5:** RM3 WEC

Six configurations combining various WTs and mooring foundations are analyzed to assess feasibility and performance. The details of these configurations are provided in Table 5. The mooring systems used to secure all different configurations are based on the mooring system in Jonkman (2010); Robertson et al. (2014). The same mooring system is used for both the 5MW and 15MW WTs with the spar platform in standalone and hybrid systems. However, a different system is employed for the 5MW and 15MW WTs with the semi-submersible platform. Additionally, the standalone WEC system utilizes a distinct mooring system with the RM3 configuration. Table 6 shows the mooring system properties for the three mooring foundations.

**Table 1**
Properties for the NREL 5 MW and IEA 15 MW WTs

| Parameter | NREL 5 MW Turbine | IEA 15 MW Turbine |
|---|---|---|
| Rating (MW) | 5 | 15 |
| Rotor Orientation | Upwind | Upwind |
| Number of Blades | 3 | 3 |
| Rotor Diameter (m) | 126 | 240 |
| Hub Diameter (m) | 3 | 7.94 |
| Hub Height (m) | 90 m | 150 m |
| Cut-In Wind Speed (m/s) | 3 | 3 |
| Rated Wind Speed (m/s) | 11.4 | 10.6 |
| Cut-Out Wind Speed (m/s) | 25 | 25 |
| Cut-In Rotor Speed (rpm) | 6.9 | 5 |
| Rated Rotor Speed (rpm) | 12.1 | 7.56 |
| Rated Tip Speed (m/s) | 80 | 95 |

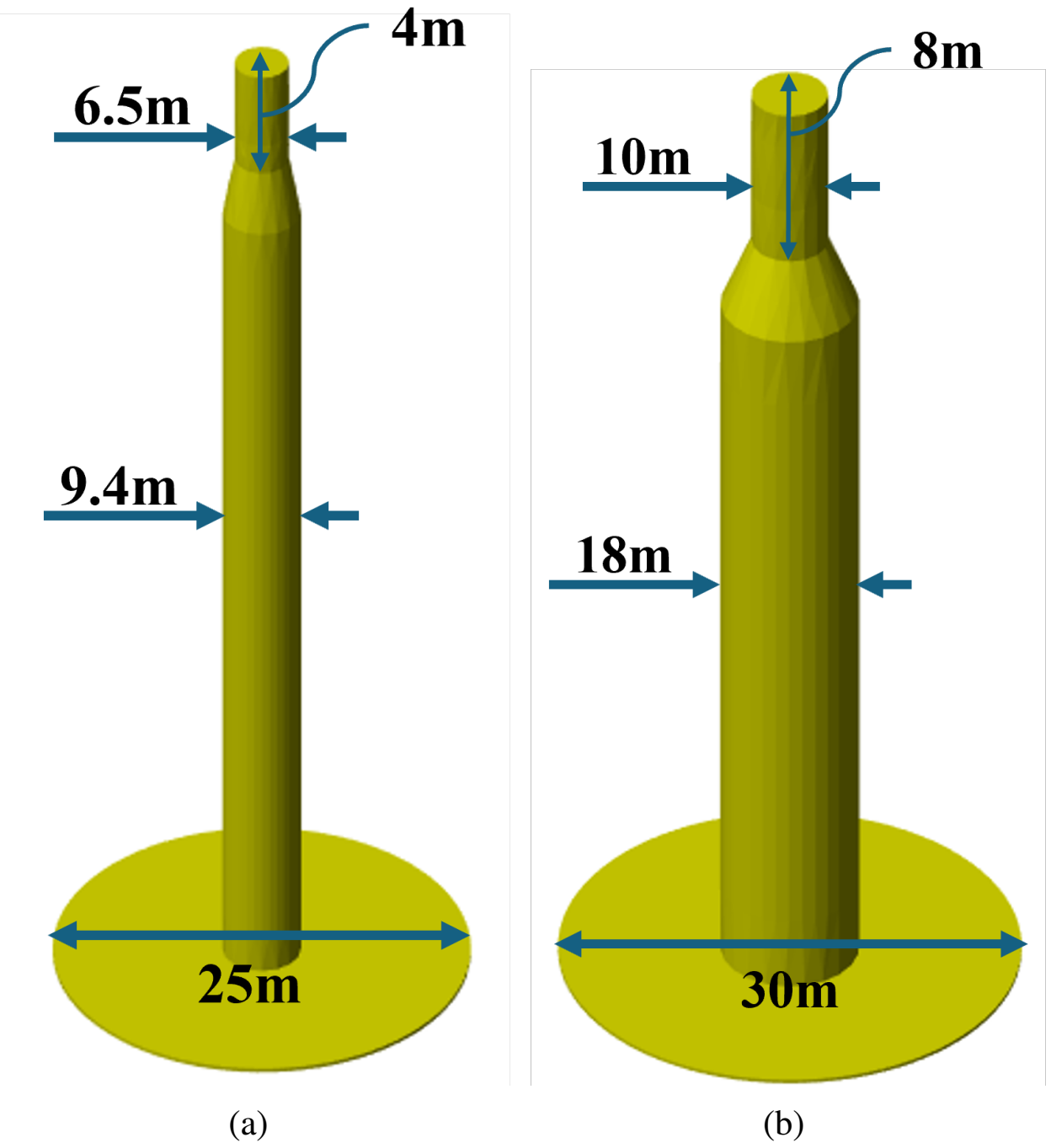

**Figure 6:** Spar platform with a reaction plate at its base for (a) 5MW and (b) 15MW WTs

## 4. System Synergy

While synergy has been conceptually discussed in the context of hybrid renewables, this study is the first to operationalize and apply it to the evaluation of real hybrid wind–wave designs. In this context, system synergy refers to the combined benefits that emerge when FWTs and WECs are integrated into a single platform. A synergistic system is achieved when the hybrid system outperforms its standalone counterparts in overall metrics and operational efficiency. Depending on the interaction, this relationship may manifest as mutualism, where both technologies benefit from hybridization; commensalism, where one benefits without affecting the other; or parasitism, where one benefits at the expense of the other Garcia-Sanz (2024). This balance of benefits, often referred to as system synergy, is crucial to justifying hybridization as a viable investment in the offshore renewable sector. Identifying and quantifying synergy is therefore critical to determining whether hybrid wind-wave systems provide genuine techno-economic advantages over standalone approaches. This study contributes by quantifying these interaction types in realistic case studies, advancing synergy from a theoretical construct to an applied evaluation framework.

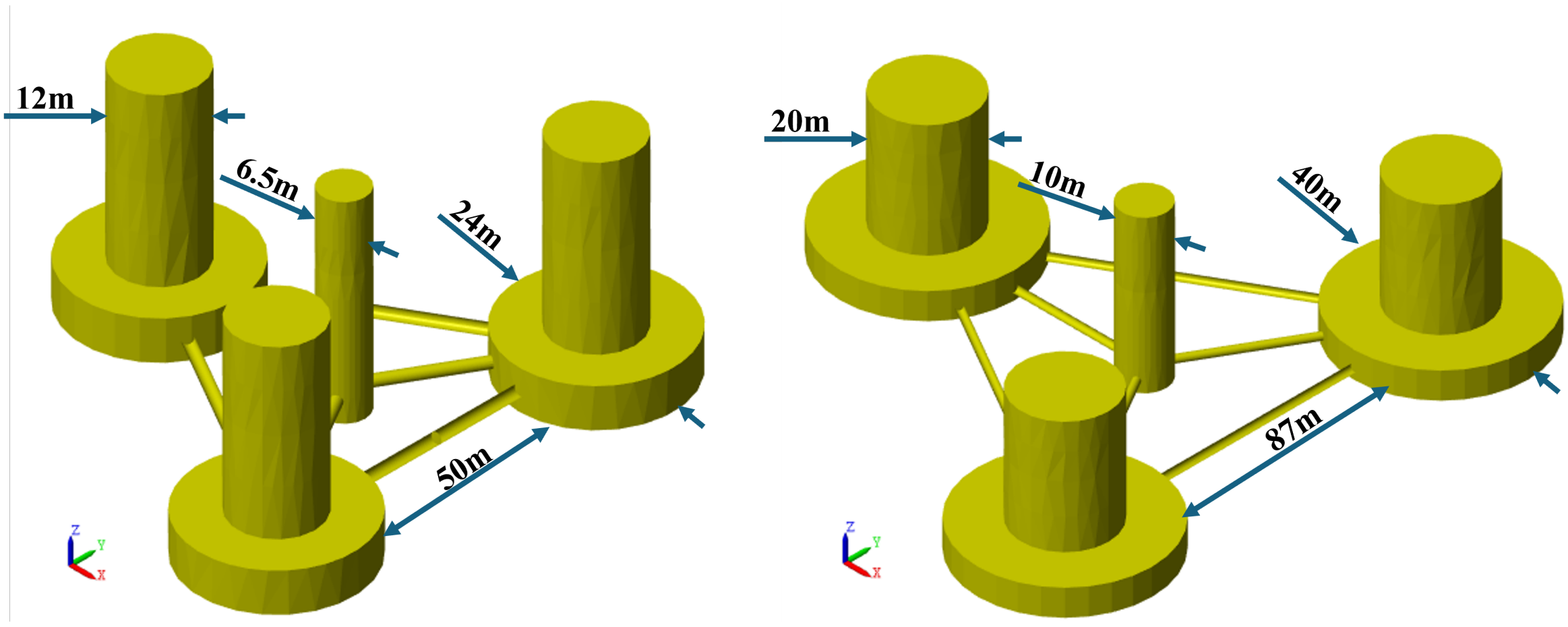

(a) (b)

**Figure 7:** Semi-submersible platform for (a) 5MW and (b) 15MW WTs

**Table 2**
Dimensions of the modified floats

| Parameter | Float 1 for 5 MW WT | Float 2 for 15 MW WT |
|---|---|---|
| Inner Diameter (m) | 6.5 | 10 |
| Outer Diameter (m) | 22 | 30 |
| Height (m) | 5 | 8 |
| Draft (m) | 3 | 5 |
| CM Location Below SWL (m) | 0.72 | 1.35 |
| PTO Damping Coefficient (MN/(m/s)) | 1.2 | 1.2 |
| Natural Period (s) | 6.23 | 6.12 |

**Table 3**
Dimensions of the spar platform for both 5MW and 15MW WTs

| Parameter | Spar platform for 5MW WT | Spar platform for 15MW WT |
|---|---|---|
| Platform Diameter Above Taper (m) | 6.5 | 10 |
| Platform Diameter Below Taper (m) | 9.4 | 18 |
| Depth to Top of Taper Below SWL (m) | 4 | 8 |
| Depth to Bottom of Taper Below SWL (m) | 12 | 18 |
| Elevation to Platform Top (Tower Base) Above SWL (m) | 10 | 15 |
| Depth to Platform Base Below SWL (Total Draft) (m) | 120 | 120 |
| CM Location Below SWL Along Platform Centerline (m) | -90 | -100 |

To evaluate the effectiveness and synergy of hybrid FWT-WEC systems, performance metrics are calculated for both hybrid and standalone configurations and then compared to determine the type of symbiosis present. As shown in Figure 9, mutualism occurs when a metric decreases for both devices after integration; commensalism when it improves for one device while remaining unchanged for the other; and parasitism when it improves for one device but worsens for the other.

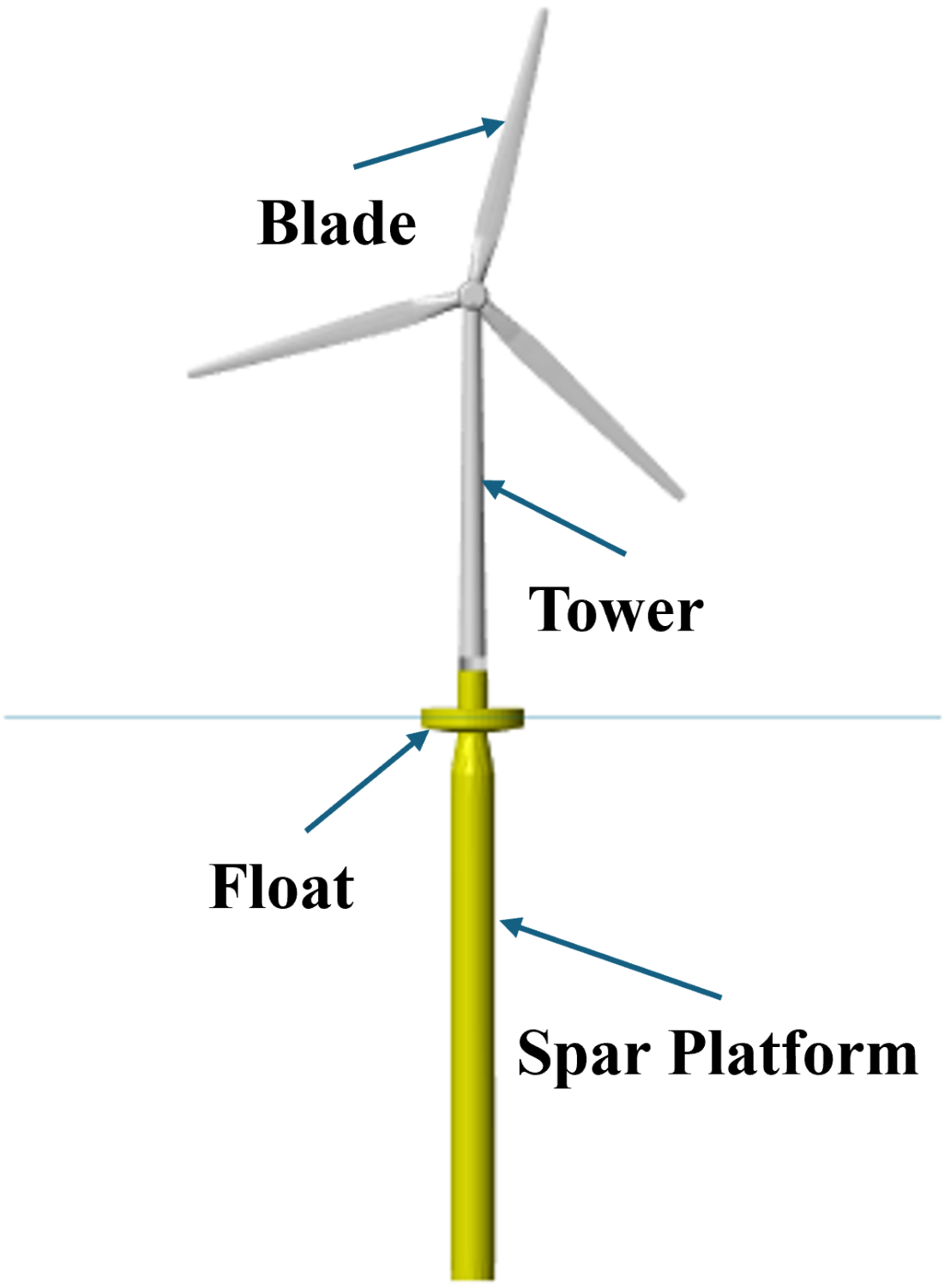

**Figure 8:** Hybrid wind-wave system

**Table 4**
Dimensions of the semi-submersible platform for both 5MW and 15MW WTs

| Parameter | Semi-submersible platform for 5MW WT | Semi-submersible platform for 15MW WT |
|---|---|---|
| Diameter of main column (m) | 6.5 | 10 |
| Diameter of offset (upper) columns (m) | 12 | 20 |
| Diameter of base columns (m) | 24 | 40 |
| Elevation of main column (tower base) above SWL (m) | 10 | 15 |
| Elevation of offset columns above SWL (m) | 12 | 12 |
| Depth to top of base columns below SWL (m) | 14 | 14 |
| Depth of platform base below SWL (total draft) (m) | 20 | 20 |
| Spacing between offset columns (m) | 50 | 87 |
| CM location below SWL (m) | -13.46 | -14 |

## 5. Methodology

The methodology adopted to assess the performance of hybrid wind–wave energy systems is outlined in this section. It presents the numerical simulation framework, the performance metrics used for evaluation, and the criteria applied to analyze system synergy. Together, these elements provide a structured approach for investigating both the dynamic behavior and the techno-economic feasibility of the proposed configurations.

### 5.1. Numerical Simulations

The coupled modeling framework, which integrates the inputs and outputs for each module, is shown in Figure 10 and is used to simulate the six hybrid configurations detailed in Table 5. The modeling framework couples the open-source software WEC-Sim (Wave Energy Converter SIMulator) Ruehl et al. (2024) with MOST (MATLAB for Offshore floating wind turbine Simulation Tool) Cottura et al. (2021) and MoorDyn Hall (2020). WEC-Sim is a time-domain modeling tool developed by Sandia National Laboratories and NREL to simulate the performance of WECs under operational wave conditions. Built in MATLAB/Simulink, WEC-Sim utilizes the multi-body dynamics solver SimMechanics to represent WECs. The modeling framework connects rigid bodies using joint or constraint blocks from the WEC-Sim library. This open-source, publicly available tool enables comprehensive simulations of WEC systems. All the hydrodynamic load calculations related to the WEC and mooring are handled in WEC-Sim, while the aerodynamic loads and the WT electrical generator calculations are performed in MOST. MOST is a wind-to-power modeling framework that integrates all key components of a FWT, including the floating platform and its hydrodynamics, the WT and its aerodynamics, as well as the electric generator and its control system. It serves as a comprehensive tool for the planning and design phases of FWT projects. To ensure adaptability and ease of use, MOST is developed in MATLAB-Simulink, a platform known for its intuitive structure and user-friendly interface, making it easily programmable and customizable. Lastly, the mooring loads and cable tensions are calculated using MoorDyn, a lumped-mass model developed by NREL to simulate the dynamics of mooring systems connected to floating offshore structures. MoorDyn can be compiled as a dynamically linked library or wrapped for use in MATLAB. It features simplified functions to facilitate easy coupling with WEC-Sim.

**Table 5**
Description of configurations integrating various WTs and mooring foundations

| Configuration | Wind Turbine | Mooring Foundation | Legend |
|---|---|---|---|
| 1 | | Spar | Hybrid 5MW Spar |
| 2 | 5MW | Spar with Reaction plate | Hybrid 5MW Spar w/ RP |
| 3 | | Semi-submersible | Hybrid 5MW SS |
| 4 | | Spar | Hybrid 15MW Spar |
| 5 | 15MW | Spar with Reaction plate | Hybrid 15MW Spar w/ RP |
| 6 | | Semi-submersible | Hybrid 15MW SS |

**Table 6**
Mooring system properties for spar and semi-submersible platforms, and RM3

| | Spar Platform | Semi-submersible Platform | RM3 |
|---|---|---|---|
| Number of Mooring Lines | 3 | 3 | 3 |
| Angle Between Adjacent Lines | 120° | 120° | 120° |
| Depth to Anchors Below SWL (m) | 320 | 200 | 70 |
| Depth to Fairleads Below SWL (m) | 70 | 14 | 11.5 |
| Unstretched Mooring Line Length (m) | 902.2 | 835.5 | 280 |
| Mooring Line Diameter (m) | 0.09 | 0.0766 | 0.144 |
| Equivalent Mooring Line Mass Density (kg/m) | 77.71 | 113.35 | 126 |
| Equivalent Mooring Line Mass in Water (kg/m) | 71.16 | 108.63 | 115 |
| Equivalent Mooring Line Extensional Stiffness (MN) | 384.2 | 753.6 | 583.4 |

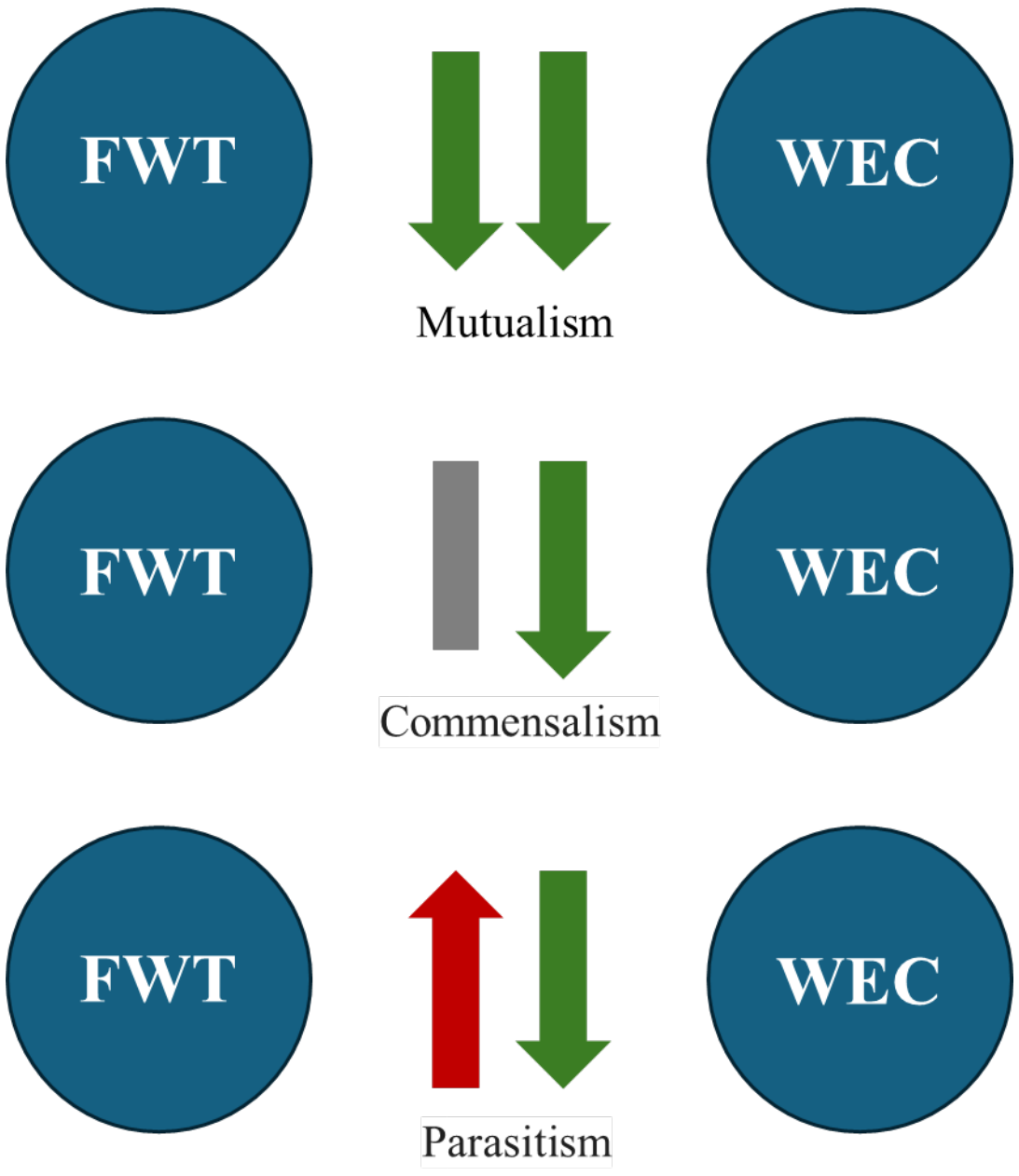

**Figure 9:** Comparison of the symbiotic effect between FWT and WEC

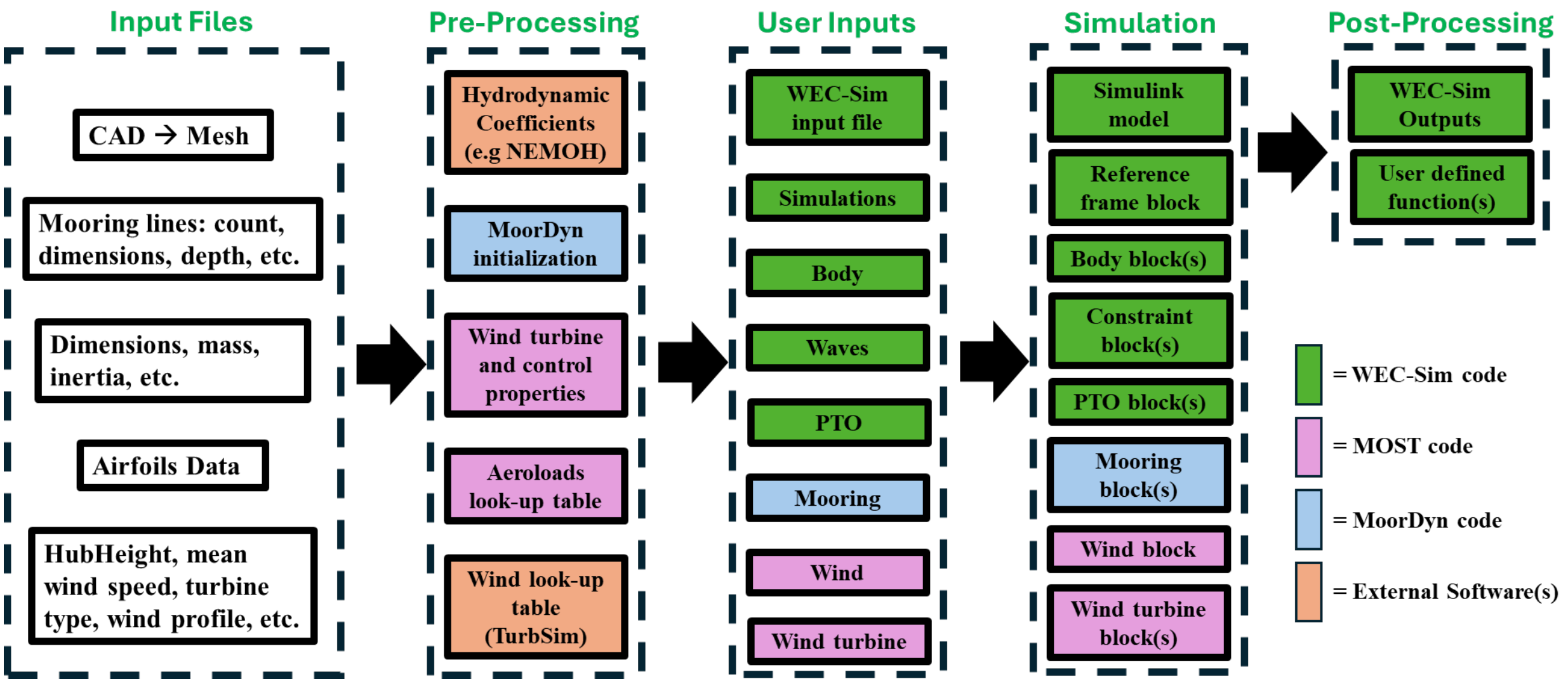

**Figure 10:** Coupling framework for hybrid system modeling

The coupling between the three simulation tools ensures that all dynamic interactions between the systems are captured, providing insights into how each system influences the other and affects overall performance. To the authors' knowledge, no published studies have utilized MOST+WEC-Sim to model and generate results for hybrid FWT-WEC systems. This study represents the first application of the coupling framework illustrated in Figure 10 for analyzing such integrated systems.

One limitation of this framework is that the numerical simulations rely on linear wave theory, which is based on potential flow theory. These assumptions do not fully capture nonlinear hydrodynamic effects, such as viscous drag or wave-breaking phenomena, which are instead incorporated as external terms in the calculations. More details about these values and how they are calculated can be found in Jonkman (2010); Robertson et al. (2014).

## 5.2. Performance Metrics

Evaluating the performance of hybrid energy systems requires the selection of some critical metrics such as AEP, LCOE, and $P_{CV}$. AEP reflects the total energy output over a year and directly depends on the availability of wind and wave resources at the selected site. In a hybrid system, energy harvesting from WT and WECs is expected to increase the overall AEP, as the two resources can complement each other, particularly since wave energy tends to be more consistent throughout the year compared to wind, which exhibits greater seasonal variability Ibáñez-Rioja et al. (2023); Tang et al. (2022). Evaluating $P_{CV}$ provides insight into the reliability and stability of the energy supply, with a lower variation indicating a more continuous generation profile. This improved consistency reduces reliance on energy storage systems, which in turn can lower capital and operational costs. Ultimately, these benefits are reflected in the LCOE, a key metric that quantifies the cost of electricity offered to the market, thereby affecting the competitiveness and affordability of hybrid renewable energy systems.

### 5.2.1. Annual Energy Production (AEP)

The AEP is determined by simulating the power output for various environmental conditions represented in the JPD of wave heights and periods, as shown in Figure 1 when simulating the standalone WEC system, and the wind rose in Figure 2 when simulating the standalone FWT. For the hybrid FWT-WEC, each cell in the JPD is associated with an averaged wind speed as computed in Figure 3 based on the study in Lou and Robertson (2024). A separate simulation is conducted for each wave condition with the associated wind speed, ensuring a comprehensive assessment of power generation across different environmental scenarios. Following the simulations, the power output for each condition is weighted by its occurrence probability, as defined in the JPD. This approach ensures that the AEP accounts for the relative frequency of different wave and wind conditions over a typical year.

The power generated by the WEC is derived from the PTO force computed in WEC-Sim, while the WT power output is obtained using the MOST model. The final set of wave conditions and their corresponding wind speeds used in the simulations are provided in Figure 3. This methodology provides an accurate estimation of the expected energy yield from the hybrid system, aiding in the assessment of its economic and operational viability.

To transfer mechanical power harvested by the WEC and the WT to electrical power, an efficiency of 80% is assumed. This lumped parameter serves as an approximation, since the achievable efficiency ultimately depends on control co-design and implementation strategy, as highlighted in prior studies on WEC PTO systems and control co-design Ströfer, Gaebele, Coe and Bacelli (2023); Ringwood (2025). The electrical power is capped based on the generator used; it was rated at 286 kW for the WEC and 5 MW and 15 MW for the WTs. Lastly, some losses due to the device availability and transmission were assumed to be 95% and 98%, respectively.

The capacity factor (CF) is calculated using Equation 1, it is defined as the ratio of actual AEP to the maximum possible AEP, the latter being achieved if the system operated at its rated power continuously throughout the year. This metric provides a normalized measure of the system's utilization and performance under realistic environmental conditions and can highlight the impact of location, as sites with higher wind and wave resources result in higher energy output and higher CF. CF is applied to both the WEC and FWT in all standalone configurations.

$$\text{CF} = \frac{\text{AEP}}{P_{rated} \times 365 \times 24} \tag{1}$$

### 5.2.2. Power Coefficient of Variation ($P_{CV}$)

The power generation is linked to and impacted directly by the wind and wave availability at the targeted site. Therefore, the variability in wind and wave resources can be assessed by quantifying the variability of the generated power using the power coefficient of variation ($P_{CV}$). It is defined as the ratio of the standard deviation to the mean power output (Equation 2), is a useful metric for assessing the stability of power generation.

$$P_{CV} = \frac{\sigma}{\mu} \tag{2}$$

A lower $P_{CV}$ indicates a more consistent and predictable power supply. In a hybrid system integrating a FWT with a WEC, the $P_{CV}$ is expected to be lower than that of a standalone FWT. This is because wave energy is available year-round with relatively steady fluctuations, whereas wind speeds exhibit significant seasonal variations, leading to periods of low or curtailed power output. By harnessing wave energy alongside wind, the hybrid system is expected to achieve a smoother and more stable power generation, reducing instances of curtailment. Consequently, this improved stability is expected to minimize the need for energy storage, leading to a lower LCOE. By demonstrating improved power stability, the $P_{CV}$ serves as a quantitative measure of the hybrid system's advantage over standalone FWTs, reinforcing its economic and operational benefits.

### 5.2.3. Levelized Cost of Energy (LCOE)

LCOE is calculated using Equation 3 where (CapEx) is capital expenditure ($), (OpEx) is operational expenditure ($/year), FCR is the fixed change rate, and AEP is the annual energy production (MWh/year). Figure 11 shows the components included in the calculations. The deployment of FWTs and WECs typically occurs in the form of large-scale farms rather than individual devices, which enables cost reductions through economies of scale and mass production. To demonstrate the feasibility of FWT-WEC integration, we consider a 100-unit farm as a case study. Three main system configurations are analyzed: (1) a farm composed of 100 standalone RM3 WECs, (2) a farm composed of 100 standalone 5 MW FWTs, and (3) a farm composed of 100 hybrid FWT-WEC systems. For each configuration that has a FWT, three different mooring foundation designs, as summarized in Table 5, are evaluated. This results in a total of seven different farm scenarios being assessed. In addition, the same seven configurations are re-evaluated using a 15 MW FWT instead of the 5 MW model to examine the impact of turbine scaling on system performance. The hybrid FWT-WEC system will utilize a shared mooring foundation, eliminating the need for a separate mooring system for the WEC and removing its associated costs from the LCOE. The system will also share electrical cables and transmission infrastructure, further reducing both capital expenditures (CapEx) and operational expenditures (OpEx). By accounting for all shared components and recalculating these expenditures, we will determine the new LCOE for each device using Equation 3. Neary et al. (2014). The costs associated with each configuration are drawn from publicly available data Neary et al. (2014); Beiter, Musial, Smith, Kilcher, Damiani, Maness, Sirnivas, Stehly, Gevorgian, Mooney et al. (2016); Myhr, Bjerkseter, Ågotnes and Nygaard (2014); Martinez and Iglesias (2022) and adjusted to 2025 USD, enabling a robust comparison of their respective LCOEs. This study will highlight the economic trade-offs and synergies in hybrid wind-wave farms, offering insights into their potential for large-scale deployment.

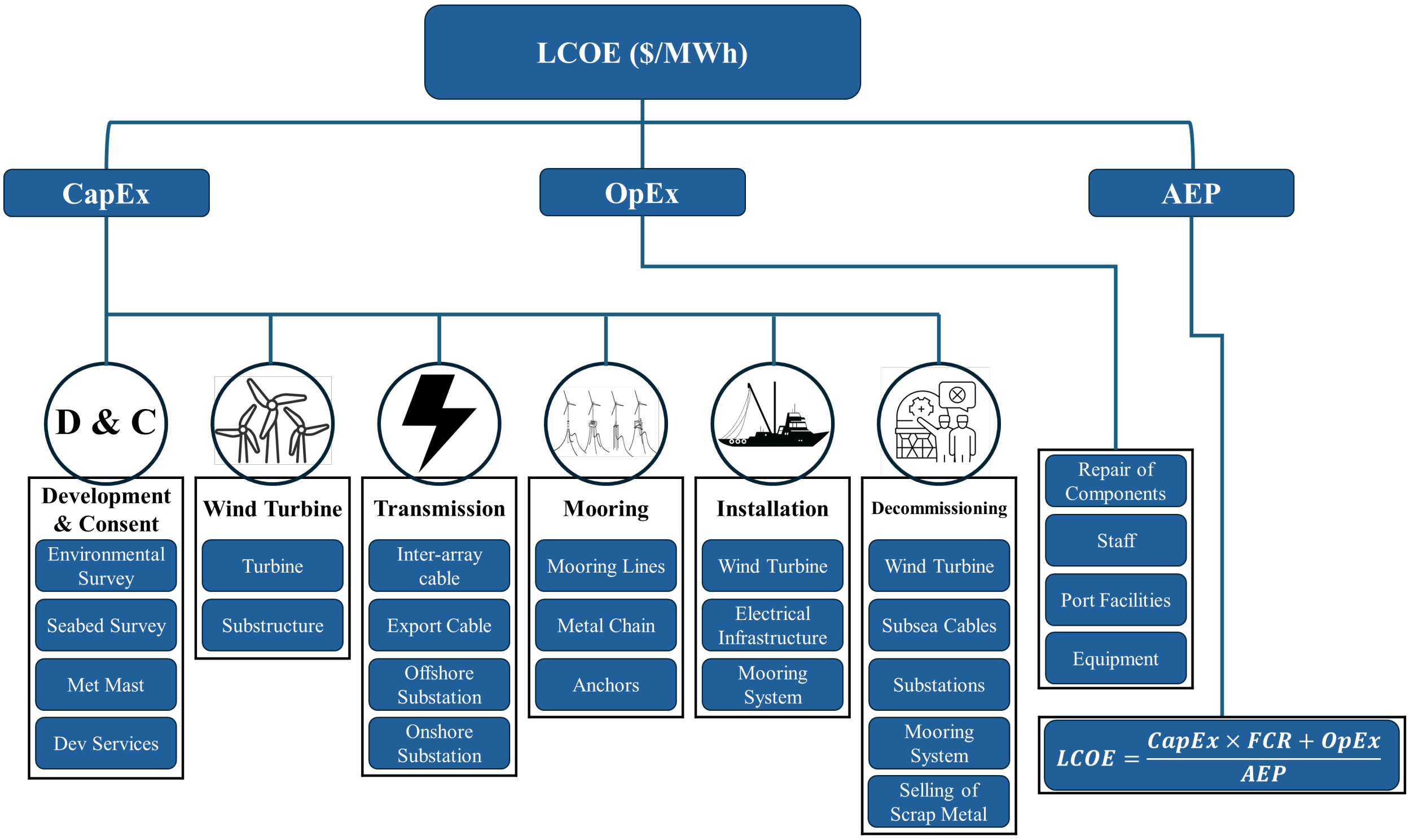

**Figure 11:** Breakdown of LCOE and components included in its calculation Martinez and Iglesias (2022)

$$\text{LCOE} = \frac{\text{CapEx} \times \text{FCR} + \text{OpEx}}{\text{AEP}} \tag{3}$$

## 5.3. Synergy Metrics

To identify the synergy of hybrid FWT-WEC systems, we propose using the LCOE and $P_{CV}$ as key performance metrics. These two metrics were chosen because hybridization is expected to both reduce the overall cost of energy and provide a more stable power supply by leveraging the complementary nature of wind and wave resources: when wind speeds are low, wave energy is often still available, and vice versa. As a result, hybrid systems have the potential to both reduce the cost of energy by sharing infrastructure and improve power stability by smoothing fluctuations. LCOE provides a direct measure of the economic impact of integration, while $P_{CV}$ captures the variability of power output and demonstrates the extent to which hybridization enhances reliability. The standalone LCOEs of FWTs and WECs are compared with the LCOE of the integrated system to determine the type of symbiosis present. As illustrated in Figure 9, mutualism occurs if the LCOE decreases for both devices after integration, commensalism occurs if the LCOE remains unchanged for one device while decreasing for the other, and parasitism occurs if the LCOE increases for one device while decreasing for the other. A similar approach is applied to evaluate synergy using the $P_{CV}$. In this case, mutualism is achieved if the $P_{CV}$ of both devices decreases, commensalism occurs when one decreases while the other remains unchanged, and parasitism arises when one decreases while the other increases. This dual-metric framework enables a more comprehensive assessment of both the economic viability and potential benefits of integrating FWTs and WECs.

## 6. Results and Discussions

The following section presents the results of the dynamic response, power performance, and techno-economic evaluation of the proposed hybrid wind–wave energy systems. Comparisons are made between hybrid configurations and their standalone counterparts to assess improvements in dynamic response, energy production, power stability, and cost-effectiveness. These results provide insight into the benefits and trade-offs of integrating FWTs with WECs across different platform designs and scales.

### 6.1. Dynamics

The impact of different components of the systems and the dynamic response of the hybrid configurations under different wave and wind conditions are discussed, with comparisons to standalone systems to assess the influence of integration on motions and stability.

#### *6.1.1. Impact of reaction plate on AEP*

The reaction plate is a critical component of the RM3 device, contributing significantly to its CapEx while playing a key role in stabilizing heave and pitch motions to maximize energy capture. However, in a hybrid FWT-WEC system, the FWT may provide additional stability to the RM3 float, potentially reducing or eliminating the need for the reaction plate. Given that the reaction plate accounts for approximately 12% of the total CapEx of RM3, its removal could lead to a substantial reduction in the LCOE, provided the AEP remains unaffected. This section investigates the impact of the reaction plate on the AEP of the RM3 device, exploring the potential for cost savings and performance optimization in hybrid FWT-WEC systems. The impact of the reaction plate on AEP was assessed through numerical simulations for two scenarios: with and without the reaction plate, under identical wave conditions at the Humboldt Site. Simulations were conducted at three water depths (40 m, 90 m, and 120 m), representing deep-water scenarios across various locations.

The results, shown in Figure 12, indicate a significant reduction in AEP when the reaction plate is not present. Specifically, at 40 m, the AEP is decreased 60% at a depth of 40 m, 45% at 90 m, and 55% at 120 m. In contrast, extending the spar length while retaining the reaction plate increased AEP by 50% at 90 m compared to 40 m and 40% at 120 m compared to 40 m. Although elongating the spar raises the material costs, it improves the LCOE by reducing it by approximately 20%. Conversely, the cost savings from the absent reaction plate failed to offset the corresponding decrease in energy production, resulting in a higher LCOE. These findings underscore the importance of retaining the reaction plate to ensure system stability and to maximize energy production, especially at greater water depths. Building on this observation, the spar mooring platform with a reaction plate at its bottom was further tested with the hybrid wind-wave system to investigate its influence on the performance in deeper water conditions.

#### *6.1.2. Dynamic interactions of hybrid system*

Three distinct wave conditions for the Humboldt site were chosen to evaluate the dynamic interactions of the hybrid and standalone systems. The first condition, with a wave period of 6.5 s and a height of 1.75 m, represents the most occurring wave. The second, with a 10.5 s period and 3.75 m height, corresponds to the highest annualized wave. The third, characterized by a 14.5 s period and 6.75 m height, represents an extreme high-amplitude scenario. Table 7 provides a summary of these conditions. The wind speeds used in the hybrid system simulations were the rated wind speeds for each WT: 11.4 m/s for the 5MW turbine and 10.6 m/s for the 15MW turbine. All six hybrid system configurations were tested and compared with standalone WECs and WTs. These configurations, detailed in Table 5, were analyzed under the three wave conditions. The system response was evaluated across six degrees of freedom: surge (X), sway (Y), heave (Z), roll (RX), pitch (RY), and yaw (RZ). Particular attention was given to surge, heave, and pitch, as these significantly impact system performance, stability, and power generation. To ensure consistency with the standalone system studies, the water depth was set to 320 m for the configurations with the spar platform and 200 m for the configurations with the semi-submersible platform.

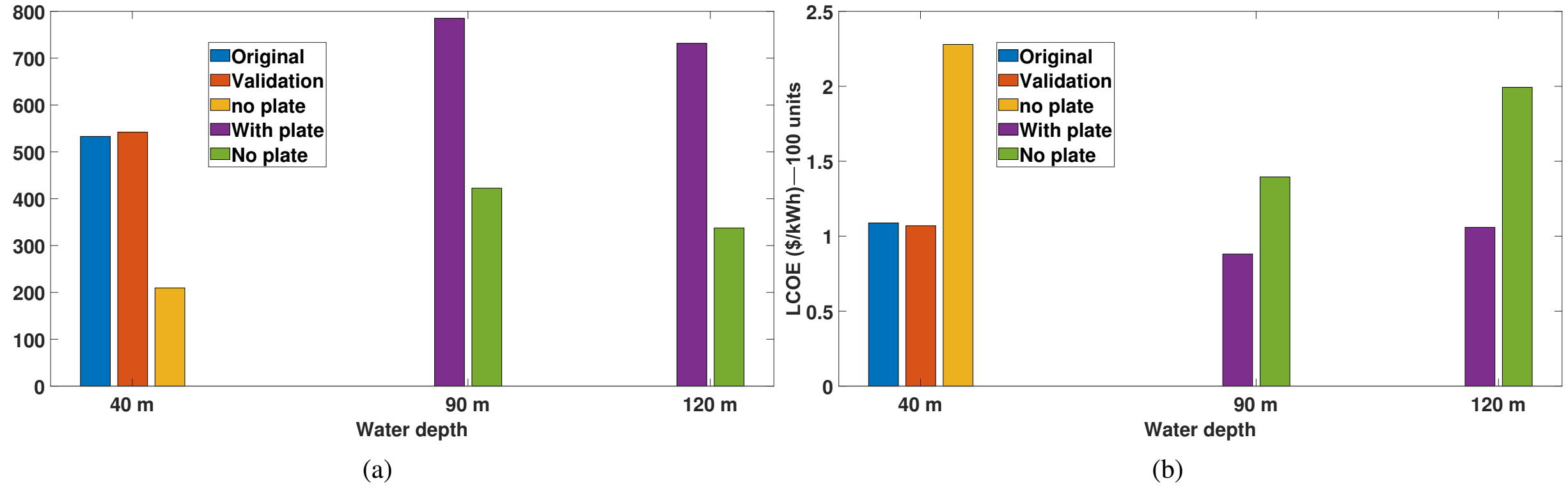

**Figure 12:** Impact of the reaction plate on (a) AEP and (b) LCOE at different water depths

**Table 7**
Selected wave conditions for testing the hybrid and standalone systems

| Runs | Wave Period (s) | Wave Height (m) |
|---|---|---|
| 1 | 6.5 | 1.75 |
| 2 | 10.5 | 3.75 |
| 3 | 14.5 | 6.75 |

**Table 8**
The RMS of the float's heave response in standalone RM3 for the six configurations under different wave conditions

| Wave condition type | Regular | | | Irregular | | |
|---|---|---|---|---|---|---|
| Run number | Run 1 | Run 2 | Run 3 | Run 1 | Run 2 | Run 3 |
| Float heave response (m) | 0.37 | 1.26 | 2.45 | 0.20 | 0.70 | 1.58 |

### *Float*

The analysis of the float's motion focuses on its heave response, as its movement is constrained along the tower of the WT. A comparative assessment was conducted between the float's heave response in different integrated WT configurations and its behavior in the RM3 standalone configuration at the Humboldt site. Table 8 presents the root mean square (RMS) of the float's heave response in the standalone RM3 system under regular and irregular wave conditions. The difference in the RMS of the time-domain dynamic response under the three different wave conditions in Table 7 is presented in Figure 13 for regular wave conditions and Figure 14 for irregular wave conditions, where negative values indicate a reduction in response relative to the standalone system. To illustrate these variations, the time series of a selected case of the float's response in the X, Z, and RY directions under excitation from a wave with a 6.5 s period and 1.75 m height are presented in Figure 15 for integration with the 5MW WT on a spar platform, with and without the reaction plate; Figure 16 for the 5MW turbine on a semi-submersible platform; Figure 17 for the 15MW turbine on a spar platform, with and without the reaction plate; and Figure 18 for the 15MW turbine on a semi-submersible platform.

The integration of the float with the WT on the spar platform, both with and without the reaction plate, results in relatively small changes in motion across most cases, with variations typically remaining below 5%. However, for specific configurations, particularly those involving the 15MW WT, the presence of the reaction plate influences the float's heave response up to 15% increase. The reaction plate enhances the hydrodynamic behavior of the system by modifying the damping and added mass characteristics, leading to an increase in the float's vertical motion, which directly contributes to greater power extraction by the WEC. On the other hand, when the float is integrated with the WTs on the semi-submersible platform and excited under regular wave conditions, the float's heave motion exhibits an increase ranging from 5% to 27%. For irregular wave conditions, the increase in heave motion varies between 9% and 20%. The most pronounced amplification in heave motion is observed in configurations involving the 15MW WT, in both regular and irregular wave conditions.

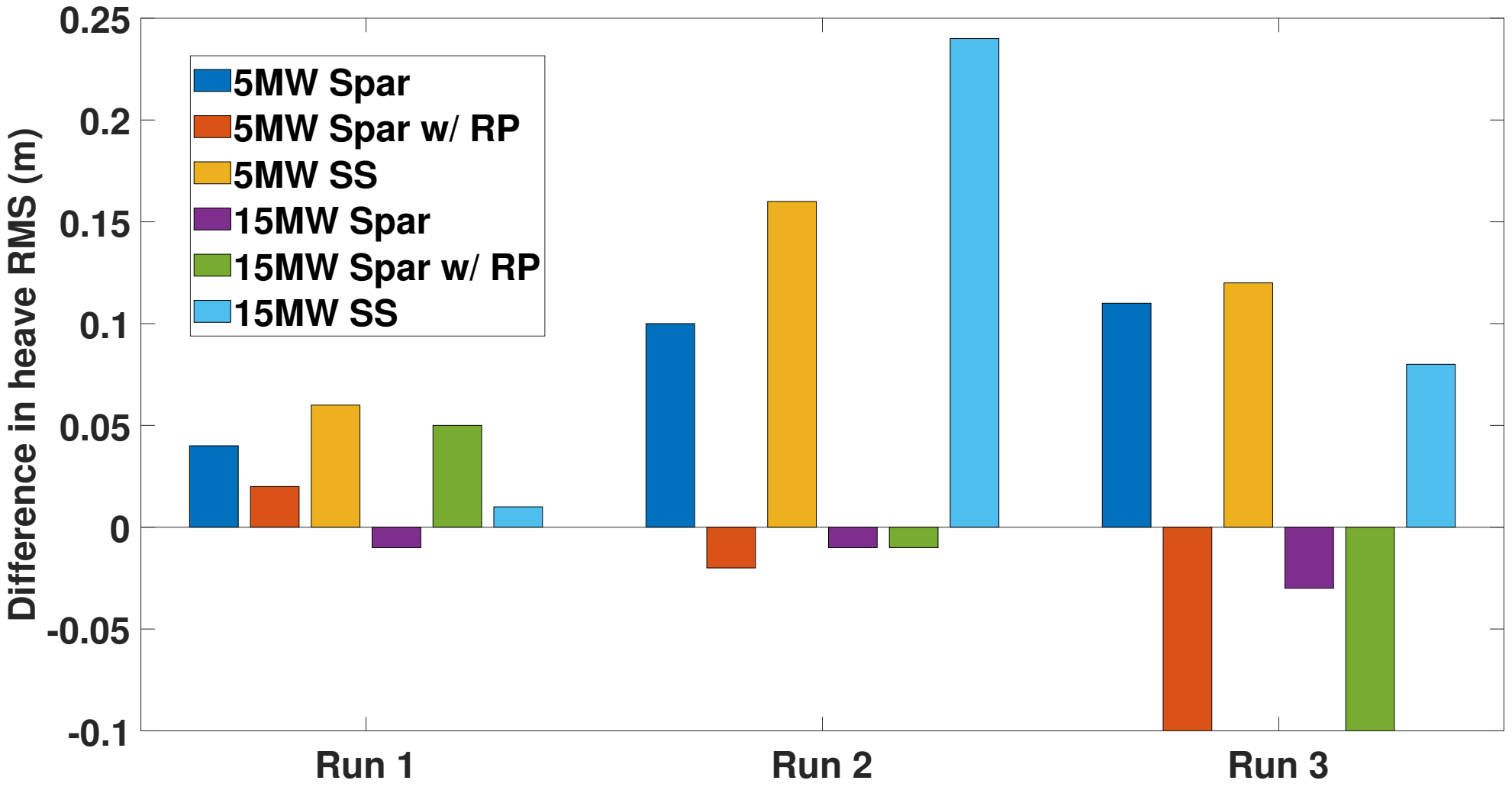

**Figure 13:** Difference in RMS of the float's heave response RMS between standalone and hybrid system for the six configurations under different regular wave conditions—Zero value means equivalent to the corresponding standalone system RMS

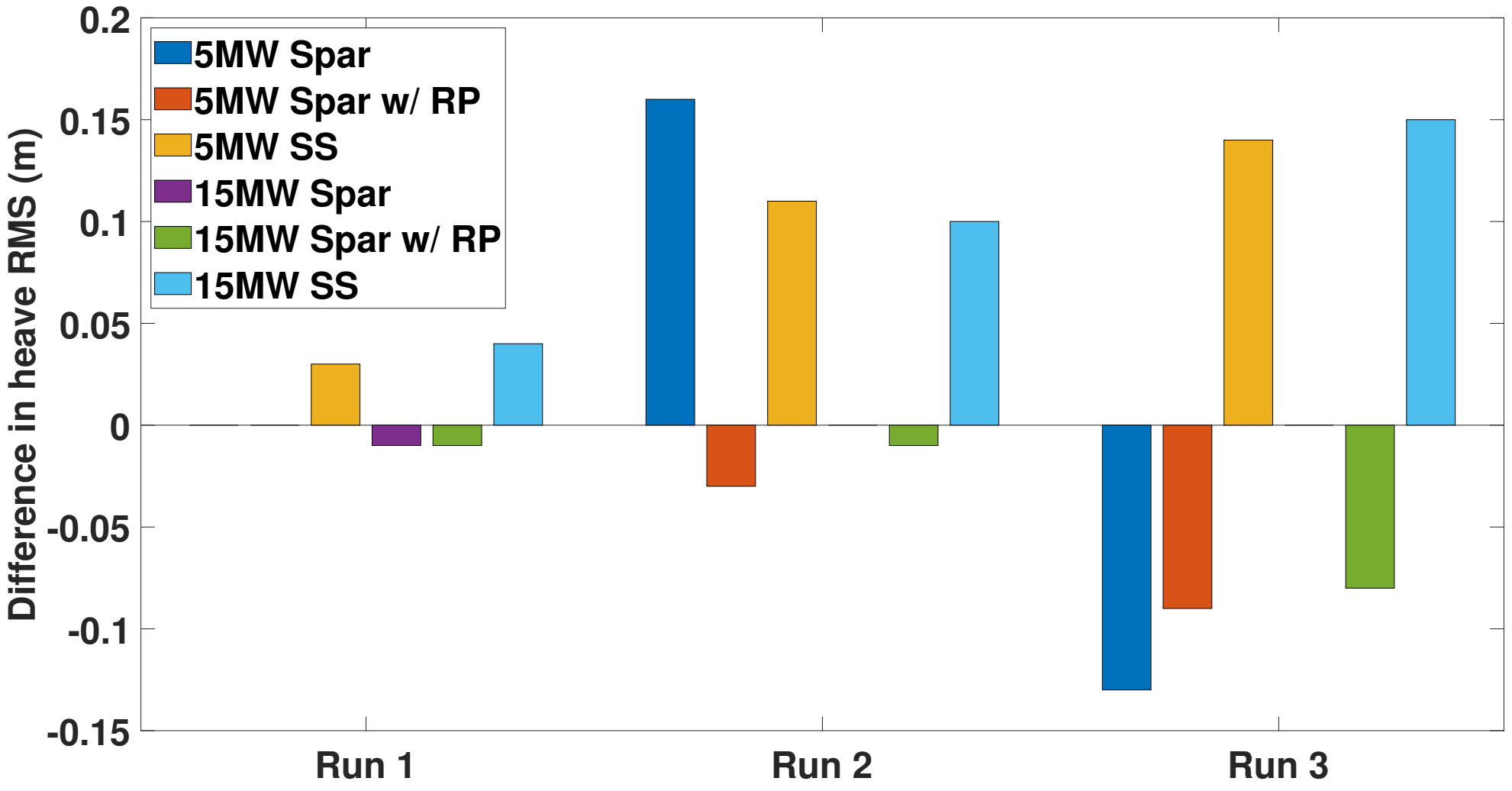

**Figure 14:** Difference in RMS of the float's heave response RMS between standalone and hybrid system for the six configurations under different irregular wave conditions—Zero value means equivalent to the corresponding standalone system RMS

The equilibrium position of the hybrid system shifts from that of the standalone WEC, as illustrated in Figures 15-18. This shift is primarily attributed to the aerodynamic forces exerted by the WT, which influence the system's response across multiple degrees of freedom. The addition of the reaction plate at the base of the spar platform mitigates this effect by providing additional hydrodynamic damping and restoring forces. As a result, the system

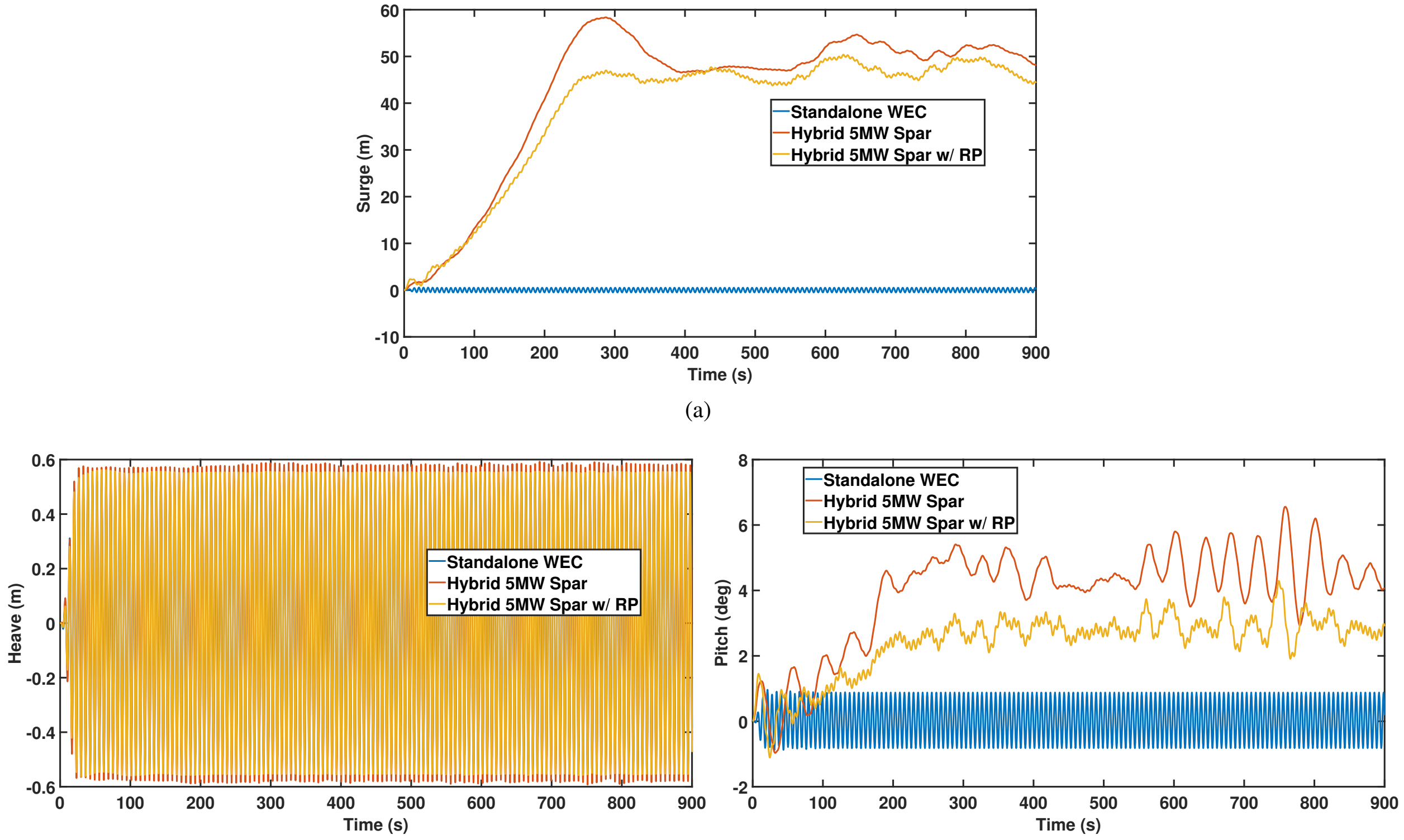

(a)

(b)

(c)

**Figure 15:** Time series for the float response in (a) surge (b) heave (c) pitch directions in a standalone system and hybrid system with 5MW WT on a spar platform with and without the reaction plate under regular wave condition of 6.5 s and 1.75 m

stabilizes more rapidly, with a smaller deviation from the equilibrium position compared to the standalone WEC. These findings underscore the significant impact of hybrid integration on the float's vertical dynamics, driven by the complex hydrodynamic and aerodynamic interactions between the WEC and the WT systems. The reaction plate plays a critical role in amplifying this effect by altering the system's damping characteristics, thereby enabling greater heave motion and, consequently, improved power generation from the WEC.

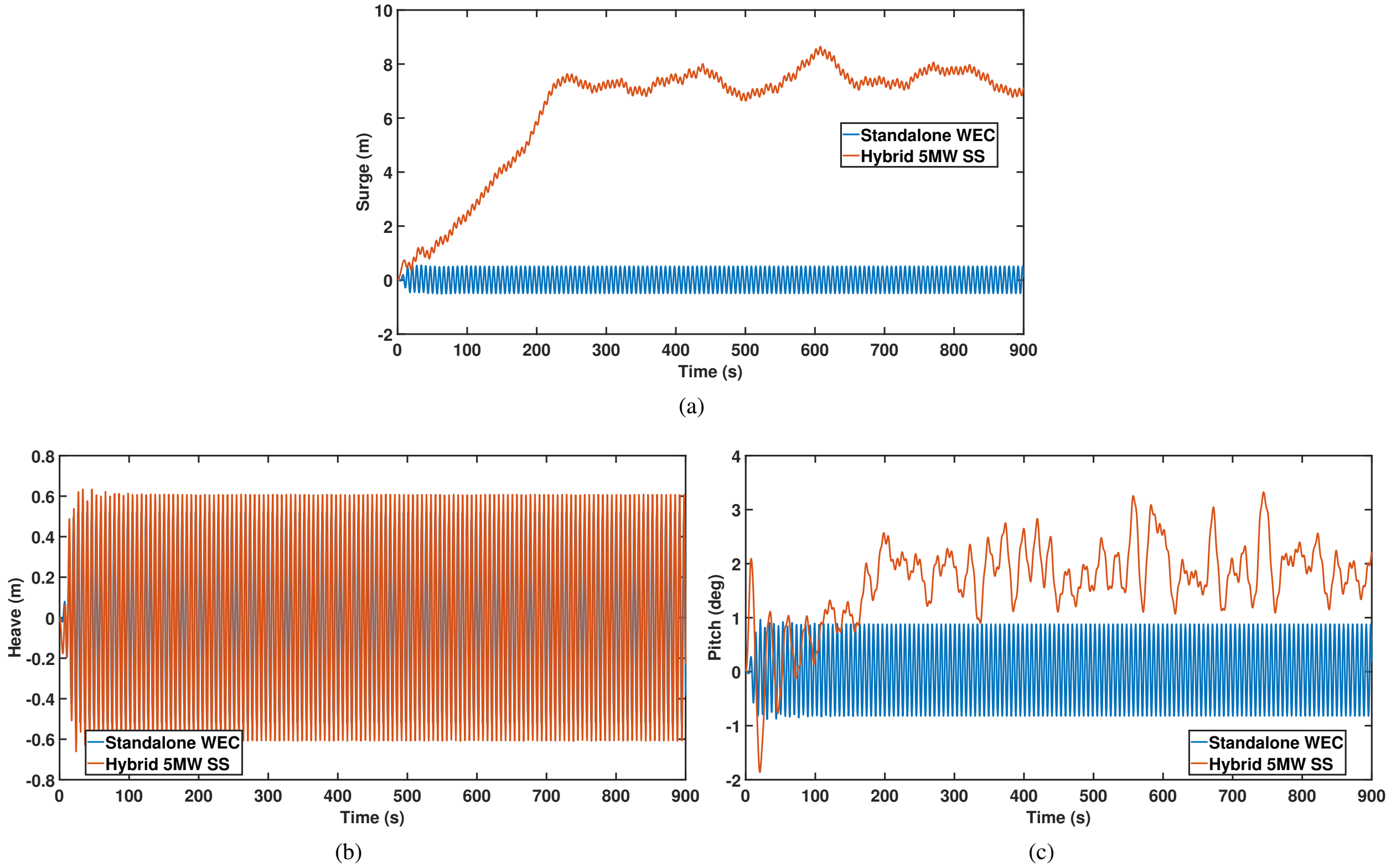

**Figure 16:** Time series for the float response in (a) surge (b) heave (c) pitch directions in a standalone system and hybrid system with 5MW WT on a semi-submersible platform under regular wave condition of 6.5 s and 1.75 m

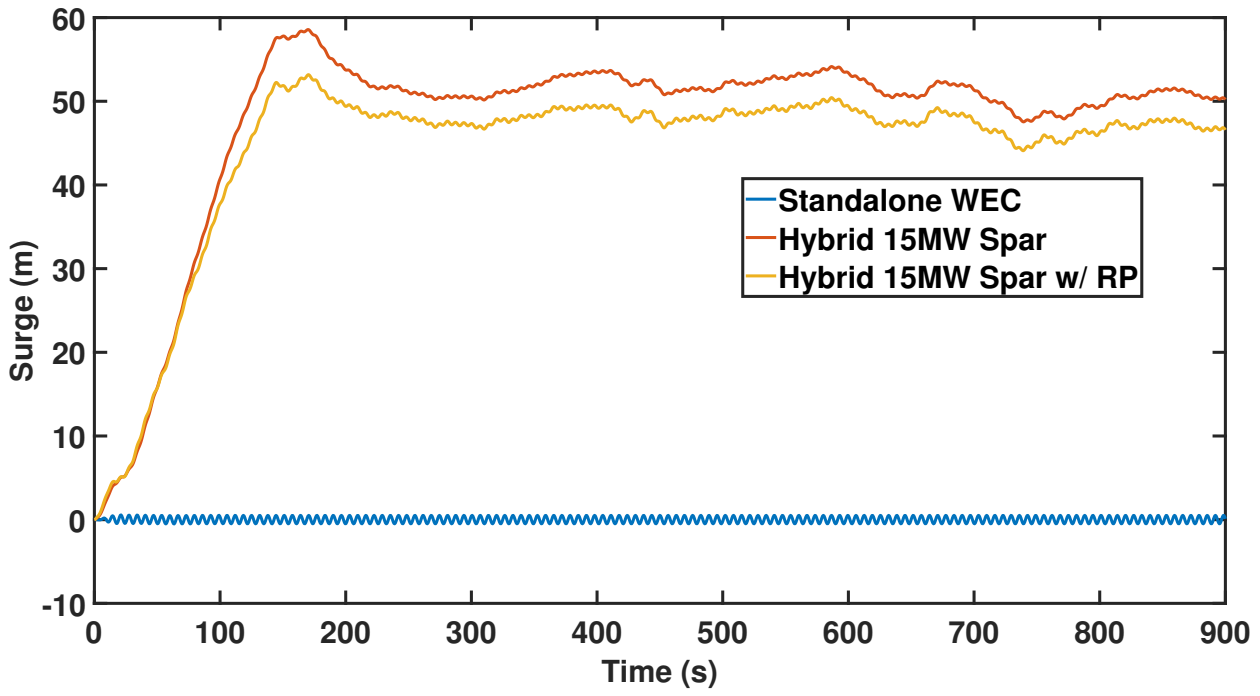

(a)

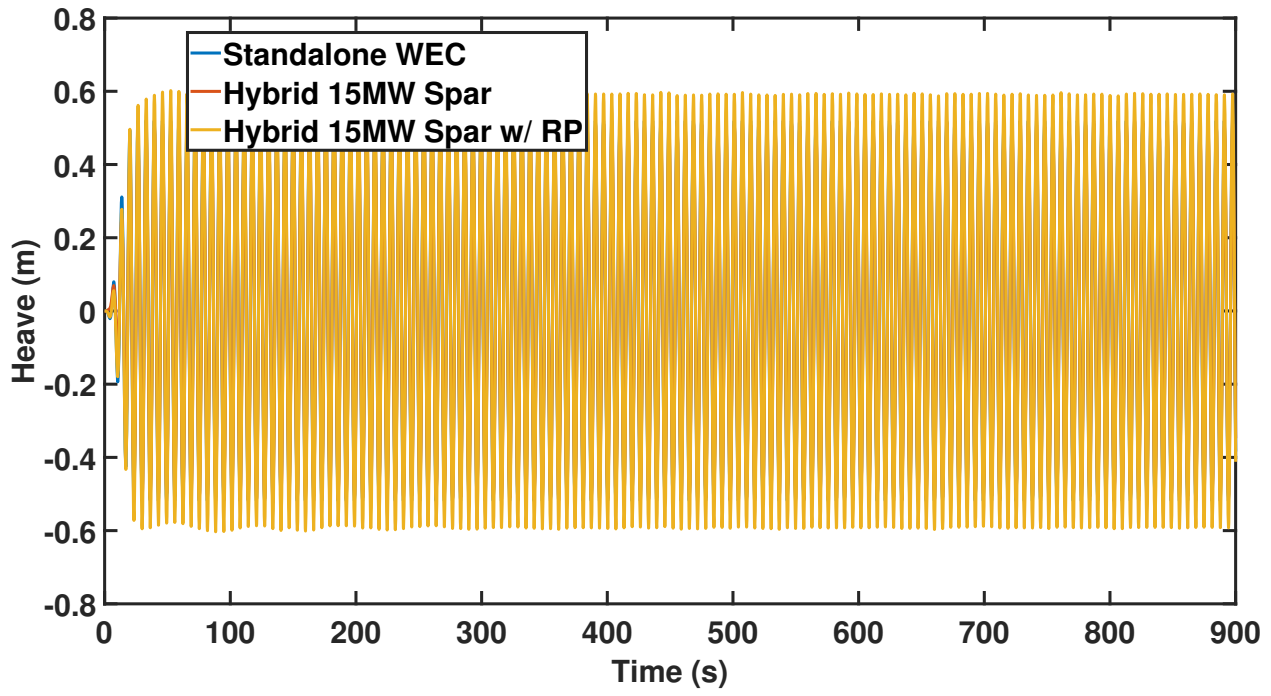

(b)

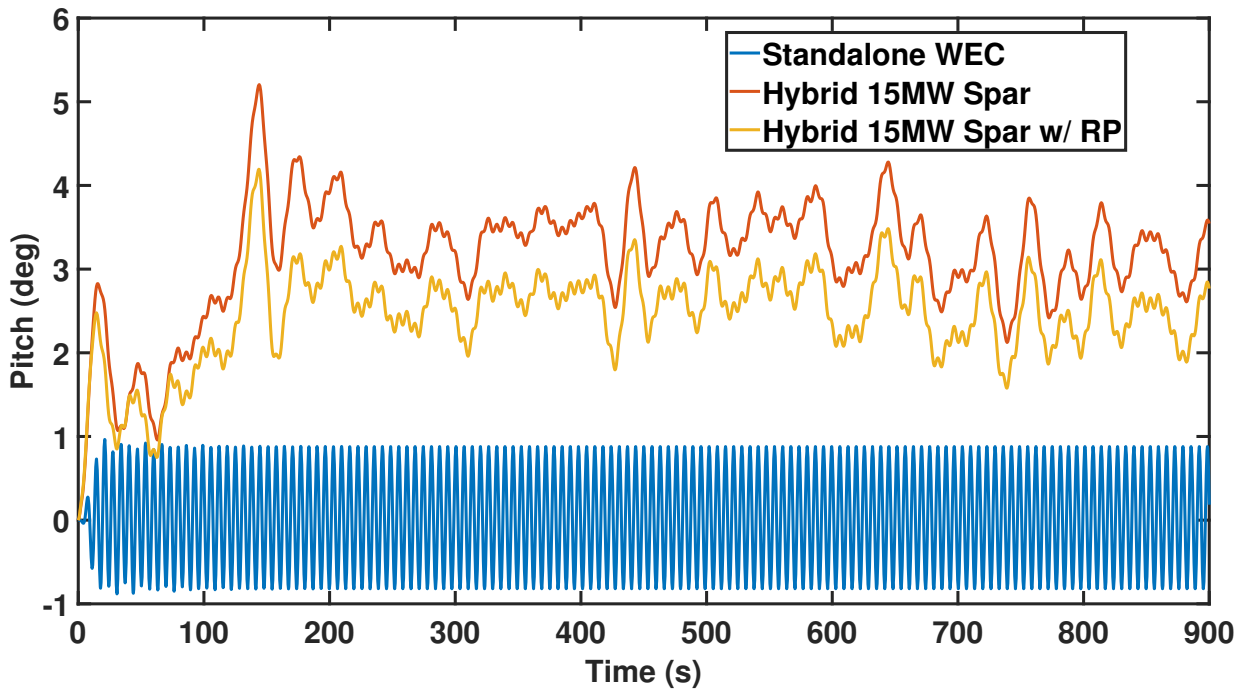

(c)

**Figure 17:** Time series for the float response in (a) surge (b) heave (c) pitch directions in a standalone system and hybrid system with 15MW WT on a spar platform with and without the reaction plate under regular wave condition of 6.5 s and 1.75 m

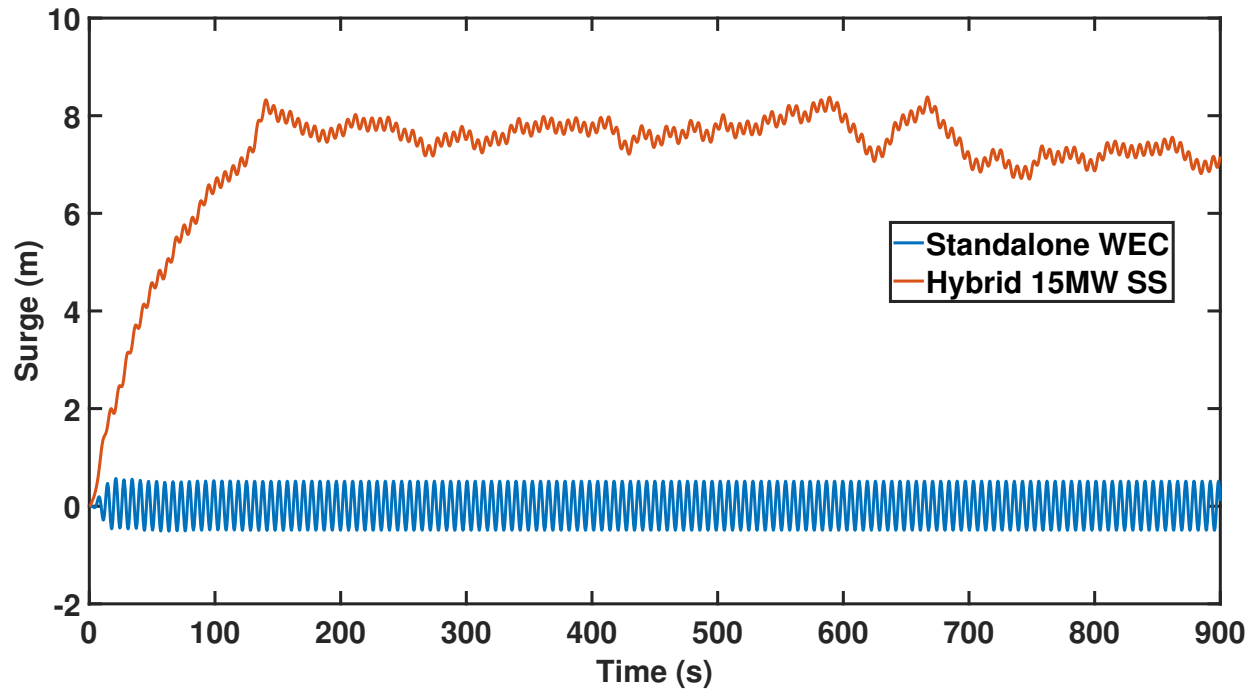

(a)

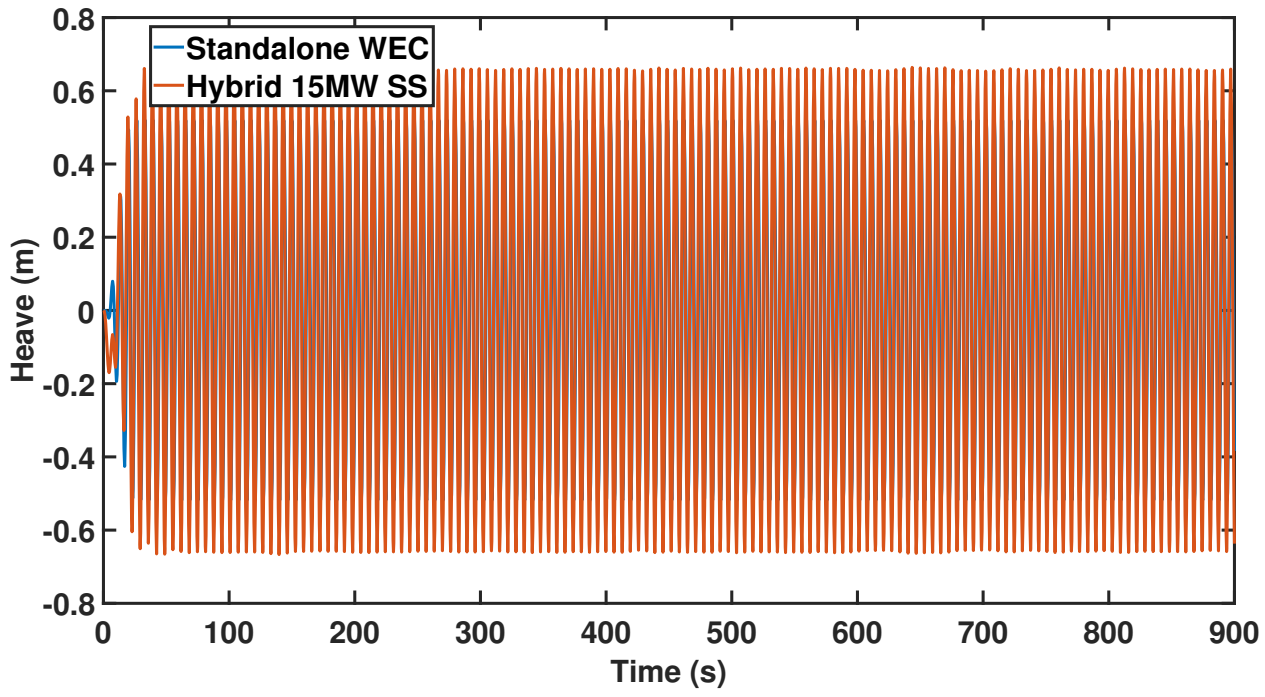

(b)

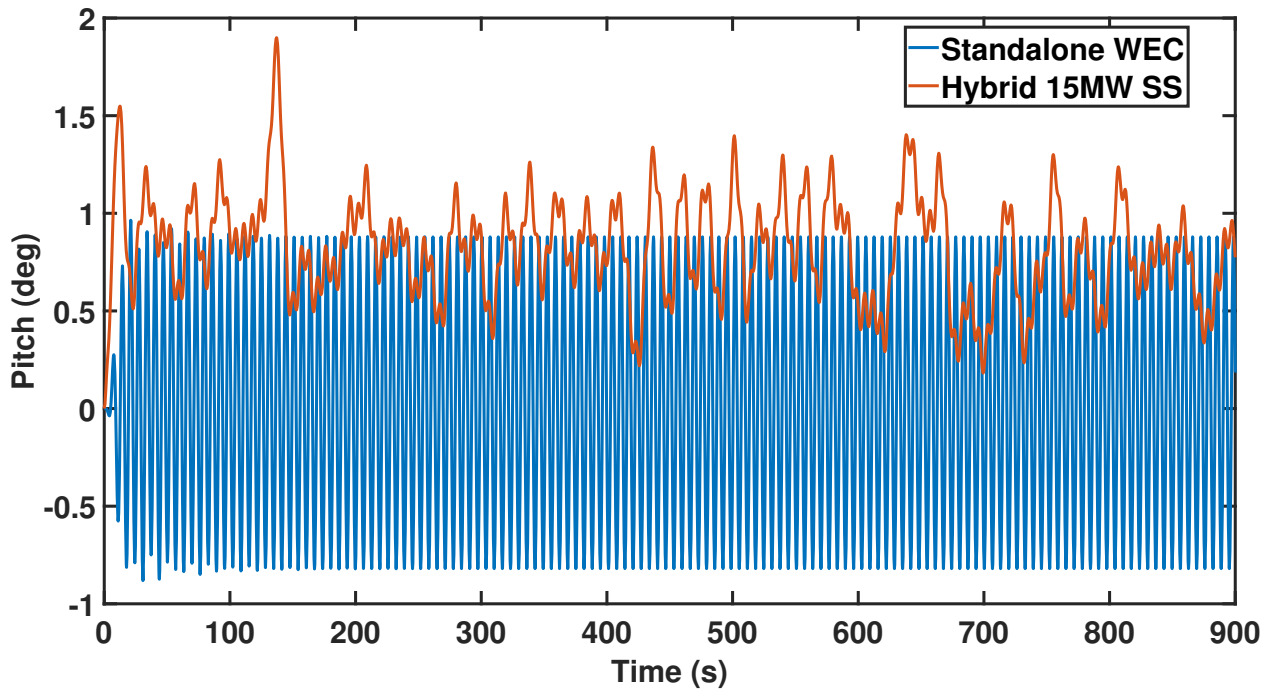

(c)

**Figure 18:** Time series for the float response in (a) surge (b) heave (c) pitch directions in a standalone system and hybrid system with 15MW WT on a semi-submersible platform under regular wave condition of 6.5 s and 1.75 m

**Table 9**
RMS of the platform response and mooring forces for a standalone system of 5 MW WT under different regular wave conditions with wind speed of 11.4 m/s

| Degree of Freedom | X | | | Y | | | Z | | | RX | | | RY | | | RZ | | |
|---|---|---|---|---|---|---|---|---|---|---|---|---|---|---|---|---|---|---|
| Run number | Run 1 | Run 2 | Run 3 | Run 1 | Run 2 | Run 3 | Run 1 | Run 2 | Run 3 | Run 1 | Run 2 | Run 3 | Run 1 | Run 2 | Run 3 | Run 1 | Run 2 | Run 3 |
| Spar mooring platform (m or deg) | 2.29 | 3.23 | 3.30 | 0.46 | 0.64 | 0.75 | 0.08 | 0,16 | 0.47 | 0.18 | 0.22 | 0.56 | 0.57 | 0.71 | 1.15 | 0.50 | 0.62 | 0.95 |
| Semi-submersible mooring platform (m or deg) | 0.42 | 0.86 | 1.89 | 0.07 | 0.18 | 0.29 | 0.03 | 0.29 | 0.27 | 0.15 | 0.24 | 0.28 | 0.53 | 0.72 | 0.79 | 0.31 | 0.79 | 0.72 |
| Mooring forces with spar platform (kN or kNm) | 39.80 | 54.00 | 55.55 | 16.12 | 21.84 | 40.77 | 19.30 | 26.45 | 28.51 | 1111.70 | 1482.1 | 3179.60 | 2768.90 | 3613.20 | 3950.70 | 153.46 | 180.59 | 508.30 |
| Mooring forces with semi-submersible platform (kN or kNm) | 44.52 | 78.01 | 228.63 | 3.59 | 8.44 | 13.91 | 10.73 | 20.78 | 64.03 | 158.97 | 278.73 | 327.70 | 811.24 | 1503.70 | 2481.90 | 513.92 | 1315.40 | 1227.00 |

**Table 10**
RMS of the platform response and mooring forces for a standalone system of 5 MW WT under different irregular wave conditions with wind speed of 11.4 m/s

| Degree of Freedom | X | | | Y | | | Z | | | RX | | | RY | | | RZ | | |
|---|---|---|---|---|---|---|---|---|---|---|---|---|---|---|---|---|---|---|
| Run number | Run 1 | Run 2 | Run 3 | Run 1 | Run 2 | Run 3 | Run 1 | Run 2 | Run 3 | Run 1 | Run 2 | Run 3 | Run 1 | Run 2 | Run 3 | Run 1 | Run 2 | Run 3 |
| Spar mooring platform (m or deg) | 3.73 | 3.89 | 4.04 | 0.89 | 0.86 | 0.82 | 0.12 | 0.15 | 0.29 | 0.27 | 0.26 | 0.26 | 0.63 | 0.66 | 0.86 | 0.48 | 0.55 | 0.65 |
| Semi-submersible mooring platform (m or deg) | 1.07 | 0.94 | 1.35 | 0.44 | 0.42 | 0.40 | 0.04 | 0.14 | 0.80 | 0.34 | 0.32 | 0.28 | 0.58 | 0.62 | 0.72 | 0.87 | 0.83 | 0.76 |
| Mooring forces with spar platform (kN or kNm) | 61.15 | 63.66 | 66.09 | 29.09 | 28.25 | 27.66 | 29.57 | 30.94 | 32.44 | 1947.90 | 1885.7 | 1857.60 | 3976.40 | 4129.90 | 4341.70 | 218.50 | 209.18 | 210.06 |
| Mooring forces with semi-submersible platform (kN or kNm) | 83.26 | 76.28 | 149.07 | 21.06 | 20.27 | 19.25 | 8.31 | 12.71 | 53.67 | 375.87 | 357.32 | 317.23 | 862.35 | 1034.40 | 1784.10 | 1452.80 | 1382.40 | 1272.00 |

**Table 11**
RMS of the platform response and mooring forces for a standalone system of 15 MW WT under different regular wave conditions with wind speed of 10.6 m/s

| Degree of Freedom | X | | | Y | | | Z | | | RX | | | RY | | | RZ | | |
|---|---|---|---|---|---|---|---|---|---|---|---|---|---|---|---|---|---|---|
| Run number | Run 1 | Run 2 | Run 3 | Run 1 | Run 2 | Run 3 | Run 1 | Run 2 | Run 3 | Run 1 | Run 2 | Run 3 | Run 1 | Run 2 | Run 3 | Run 1 | Run 2 | Run 3 |
| Spar mooring platform (m or deg) | 1.48 | 1.74 | 2.08 | 0.55 | 0.59 | 1.44 | 0.25 | 0.26 | 0.68 | 0.09 | 0.15 | 0.27 | 0.47 | 0.57 | 0.94 | 0.32 | 0.34 | 0.37 |
| Semi-submersible mooring platform (m or deg) | 0.36 | 0.81 | 1.73 | 0.05 | 0.13 | 0.16 | 0.07 | 0.17 | 0.69 | 0.04 | 0.04 | 0.16 | 0.26 | 0.71 | 1.26 | 0.32 | 0.56 | 0.83 |
| Mooring forces with spar platform (kN or kNm) | 168.46 | 197.52 | 231.24 | 2.61 | 3.20 | 7.63 | 31.71 | 38.72 | 45.59 | 222.23 | 288.74 | 612.50 | 11380.71 | 13319.47 | 15646.95 | 175.10 | 266.36 | 469.70 |
| Mooring forces with semi-submersible platform (kN or kNm) | 90.67 | 195.74 | 312.38 | 6.85 | 19.72 | 24.44 | 45.23 | 51.44 | 117.61 | 385.97 | 629.69 | 1272.39 | 4230.11 | 3085.29 | 8467.45 | 3122.09 | 5461.90 | 8383.04 |

**Table 12**
RMS of the platform response and mooring forces for a standalone system of 15 MW WT under different irregular wave conditions with wind speed of 10.6 m/s

| Degree of Freedom | X | | | Y | | | Z | | | RX | | | RY | | | RZ | | |
|---|---|---|---|---|---|---|---|---|---|---|---|---|---|---|---|---|---|---|
| Run number | Run 1 | Run 2 | Run 3 | Run 1 | Run 2 | Run 3 | Run 1 | Run 2 | Run 3 | Run 1 | Run 2 | Run 3 | Run 1 | Run 2 | Run 3 | Run 1 | Run 2 | Run 3 |
| Spar mooring platform (m or deg) | 3.52 | 3.48 | 3.39 | 2.78 | 2.82 | 2.86 | 0.42 | 0.43 | 0.62 | 0.27 | 0.27 | 0.19 | 0.61 | 0.63 | 0.76 | 0.32 | 0.34 | 0.35 |
| Semi-submersible mooring platform (m or deg) | 2.48 | 2.08 | 1.76 | 0.59 | 0.52 | 0.41 | 0.15 | 0.16 | 0.87 | 0.15 | 0.16 | 0.13 | 0.35 | 0.46 | 0.79 | 1.09 | 1.02 | 0.92 |
| Mooring forces with spar platform (kN or kNm) | 369.98 | 366.39 | 359.51 | 14.34 | 14.55 | 14.68 | 650.84 | 65.17 | 65.33 | 983.63 | 997.55 | 970.09 | 24359.87 | 24151.31 | 23738.47 | 469.47 | 474.85 | 342.82 |
| Mooring forces with semi-submersible platform (kN or kNm) | 526.88 | 444.77 | 384.11 | 82.36 | 75.13 | 59.61 | 50.35 | 58.37 | 123.86 | 1604.87 | 1559.82 | 1384.83 | 6979.16 | 6967.05 | 7162.22 | 10311.15 | 9764.81 | 9012.48 |

### *Mooring Platform*

Tables 9-12 show the RMS of the response of the mooring platform and the RMS of the mooring forces in a standalone system. The difference in RMS of the mooring platform's response in the six degrees of freedom between the hybrid system and the standalone system is computed and presented in Figures 19 and 22. Additionally, the tension in the mooring lines was analyzed, comparing the standalone WT system to the hybrid system with the float attached. Notably, the mooring line stiffness, fixed points, and lengths remained unchanged in both cases for each mooring foundation. The RMS differences are presented in Figures 19 and 22. A negative sign in all tables represents a reduction in response compared to the standalone system. The time series for the response of the mooring platform for all configurations are shown in Figures 23-26, while the mooring forces are shown in Figures 27-30

The results indicate that integrating the float with the WT in the hybrid system led to variations in the response across different degrees of freedom, with distinct effects depending on the type of mooring platform.

For the spar platform integrated with the 5MW WT and the float, the platform exhibited varying responses across different wave conditions. Without the reaction plate at the bottom, most degrees of freedom experienced an increase in motion. However, the addition of the reaction plate significantly altered this trend, leading to a substantial reduction in platform motion across all degrees of freedom. The reduction reached up to 30% in surge, 37% in heave, and 33% in pitch, highlighting the plate's effectiveness in mitigating excessive motion. Similarly, for the 15MW WT, the integration of the float led to a notable reduction in platform motion under both regular and irregular wave conditions, with the effect becoming more pronounced when the reaction plate was incorporated. The average reduction observed was 16% in surge, 60% in heave, and 12% in pitch. This demonstrates the reaction plate's critical role in enhancing system stability by increasing added mass and damping characteristics, thereby effectively suppressing undesirable platform motions.

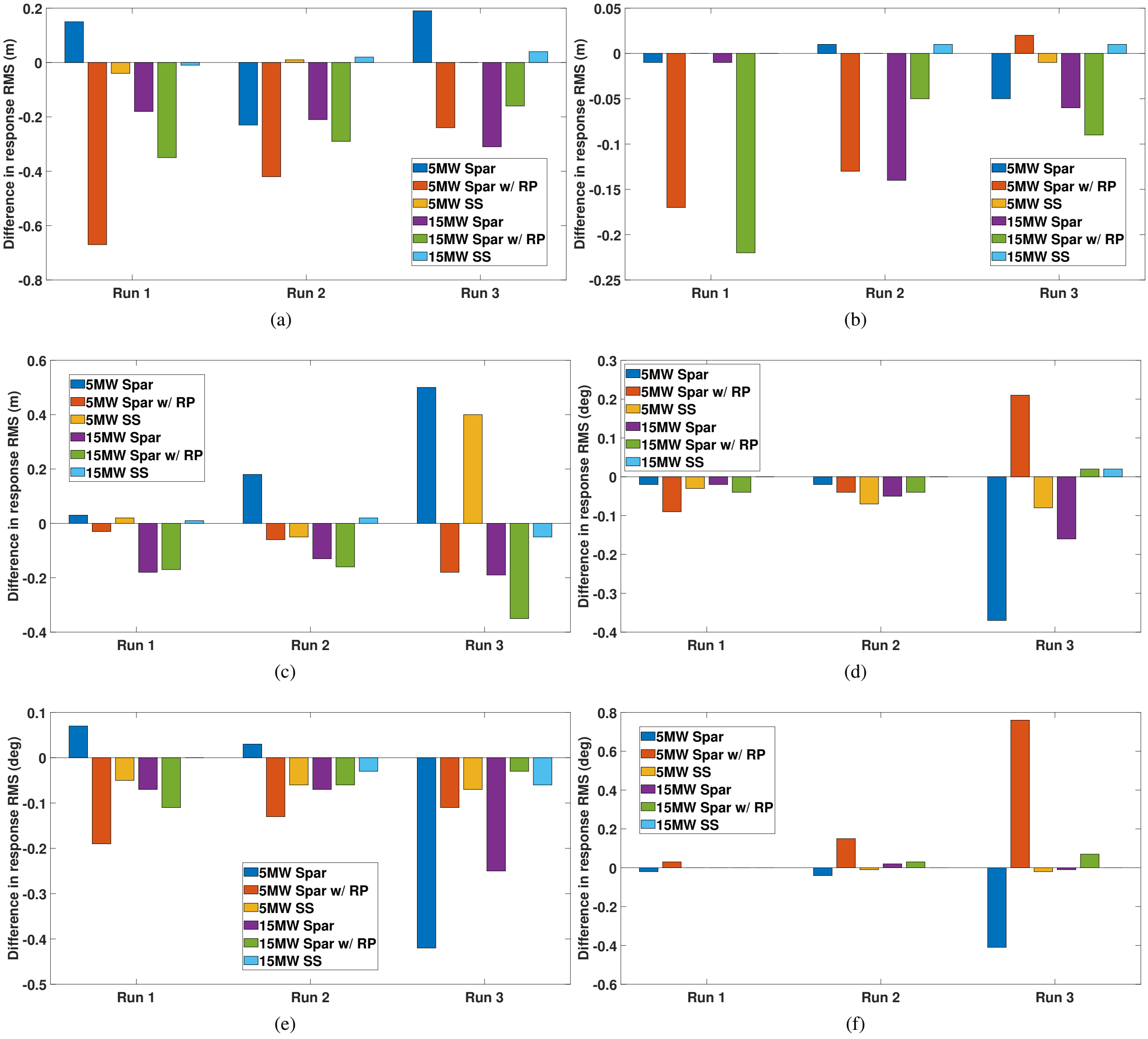

**Figure 19:** Difference in response RMS of the spar platform and sem-submersible platform between standalone and all configuration of hybrid system in (a) surge, (b) sway, (c) heave, (d) roll, (e) pitch, and (f) yaw of 5MW WT with wind speed of 11.4 m/s and 15MW WT with wind speed of 10.6 m/s under different regular wave conditions—Zero value means equivalent to the corresponding standalone system RMS

Similar to the float, the equilibrium position of the platform in the hybrid system shifts from that of the standalone system, as shown in Figures 23-26. This shift is still due to the aerodynamic forces exerted by the WT, which influence the system's response across multiple degrees of freedom. The addition of the reaction plate at the base of the spar platform enhanced the system's stability by providing additional hydrodynamic damping and restoring forces and reduced the change in the equilibrium position of the platform in the hybrid system compared to the standalone system. The reaction plate has an important role in altering the system's response and improving the stability of the mooring platform.

The semi-submersible platform exhibited a slight change in surge motion, with variations ranging from 3% to 5% in most cases. Additionally, it demonstrated a notable reduction in pitch across all wave conditions for both WTs, with reductions ranging from 5% to 10% under regular wave conditions and 3% to 7% under irregular wave conditions. However, the heave response varied between cases, most likely due to the added wave-induced vertical

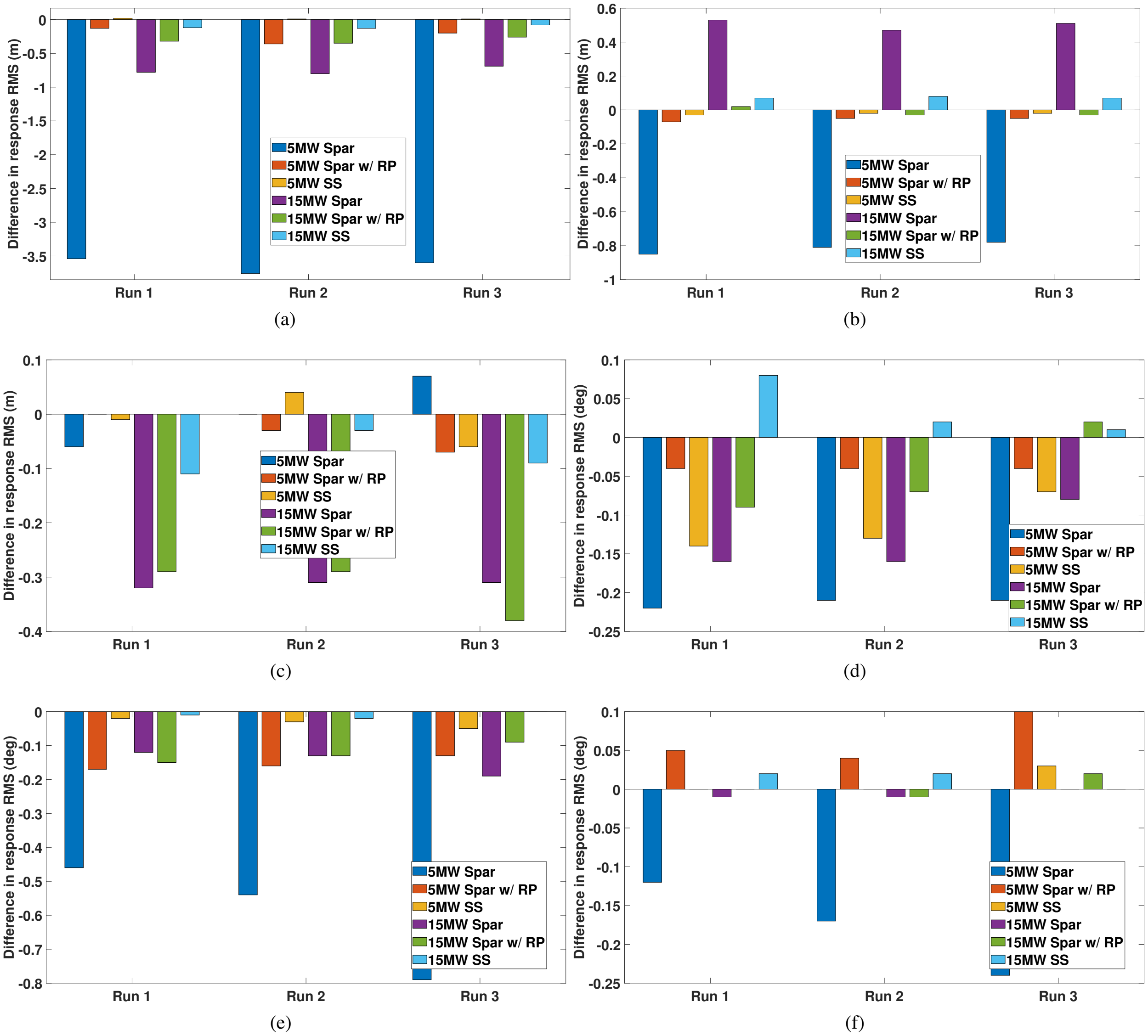

**Figure 20:** Difference in response RMS of the spar platform and sem-submersible platform between standalone and all configuration of hybrid system in (a) surge, (b) sway, (c) heave, (d) roll, (e) pitch, and (f) yaw of 5MW WT with wind speed of 11.4 m/s and 15MW WT with wind speed of 10.6 m/s under different irregular wave conditions—Zero value means equivalent to the corresponding standalone system RMS

motion introduced by the float. Under regular wave conditions, the platform generally experienced an increase in heave motion, whereas under irregular wave conditions, a significant reduction was observed. This variation is attributed to the interaction between the float's heave response and the platform's motion. Since the float exhibits a greater heave response under regular wave conditions compared to irregular ones, its influence on the platform follows a similar trend, amplifying the platform's motion in regular waves while having a dampening impact on it under irregular waves.

The analysis of mooring forces in the hybrid system highlights the significant influence of both platform type and the inclusion of the reaction plate on system stability. For the semi-submersible platform with the 5 MW WT under regular wave conditions, mooring forces exhibited moderate variations ranging from -26% to 25% across all degrees of freedom, with negative values indicating reductions in loads and positive values indicating increases. The corresponding 15 MW semi-submersible configuration showed smaller variations, between -9% and 4%, reflecting the greater inherent stability of the larger system. In contrast, the spar platform displayed more pronounced changes: the 5 MW spar experienced reductions in mooring forces of up to 67% and increases up to 17% without the plate, whereas the 15 MW spar ranged from -30% to 7%. These results underscore the spar's higher sensitivity to hydrodynamic interactions compared to the semi-submersible platform.

Under irregular wave conditions, mooring force fluctuations were amplified due to stochastic wave excitations. For the semi-submersible platform, reductions ranged from 37% for the 5 MW system and up to 7% for the 15 MW system, indicating consistent mitigation of peak loads. The 5 MW spar platform exhibited substantial reductions without the plate, ranging from 92% to 99%, while the addition of the reaction plate limited variations to between -5% and 2%. Similarly, the 15 MW spar platform experienced reductions of 22% without the plate and 10% with the plate. It is observed that the addition of the reaction plate to the 5 MW system significantly enhances stability and reduces mooring loads, particularly in surge, heave, and pitch under regular wave conditions. However, its effect under irregular waves is limited, providing minimal improvement. For the 15 MW system, the reaction plate contributes little to no additional stability under both regular and irregular wave conditions, indicating that larger, inherently more stable platforms are less sensitive to this modification.

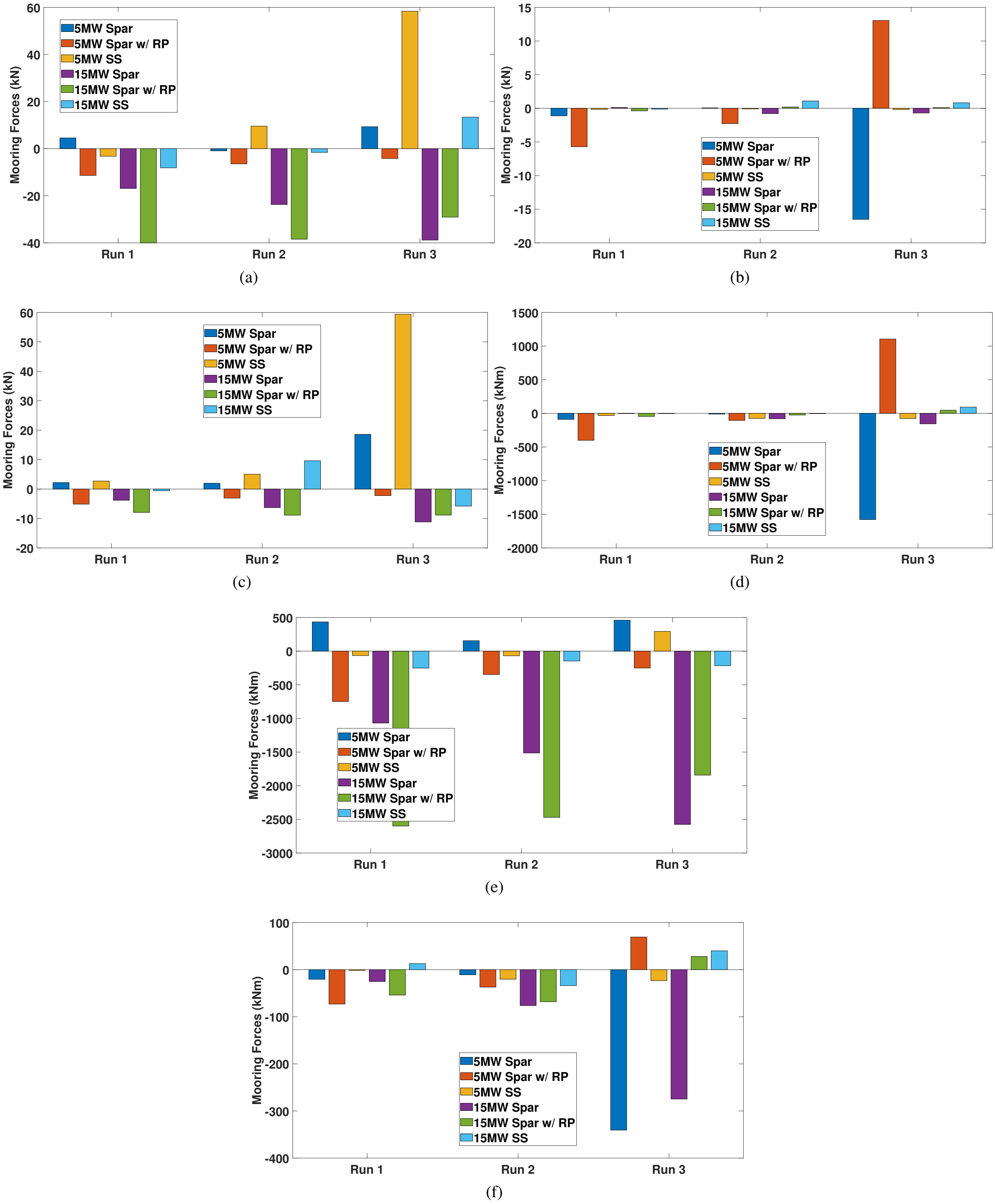

**Figure 21:** Difference in mooring forces RMS of the spar platform and sem-submersible platform between standalone and all configuration of hybrid system in (a) surge, (b) sway, (c) heave, (d) roll, (e) pitch, and (f) yaw of 5MW WT with wind speed of 11.4 m/s and 15MW WT with wind speed of 10.6 m/s under different regular wave conditions—Zero value means equivalent to the corresponding standalone system RMS

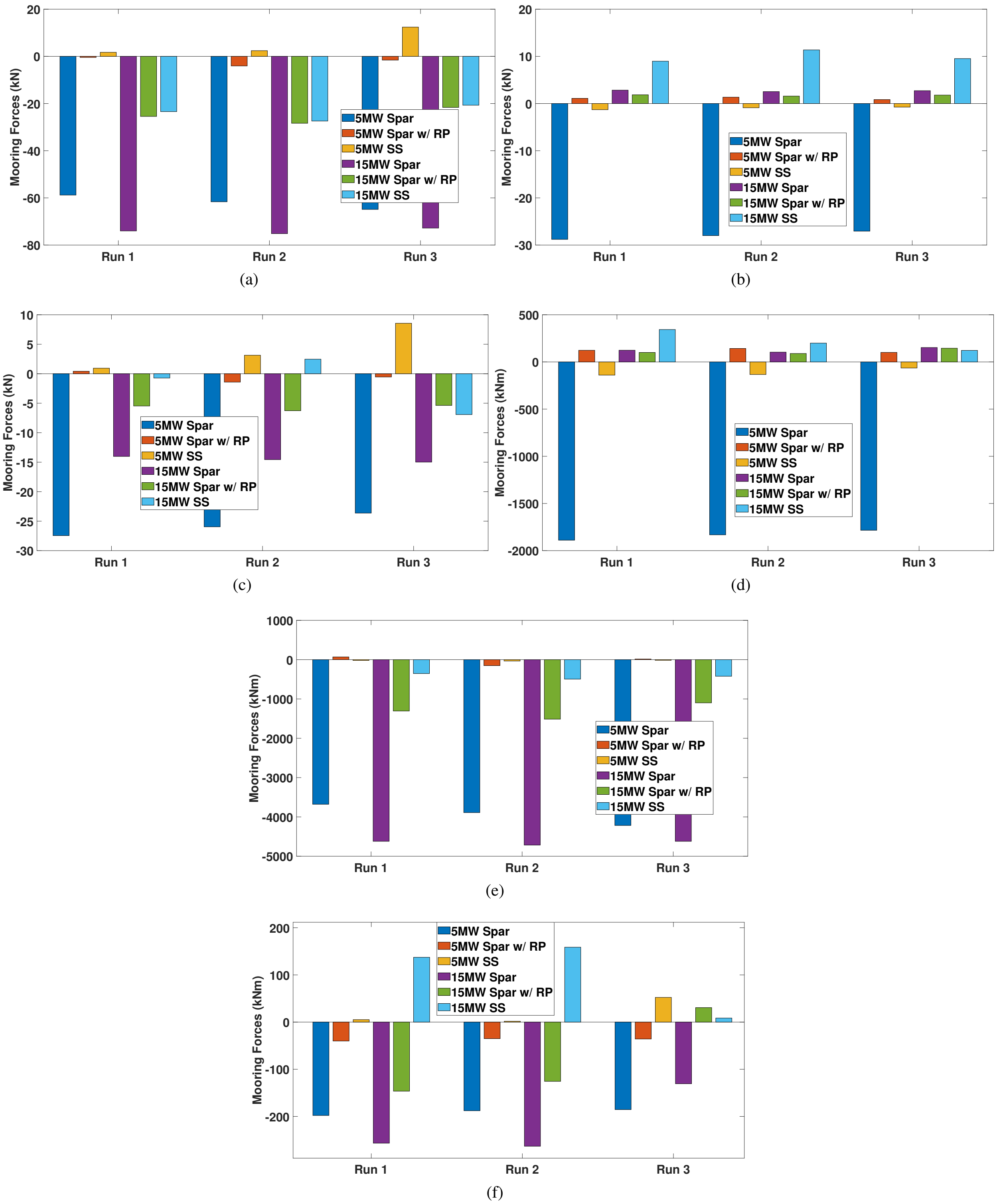

**Figure 22:** Difference in mooring forces RMS of the spar platform and sem-submersible platform between standalone and all configuration of hybrid system in (a) surge, (b) sway, (c) heave, (d) roll, (e) pitch, and (f) yaw of 5MW WT with wind speed of 11.4 m/s and 15MW WT with wind speed of 10.6 m/s under different irregular wave conditions—Zero value means equivalent to the corresponding standalone system RMS

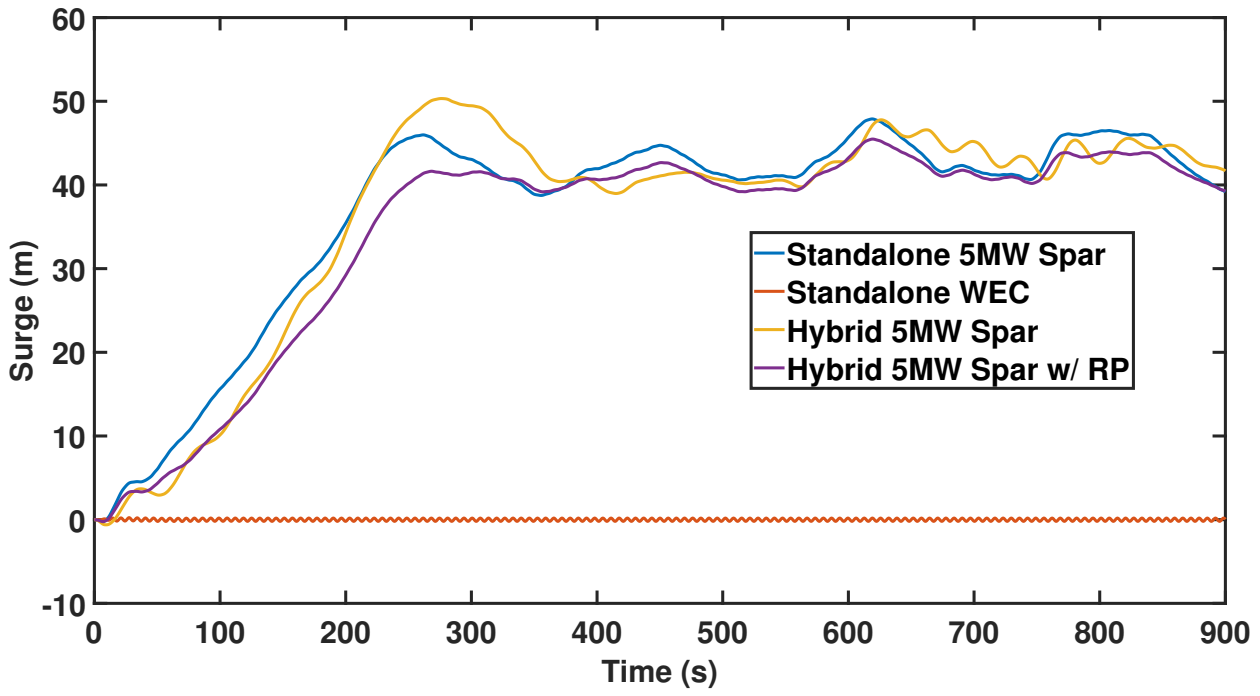

(a)

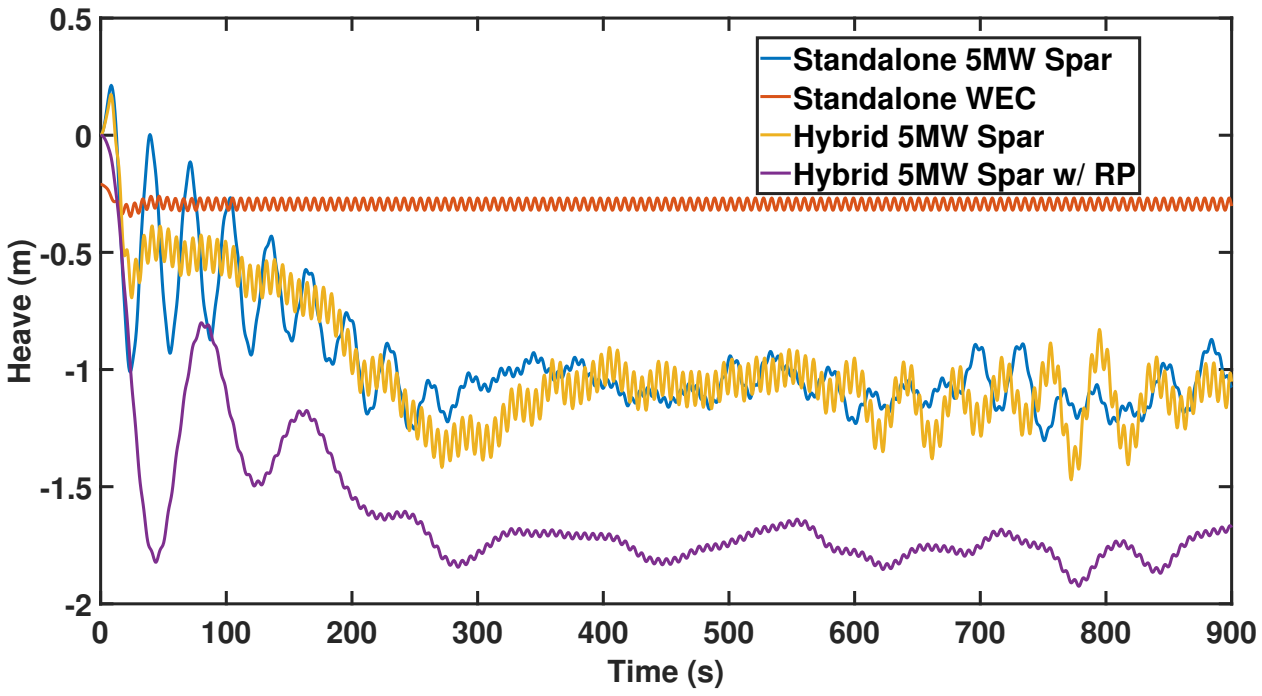

(b)

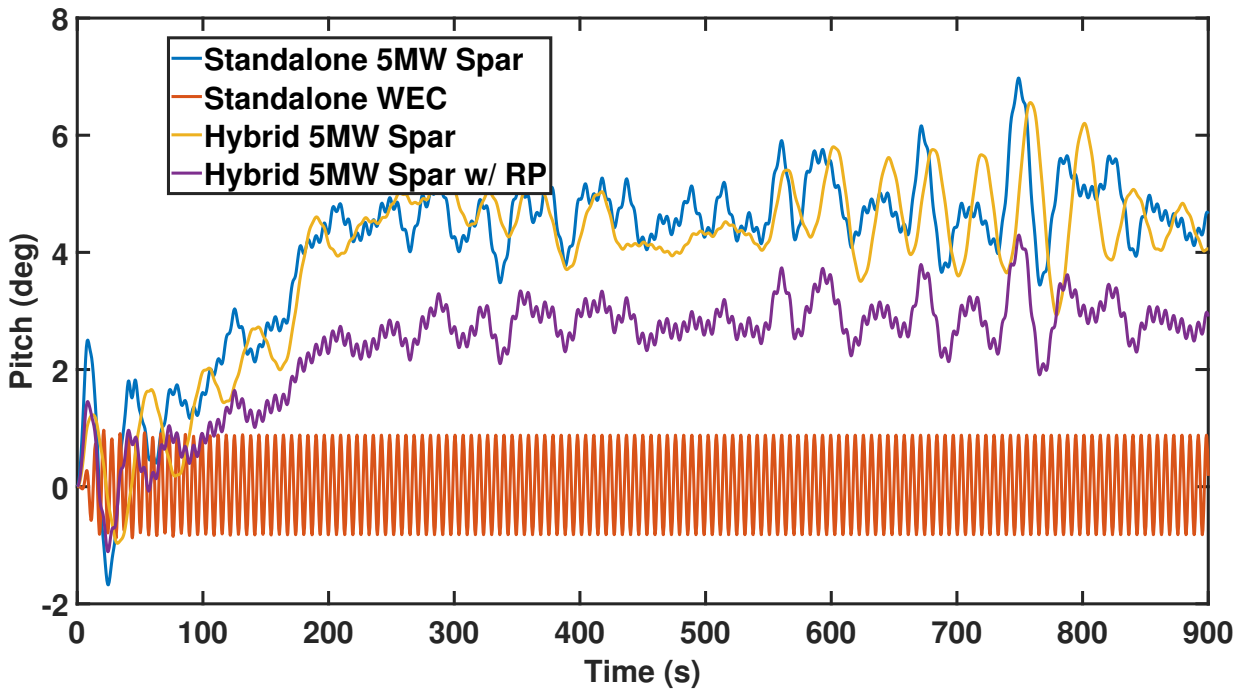

(c)

**Figure 23:** Time series for the spar response in (a) surge (b) heave (c) pitch directions in a standalone system and hybrid system with 5MW WT on a spar platform with and without the reaction plate under regular wave condition of 6.5s and 1.75m

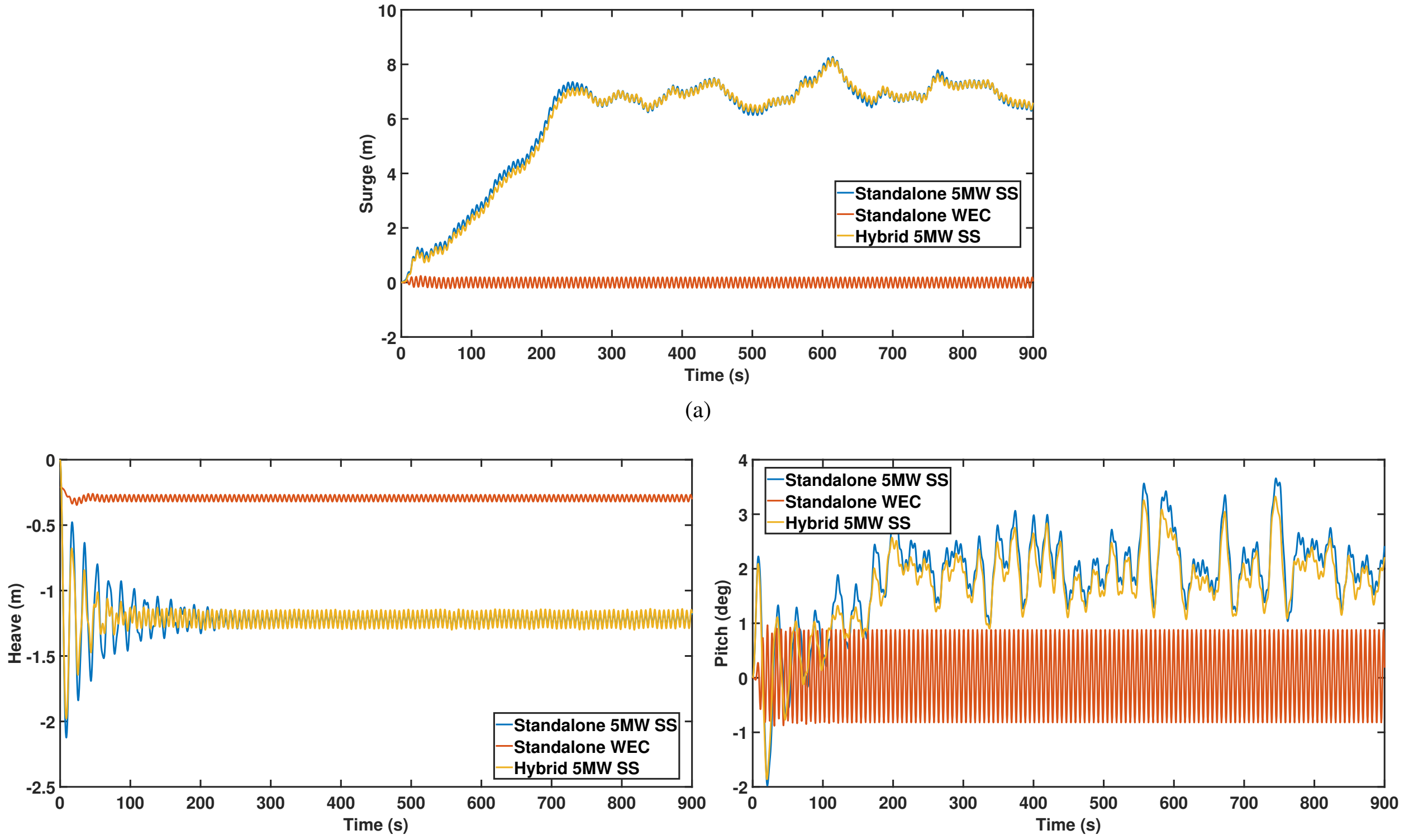

(a)

(b) (c)

**Figure 24:** Time series for the semi-submersible platform response in (a) surge (b) heave (c) pitch directions in a standalone system and hybrid system with 5MW WT on a semi-submersible platform with and without the reaction plate under regular wave condition of 6.5 s and 1.75 m

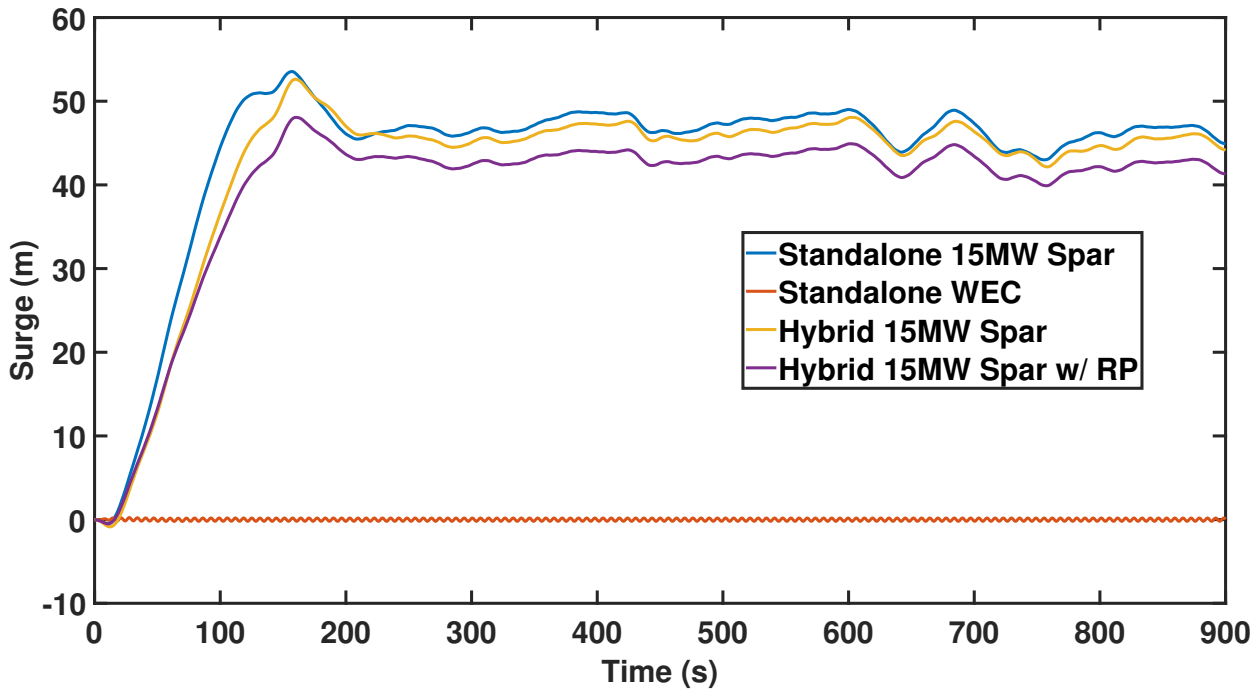

(a)

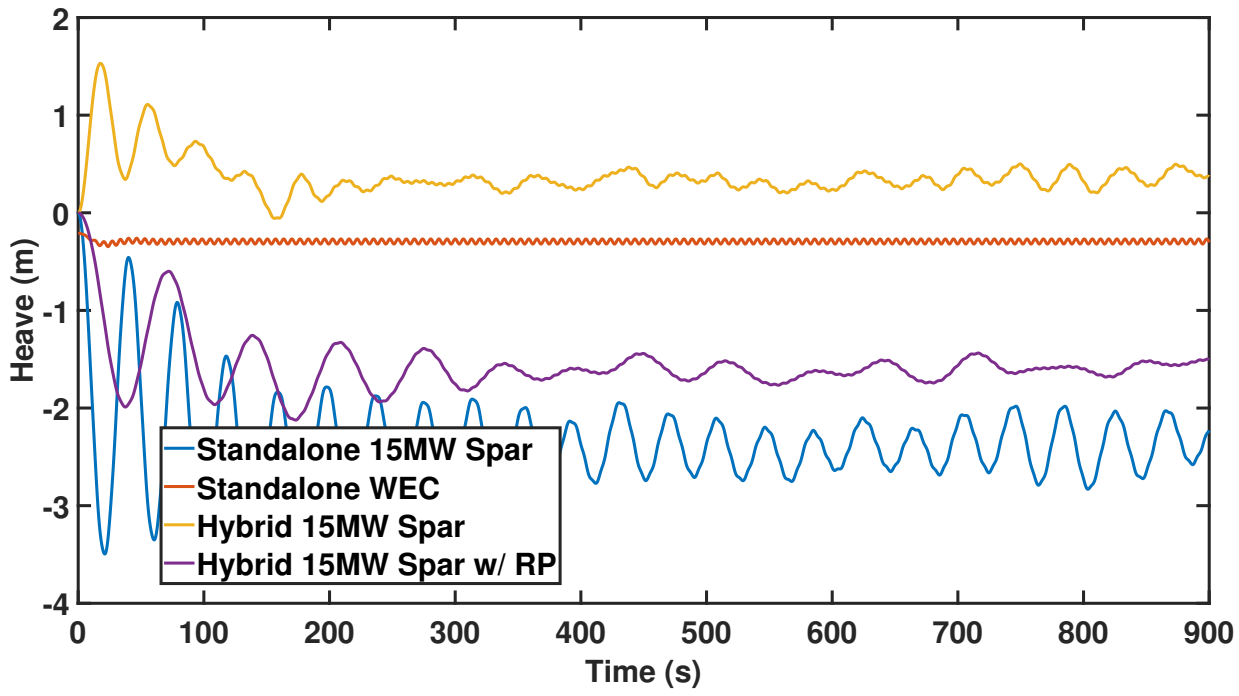

(b)

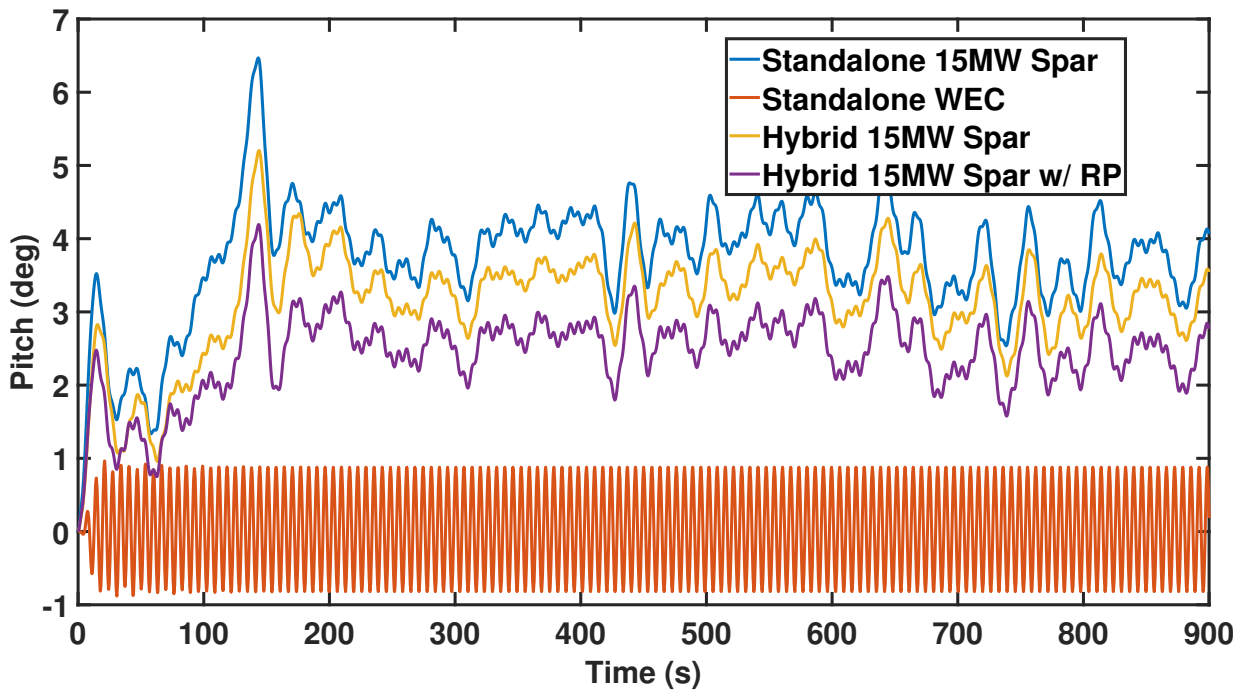

(c)

**Figure 25:** Time series for the spar response in (a) surge (b) heave (c) pitch directions in a standalone system and hybrid system with 15MW WT on a spar platform with and without the reaction plate under regular wave condition of 6.5 s and 1.75 m

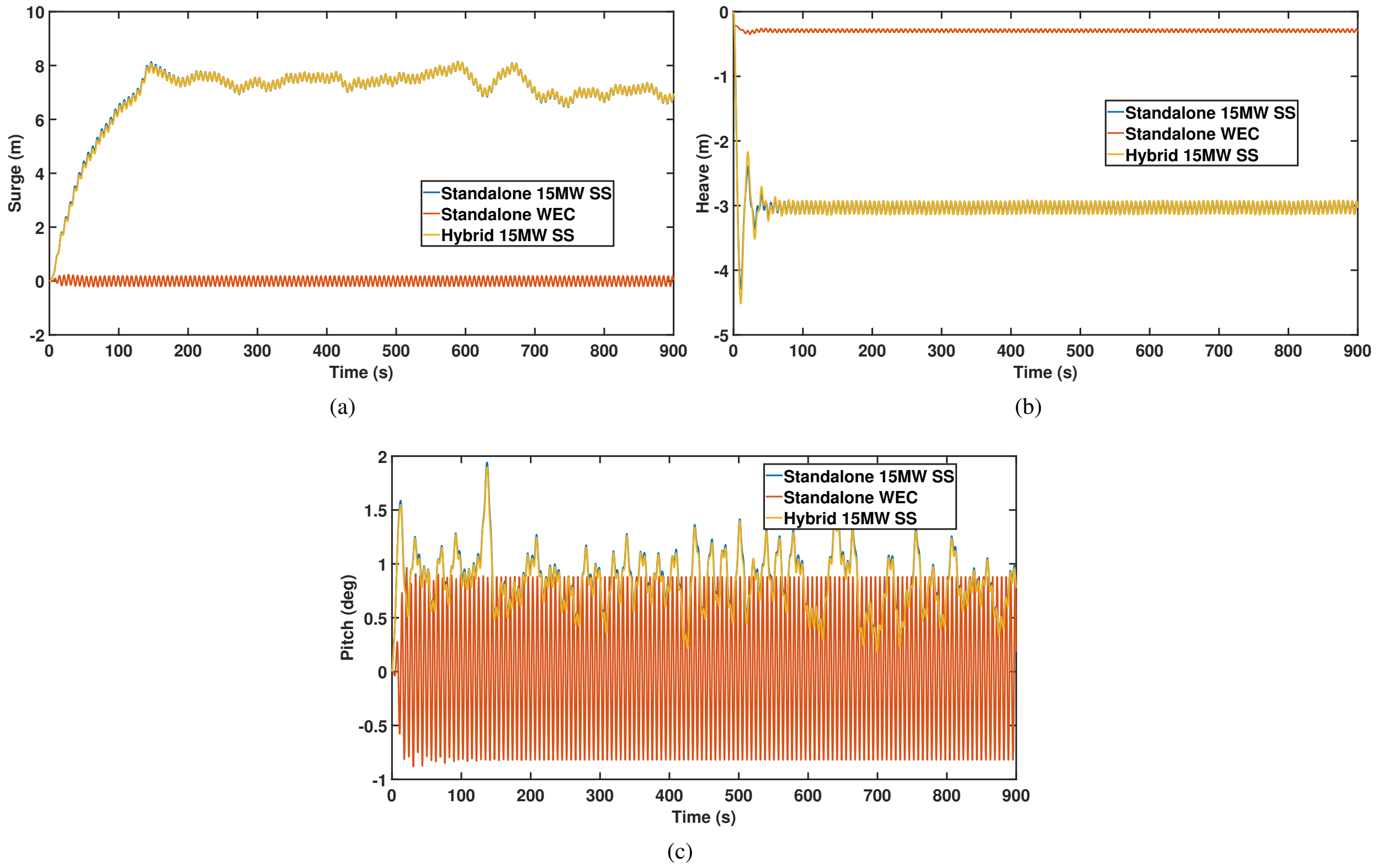

**Figure 26:** Time series for the semi-submersible response in (a) surge (b) heave (c) pitch directions in a standalone system and hybrid system with 15MW WT on a semi-submersible platform under regular wave condition of 6.5s and 1.75m

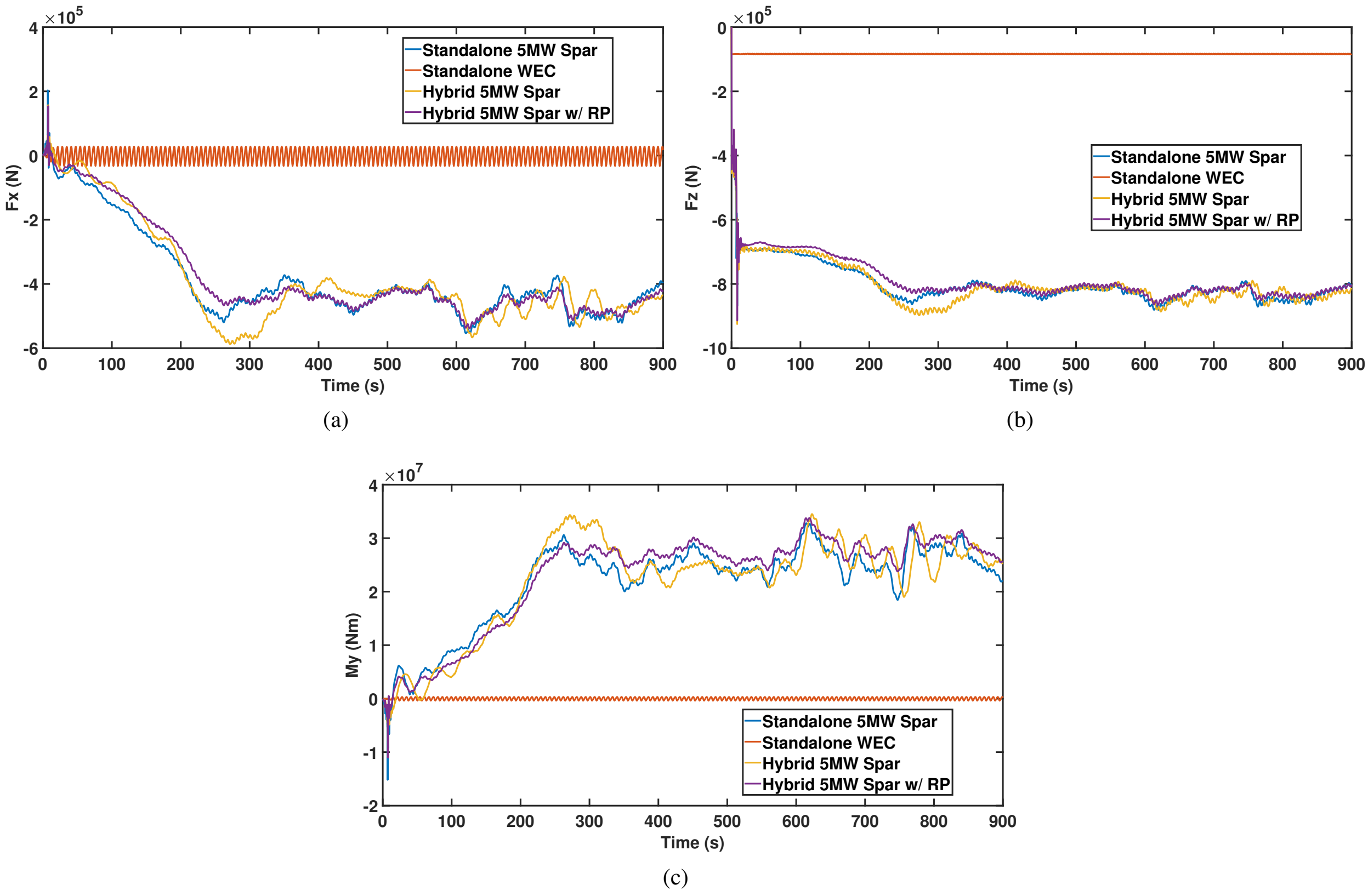

**Figure 27:** Time series for the mooring forces in (a) surge (b) heave (c) pitch directions in a standalone system and hybrid system with 5MW WT on a spar platform with and without the reaction plate under regular wave condition of 6.5s and 1.75m

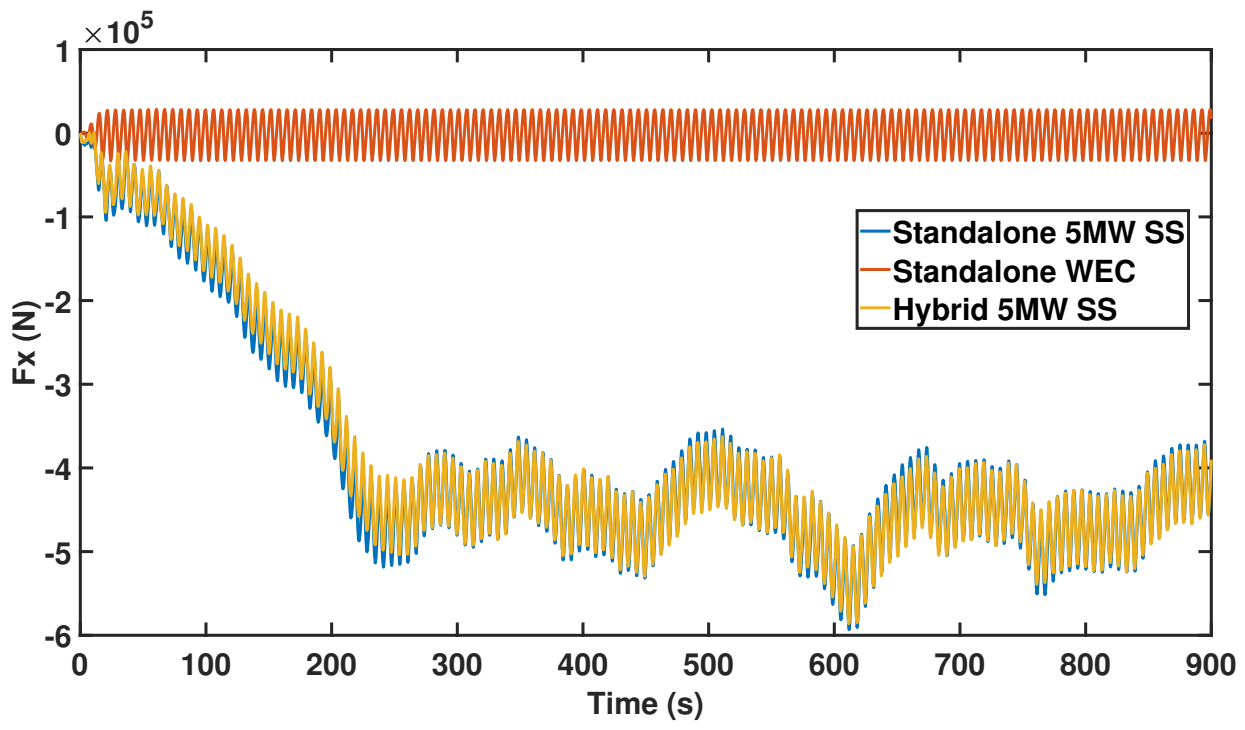

(a)

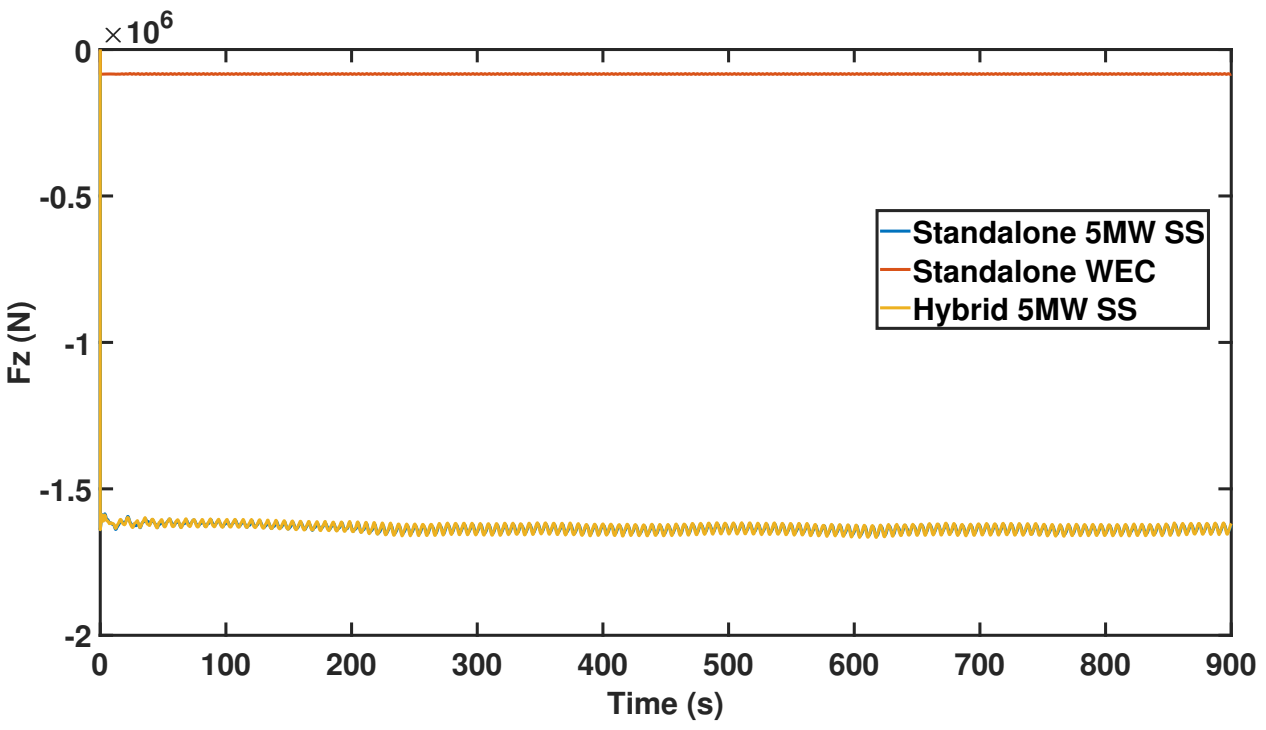

(b)

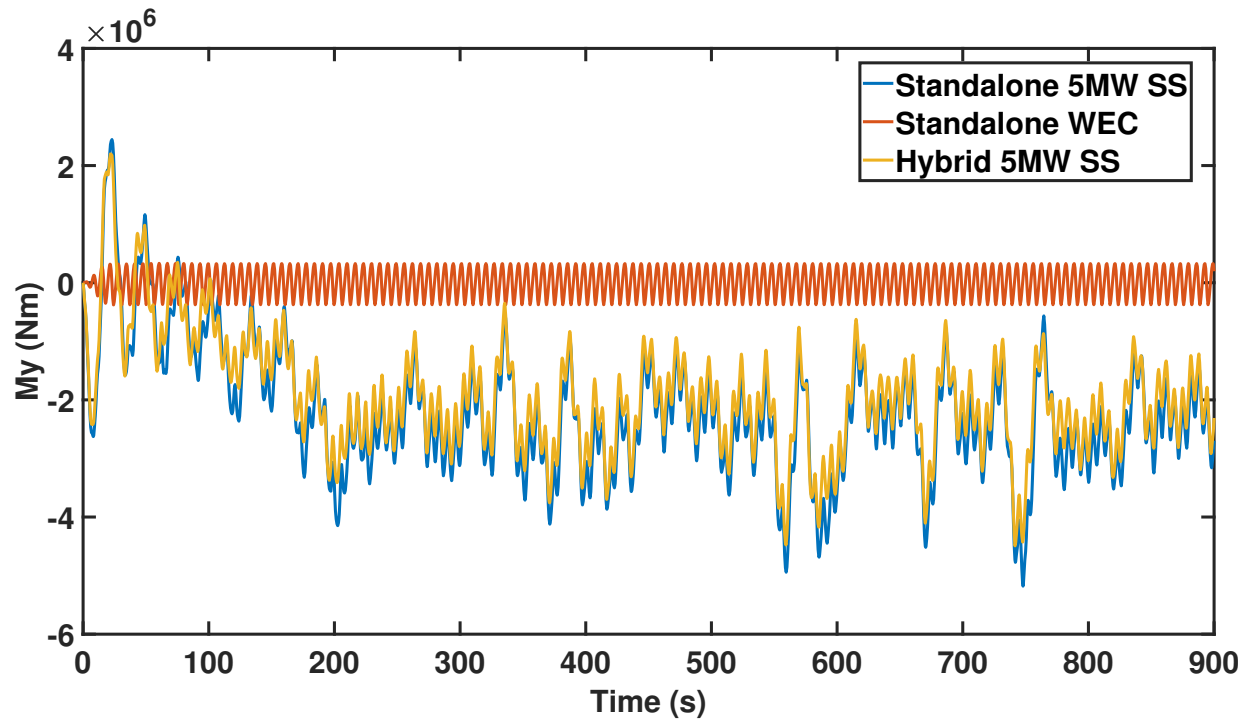

(c)

**Figure 28:** Time series for the mooring forces in (a) surge (b) heave (c) pitch directions in a standalone system and hybrid system with 5MW WT on a semi-submersible platform with and without the reaction plate under regular wave condition of 6.5 s and 1.75 m

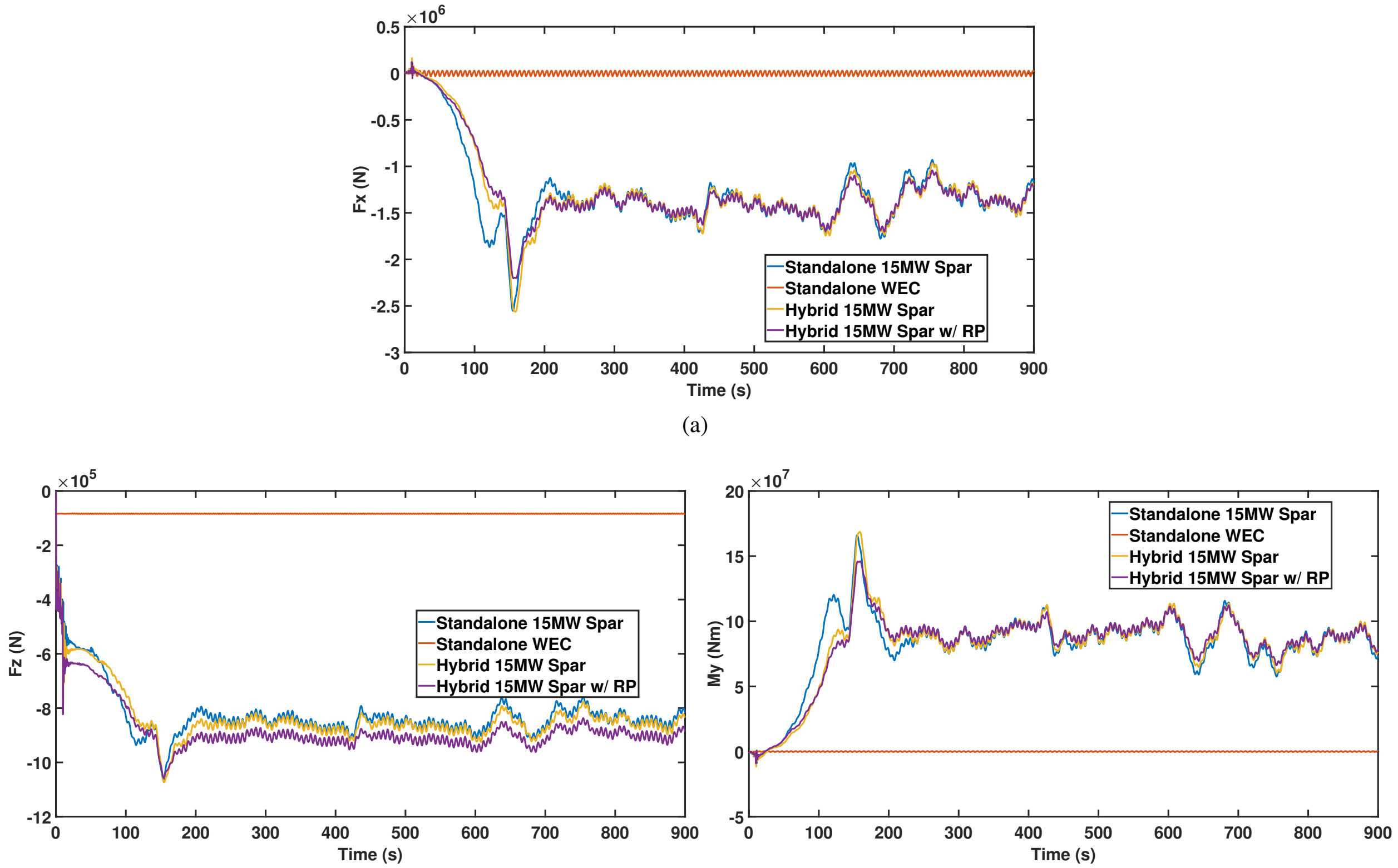

**Figure 29:** Time series for the mooring forces in (a) surge (b) heave (c) pitch directions in a standalone system and hybrid system with 15MW WT on a spar platform with and without the reaction plate under regular wave condition of 6.5 s and 1.75 m

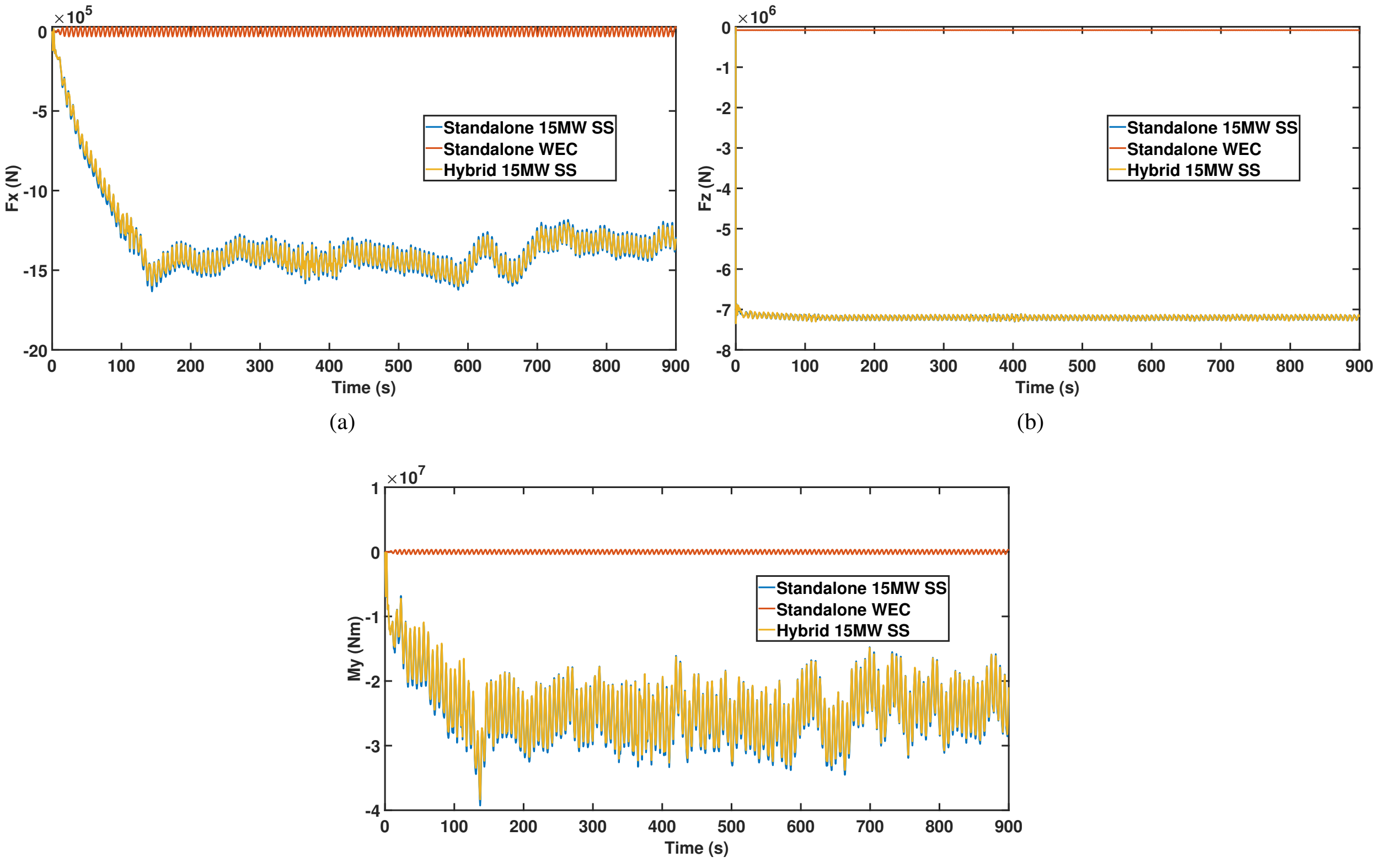

(c)

**Figure 30:** Time series for the mooring forces in (a) surge (b) heave (c) pitch directions in a standalone system and hybrid system with 15MW WT on a semi-submersible platform under regular wave condition of 6.5 s and 1.75 m

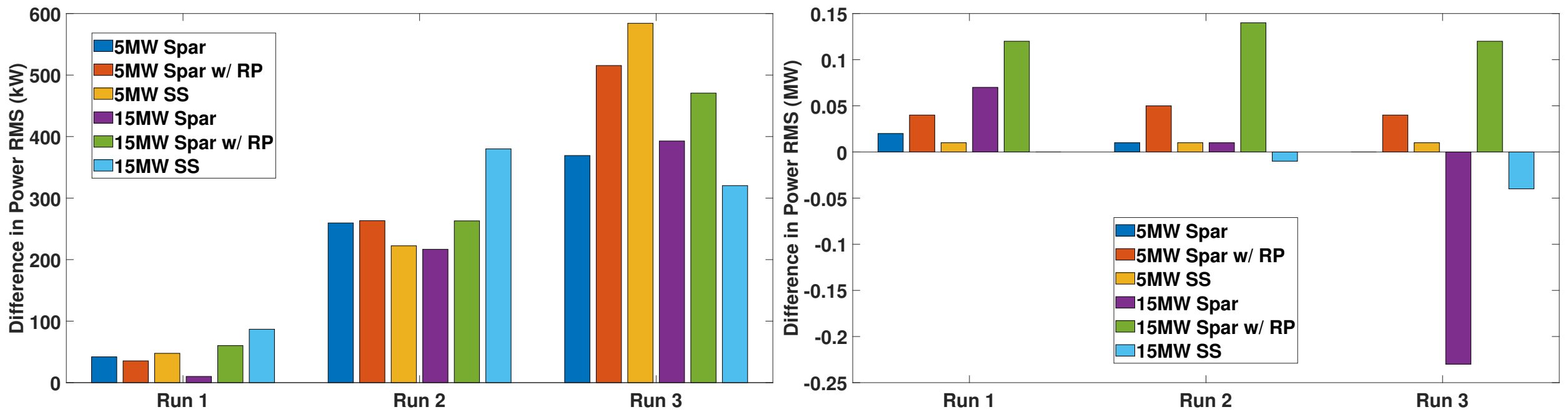

(a) (b)

**Figure 31:** Difference in RMS of power of (a) WEC (b) WT between standalone and hybrid system for different configurations under different regular wave conditions—Zero value means equivalent to the corresponding standalone system RMS

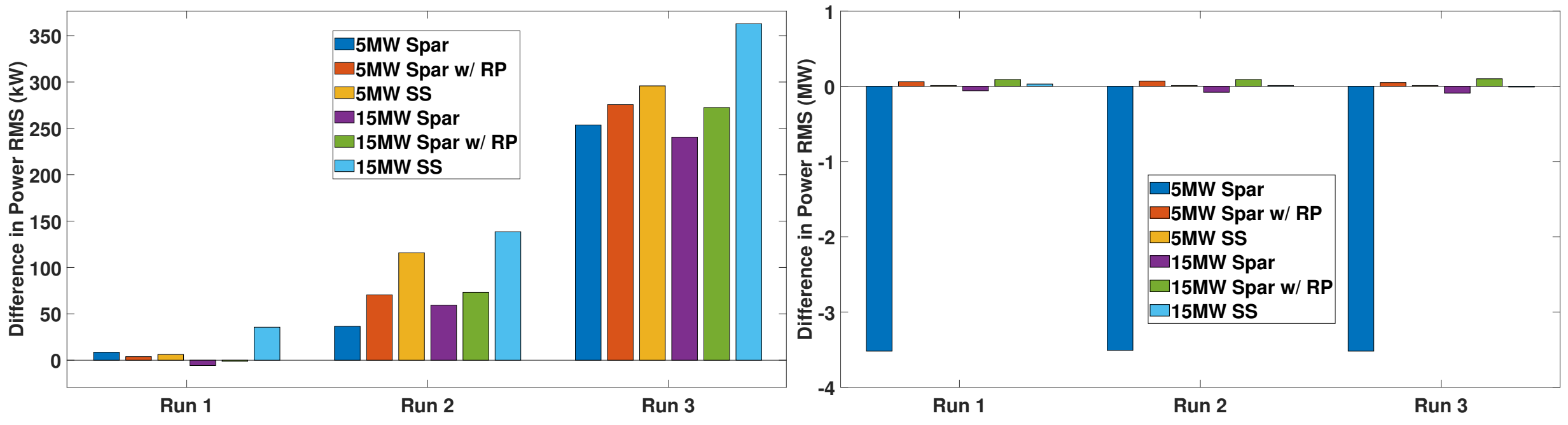

(a) (b)

**Figure 32:** Difference in RMS of power of (a) WEC (b) WT between standalone and hybrid system for different configurations under different irregular wave conditions—Zero value means equivalent to the corresponding standalone system RMS

## 6.2. Power, Coefficient of variation, and AEP

This section reports the power output of the WEC and the WT within the hybrid configurations, along with a comparison to their corresponding standalone systems. Estimates of annual energy generation and variations in power output are also presented.

### *6.2.1. Power*

A comparison of the power generated by the WT in the standalone and hybrid configurations was performed, based on the RMS of the power output, with the results also included in Figure 31 for all six configurations under regular wave conditions and Figure 32 under irregular wave conditions. A negative sign in the table represents a reduction in response compared to the standalone system.

The results reveal that the integration of the float had an insignificant impact on the power generation of the WT, with variations remaining below 1% across all configurations and wave conditions. In contrast, the power generated by the float increased consistently for all configurations by at least 7% and up to 70%. Notably, in configurations 2 and 4, the addition of the reaction plate effectively reduced the motion of the spar, allowing more energy to be transferred and harvested by the float. As a result, these configurations exhibited a greater increase in power generation compared to both the standalone system and configurations 1 and 3. These findings demonstrate that integrating the float with the WT enhances the overall power generation of the hybrid system without requiring additional mooring modifications or cost increases, making it a more efficient alternative to standalone systems.

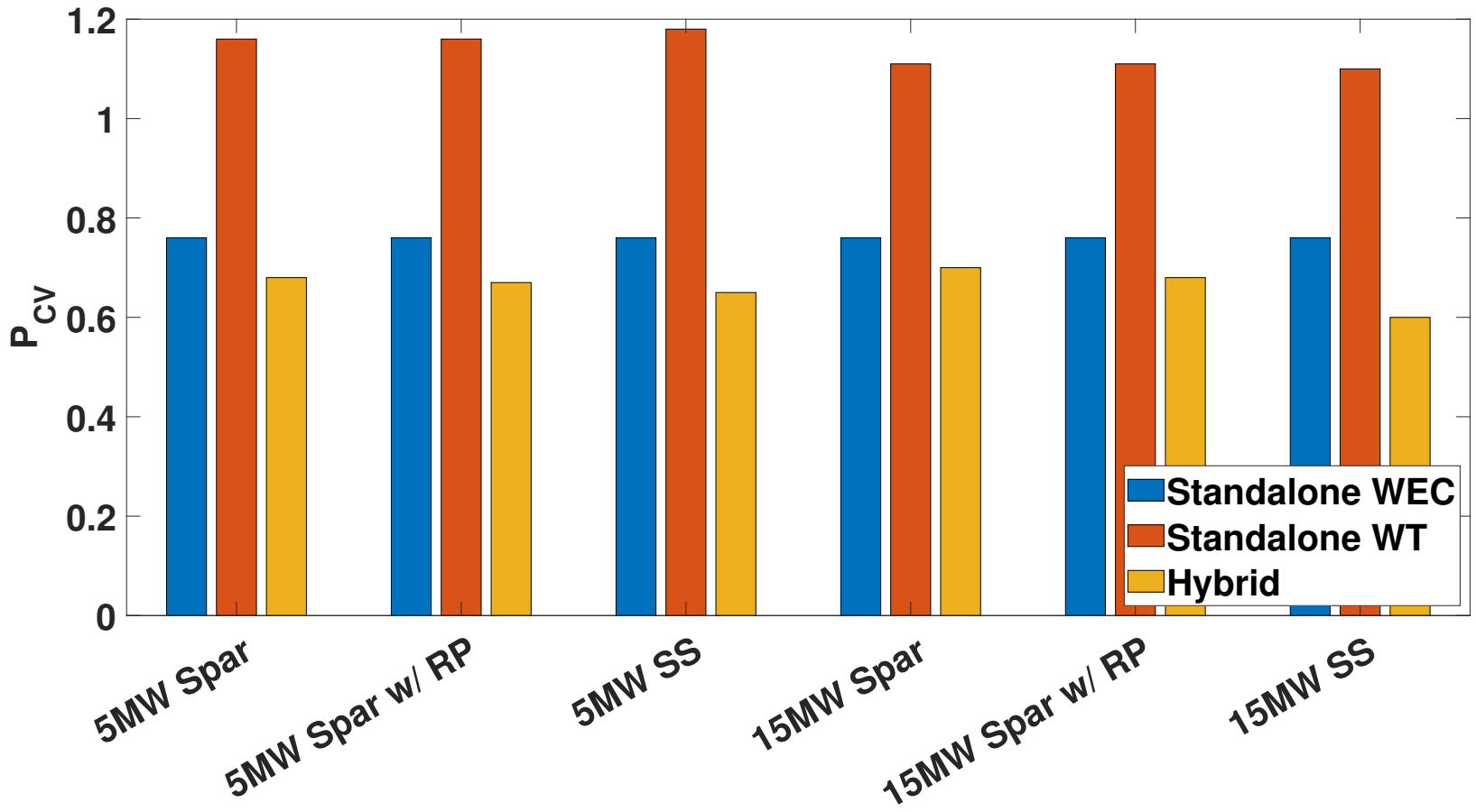

**Figure 33:** Power coefficient of variation for standalone and hybrid systems of different configurations

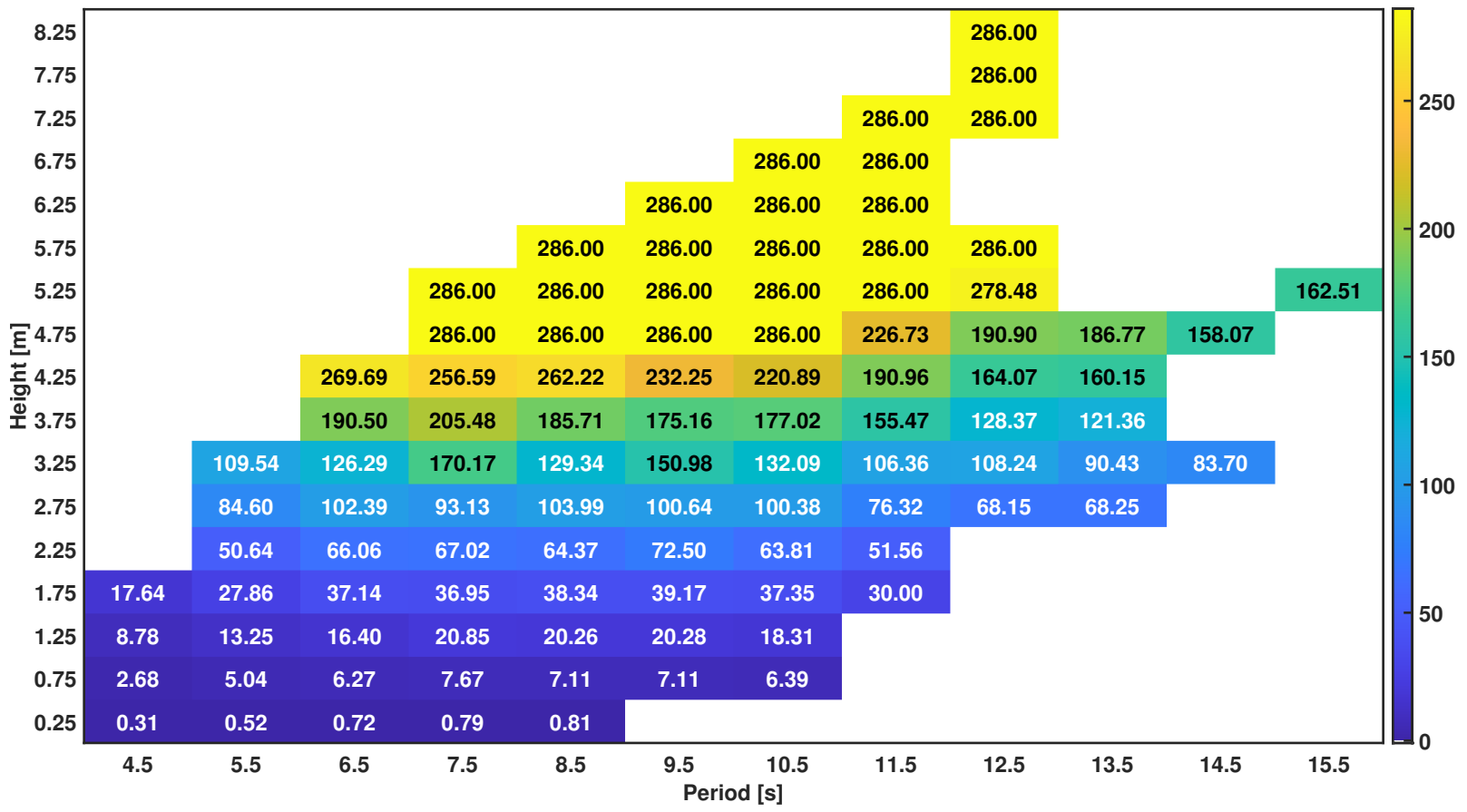

**Figure 34:** AEP of standalone RM3 WEC

### *6.2.2. Power coefficient of variation*

As much as maximizing power generation is important, ensuring the stability of that power with minimal variability is equally critical. Figure 33 presents the $P_{CV}$ for both standalone and hybrid systems with different configurations. The standalone WEC exhibits a lower $P_{CV}$ than that of both standalone FWTs by more than 30%, which proves that wave energy offers inherently greater stability and availability than wind energy. After integrating the WEC with the FWT, power variability is significantly reduced, with a reduction down to about 50% across all different configurations. This substantial reduction in power fluctuation demonstrates that hybrid systems can provide a more stable power supply all year round. Such improvements in power stability are considered a great advantage of hybrid systems over standalone systems. Ultimately, less power variation results in lower requirements for power and energy storage and eventually leads to lower LCOE.

### *6.2.3. AEP*

The AEP was calculated for all eleven configurations: a standalone RM3 WEC, standalone 5 MW WTs with spar and semi-submersible mooring foundations, standalone 15 MW WTs with spar and semi-submersible foundations, and three hybrid configurations for each WT size. All simulations were performed under irregular wave conditions to provide more accurate estimates of system performance in realistic ocean environments. Figures 34–42 present the electrical power matrices for each configuration after accounting for all losses, while Figures 43 and 44 present the corresponding AEP and CF values.

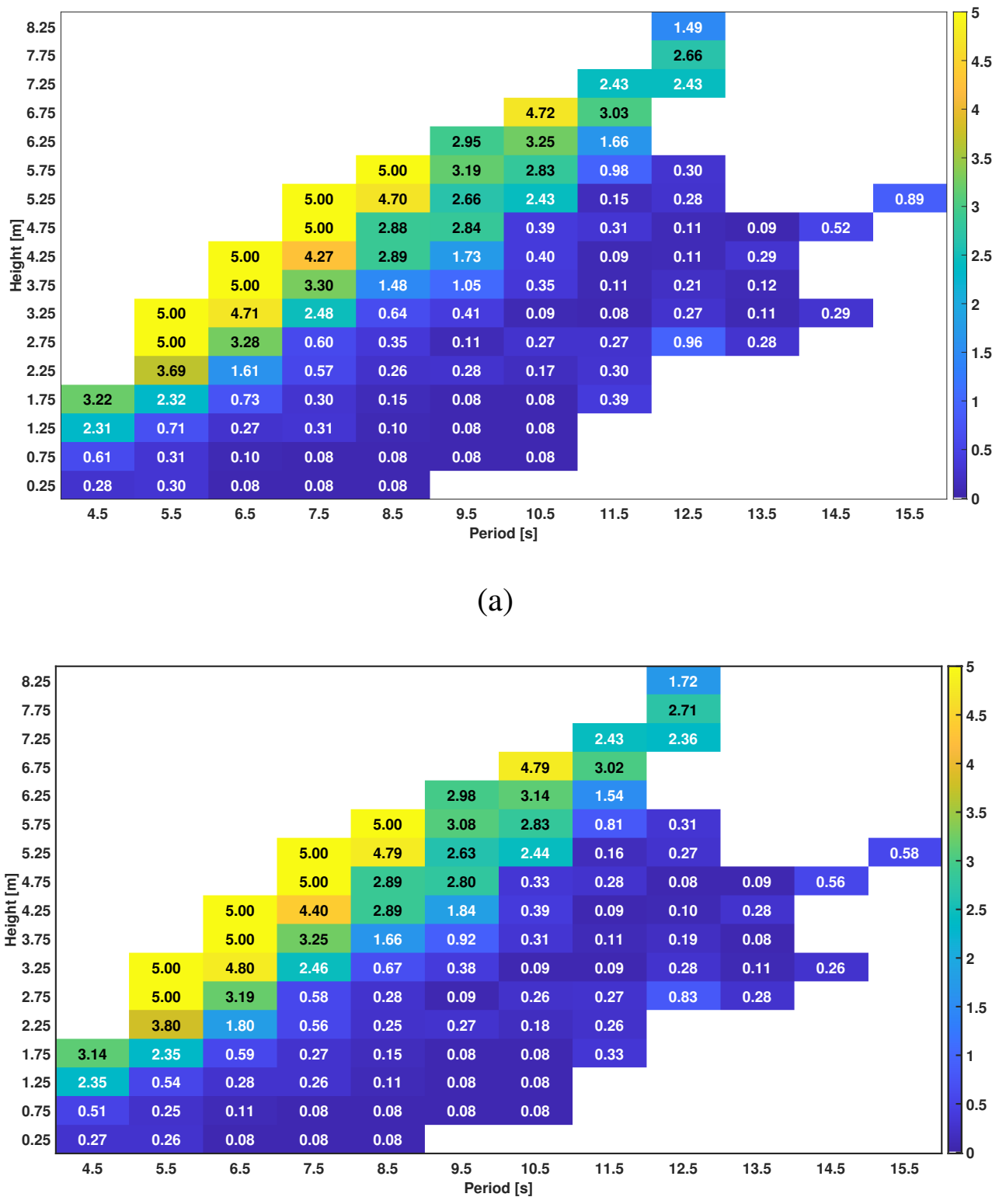

**Figure 35:** AEP of standalone 5 MW WT on a (a) spar platform and (b) semi-submersible platform at the Humboldt site

It is noticed that the AEP of the WEC is higher with the semi-submersible platform than the spar platform for both WT sizes. In contrast, the AEP of the WT remains largely unchanged across all standalone and hybrid configurations. Additionally, the results show a CF of approximately 25% for all configurations at the Humboldt Site, which is lower than typically desired for commercial viability. This low CF is primarily attributed to the relatively low average wind speeds at Humboldt, underscoring the importance of site selection for hybrid energy systems.

To examine the impact of site conditions on energy production, both standalone FWTs were tested at a second location, Kodiak. Assuming that the CF of the FWT in the hybrid system behaves similarly to that in the standalone system due to identical wind conditions, simulations were conducted to calculate only the AEP for only standalone configurations at Kodiak. Figures 45-46 present the electrical power matrices for each configuration at the Kodiak site after accounting for all losses, while Figure 47 demonstrates the AEP and CF values at the same site. The CF is approximately 37% for the 5 MW WT and 43% for the 15 MW WT. The results show a significant increase in CF, about 85% higher for both WTS, thus confirming the critical role of wind resource quality in maximizing energy production and highlighting the importance of strategic site selection for hybrid renewable energy farms.

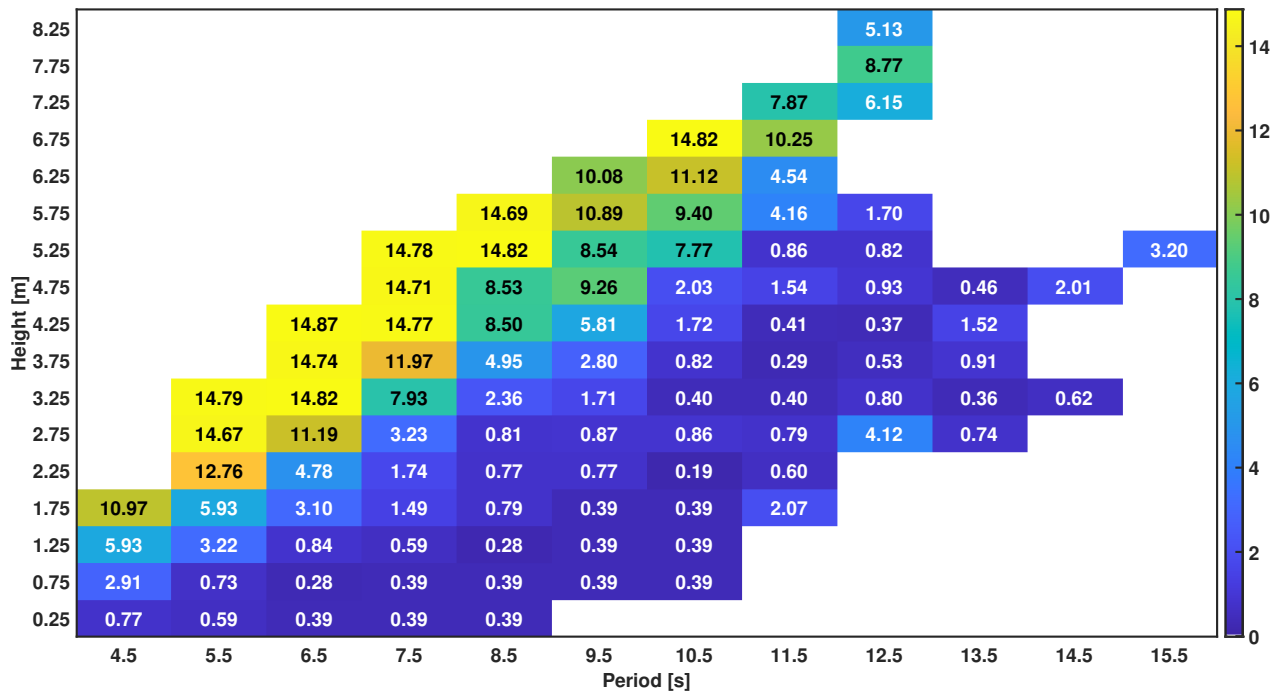

(a)

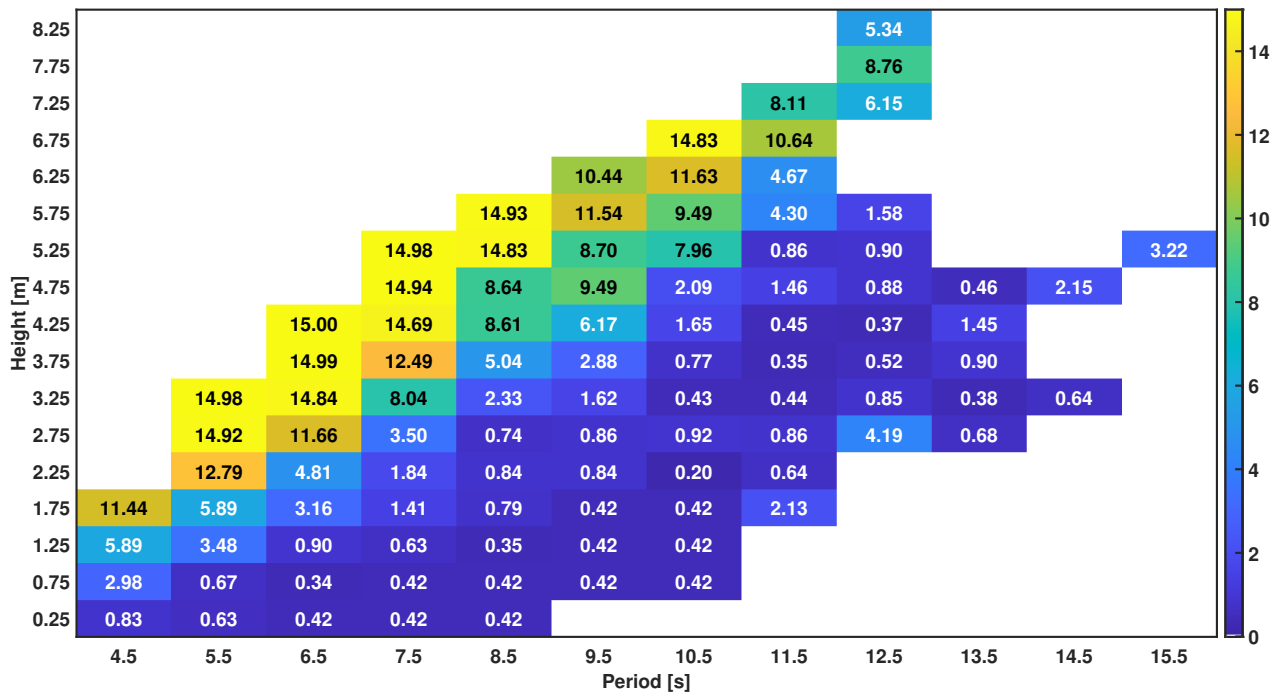

(b)

**Figure 36:** AEP of standalone 15 MW WT on a (a) spar platform and (b) semi-submersible platform at the Humboldt site

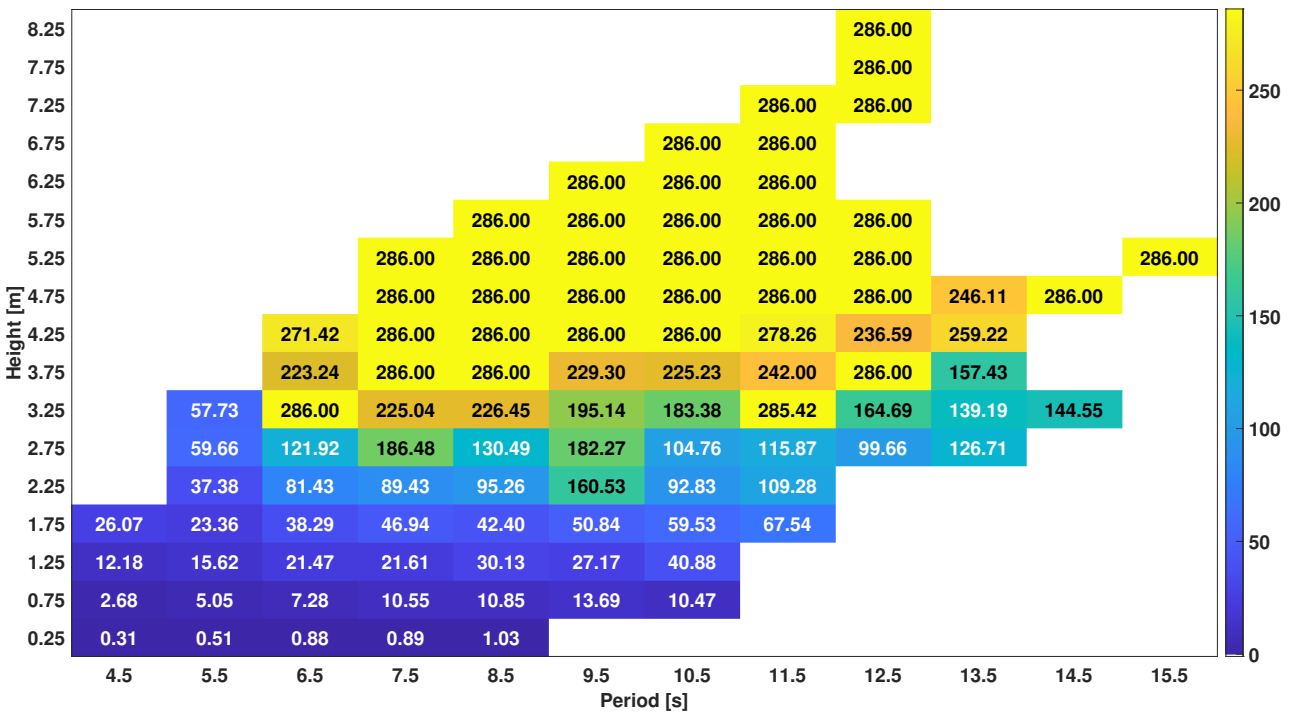

(a)

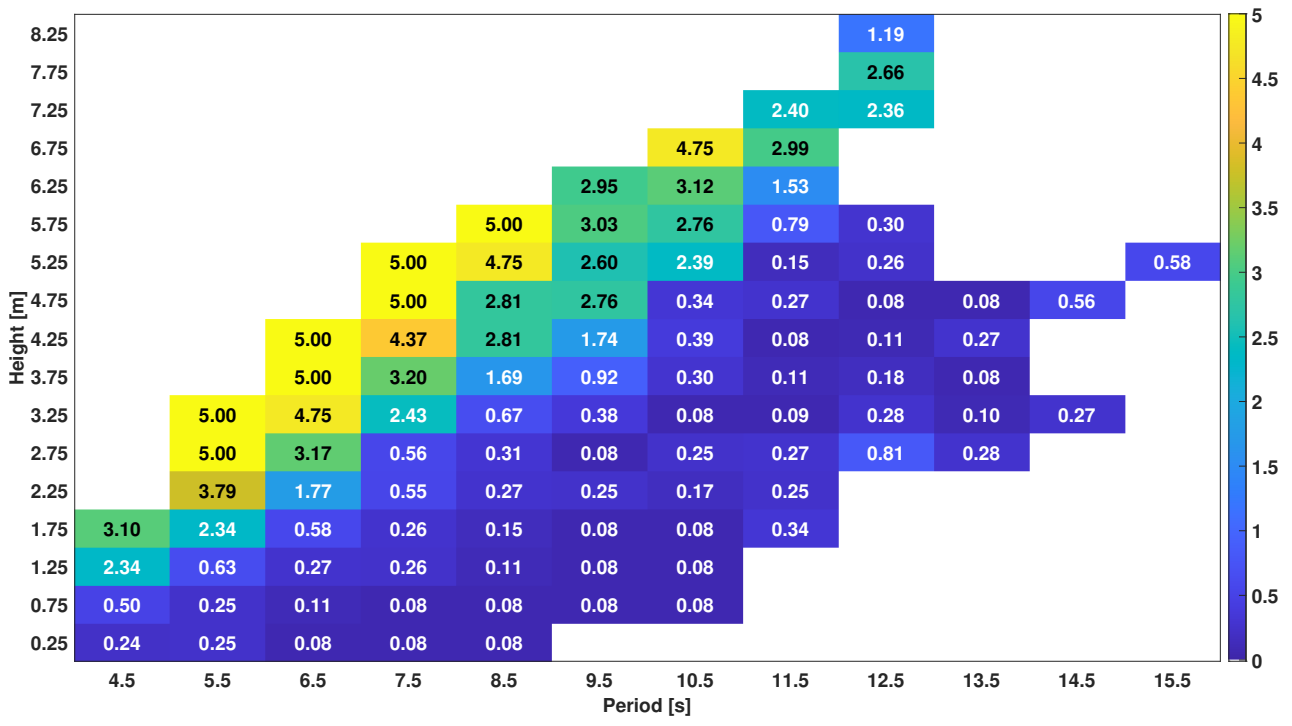

(b)

**Figure 37:** AEP of (a) WEC (b) WT in a hybrid system with 5 MW WT on a spar platform at the Humboldt site

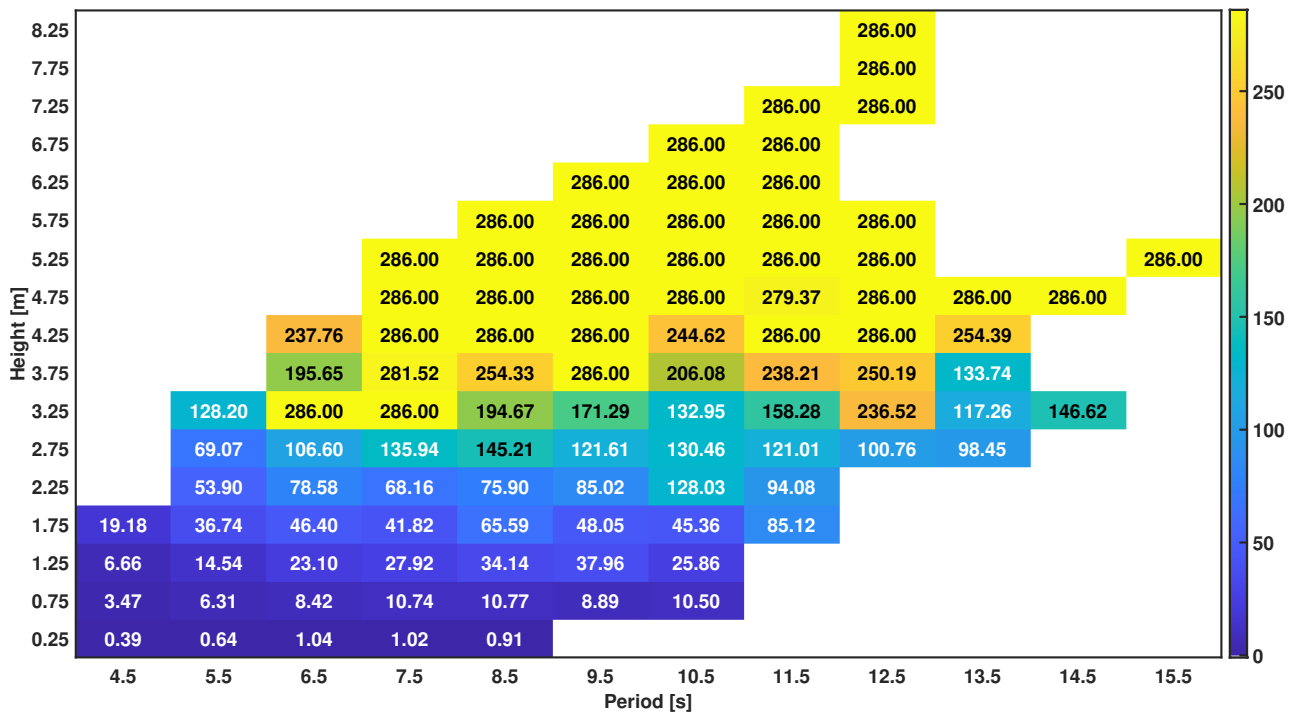

(a)

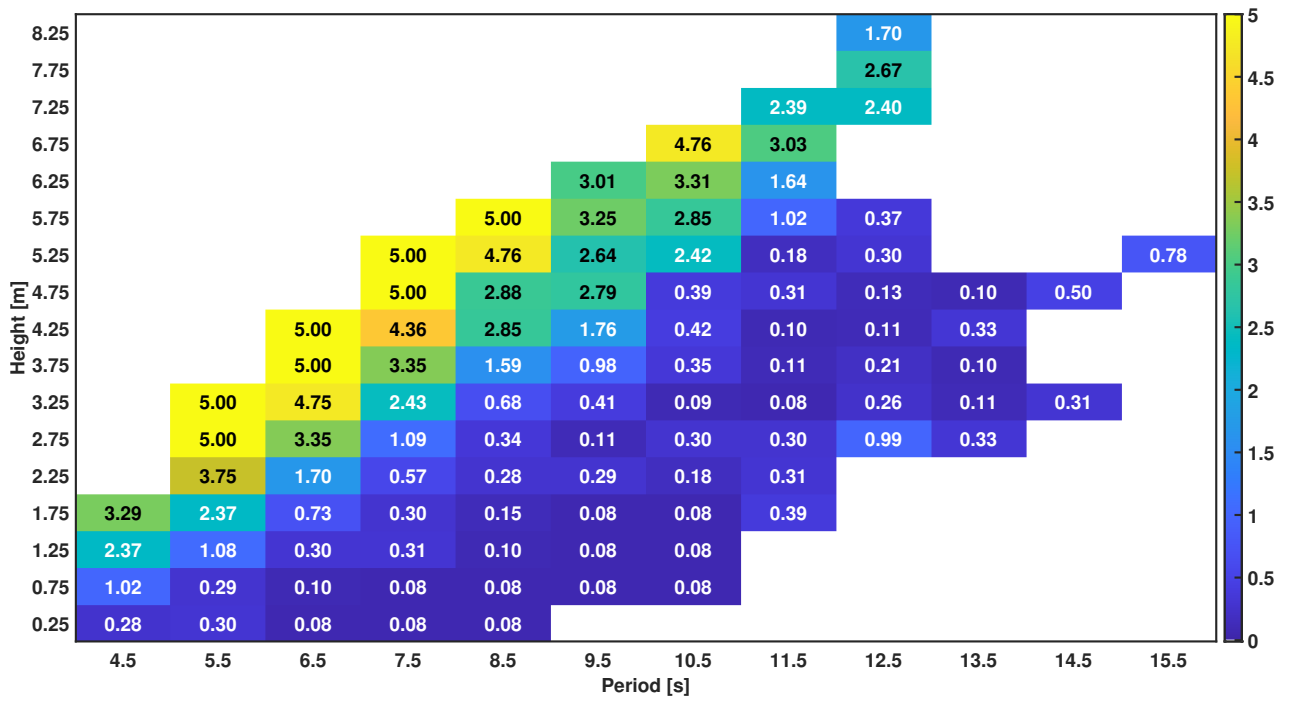

(b)

**Figure 38:** AEP of (a) WEC (b) WT in a hybrid system with 5 MW WT on a spar platform with a reaction plate at its bottom at the Humboldt site

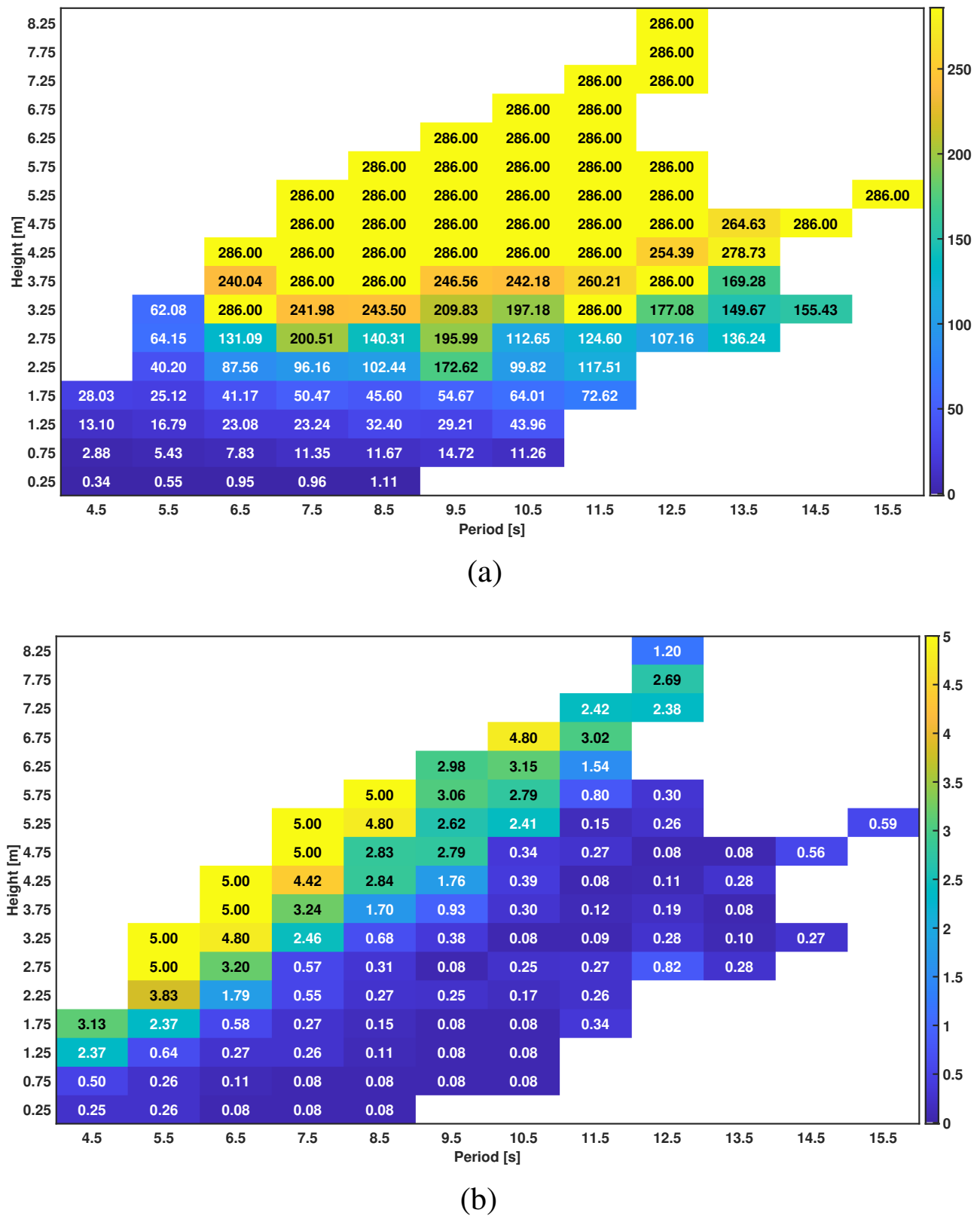

**Figure 39:** AEP of (a) WEC (b) WT in a hybrid system with 5 MW WT on a semi-submersible platform at the Humboldt site

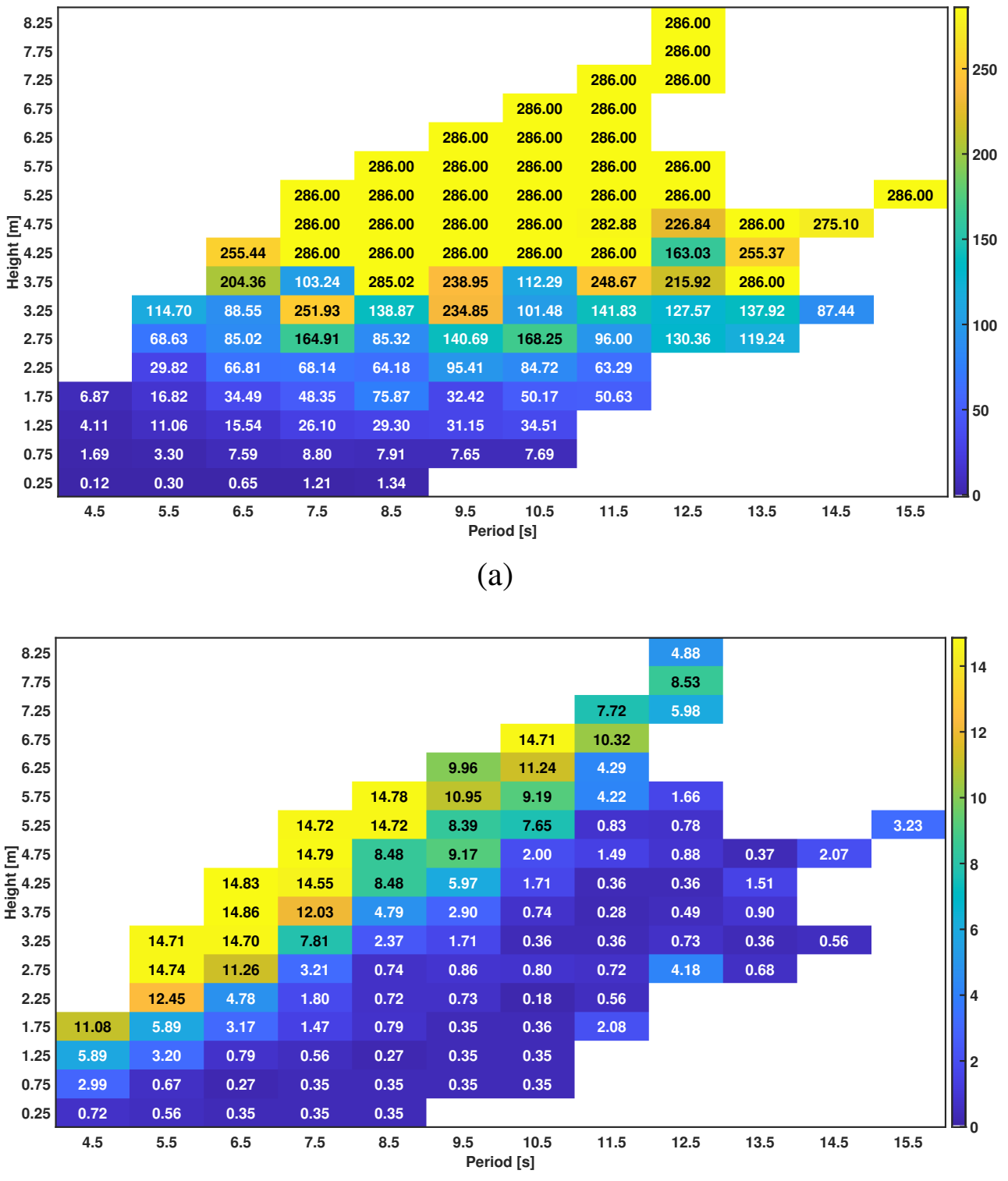

(a)

(b)

**Figure 40:** AEP of (a) WEC (b) WT in a hybrid system with 15 MW WT on a spar platform at the Humboldt site

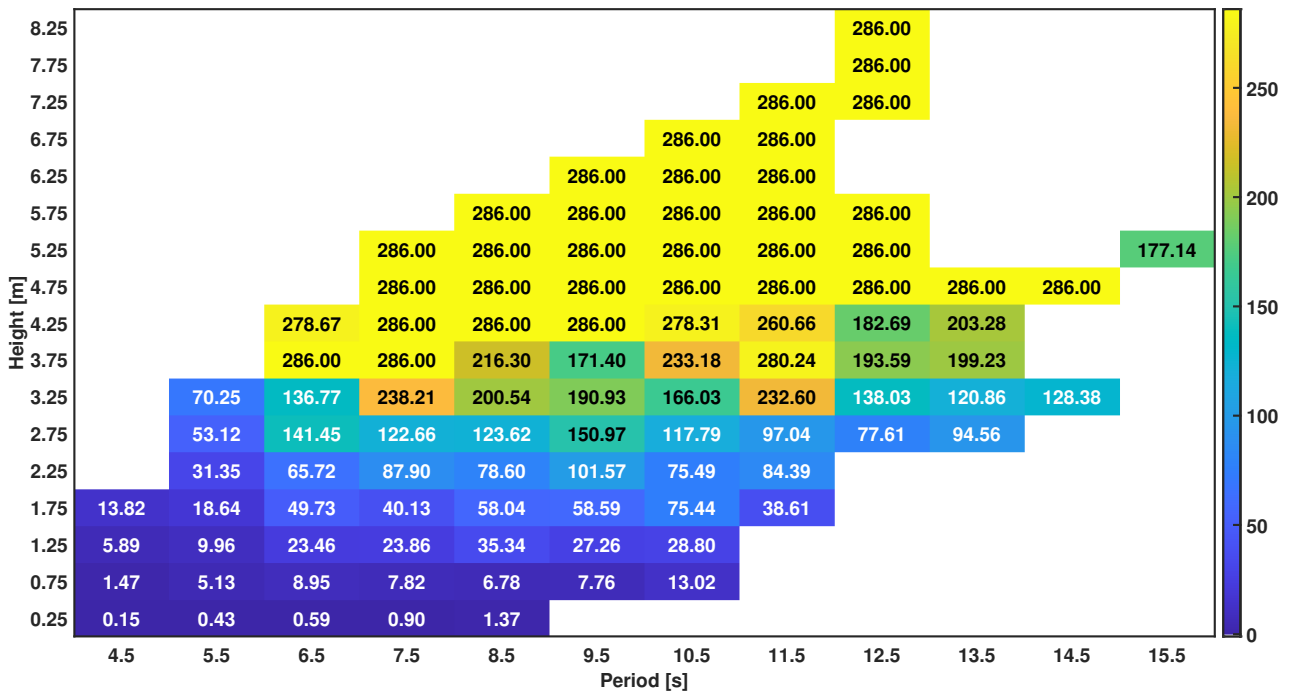

(a)

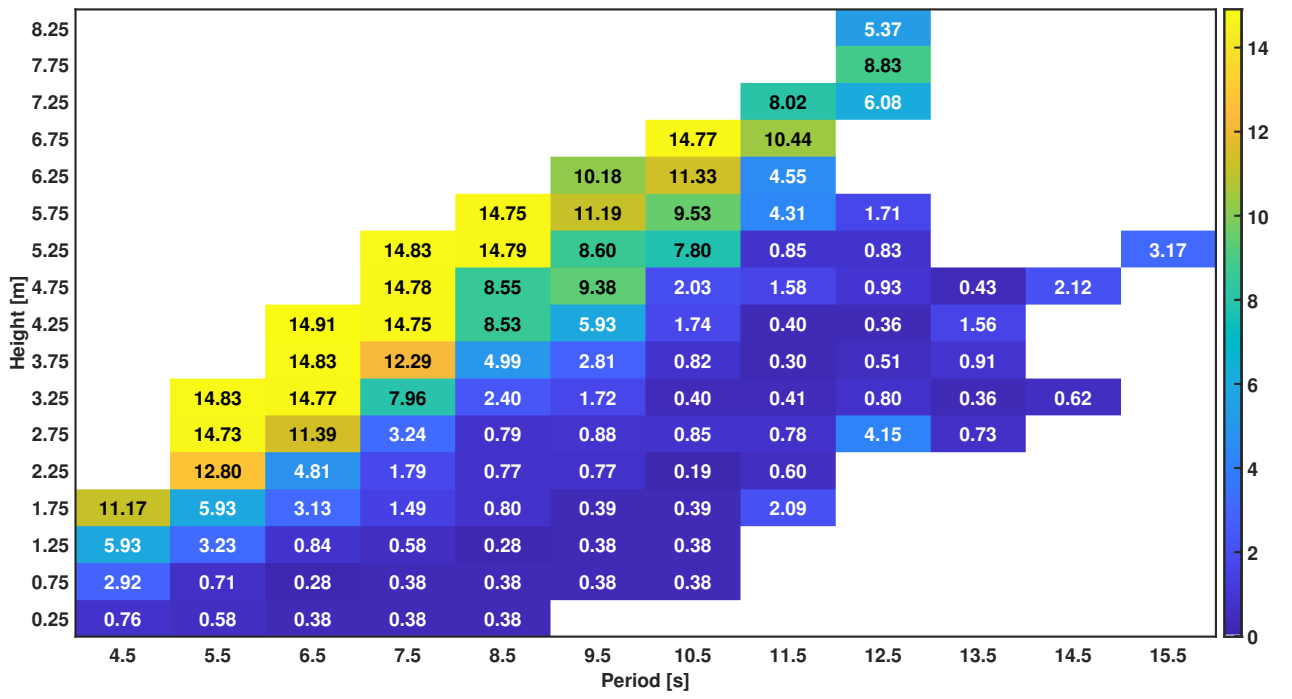

(b)

**Figure 41:** AEP of (a) WEC (b) WT in a hybrid system with 15 MW WT on a spar platform with a reaction plate at its bottom at the Humboldt site

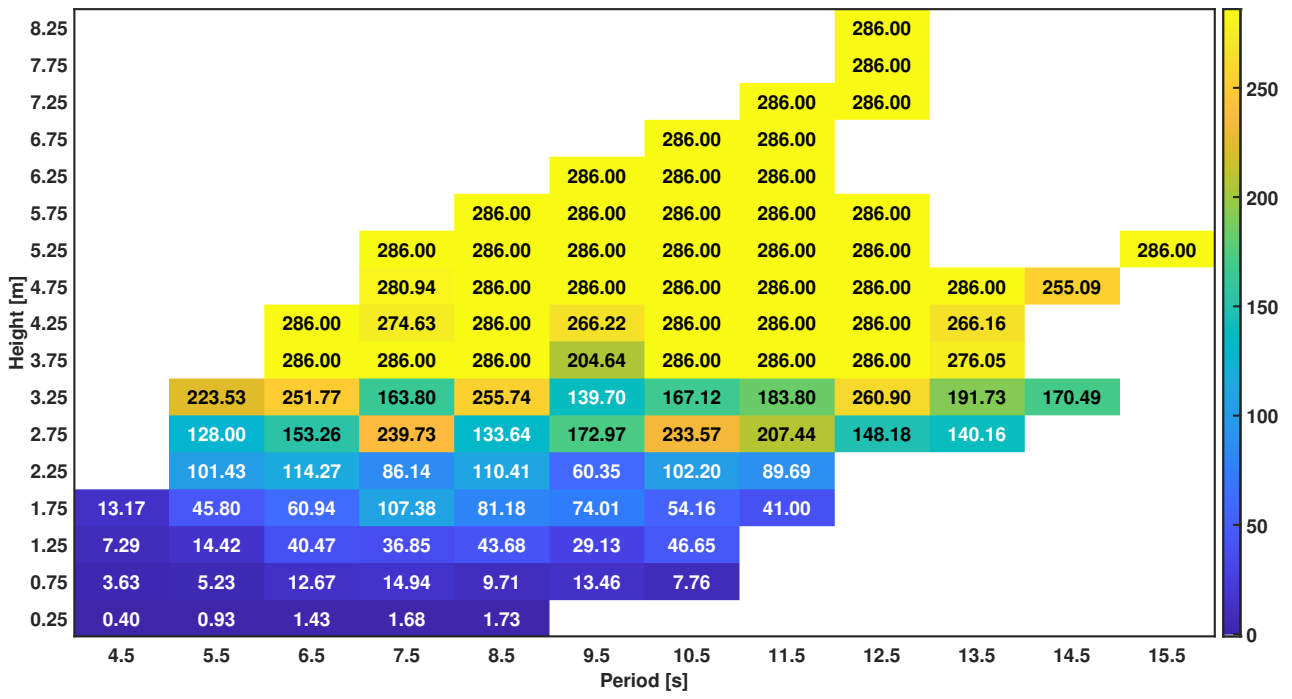

(a)

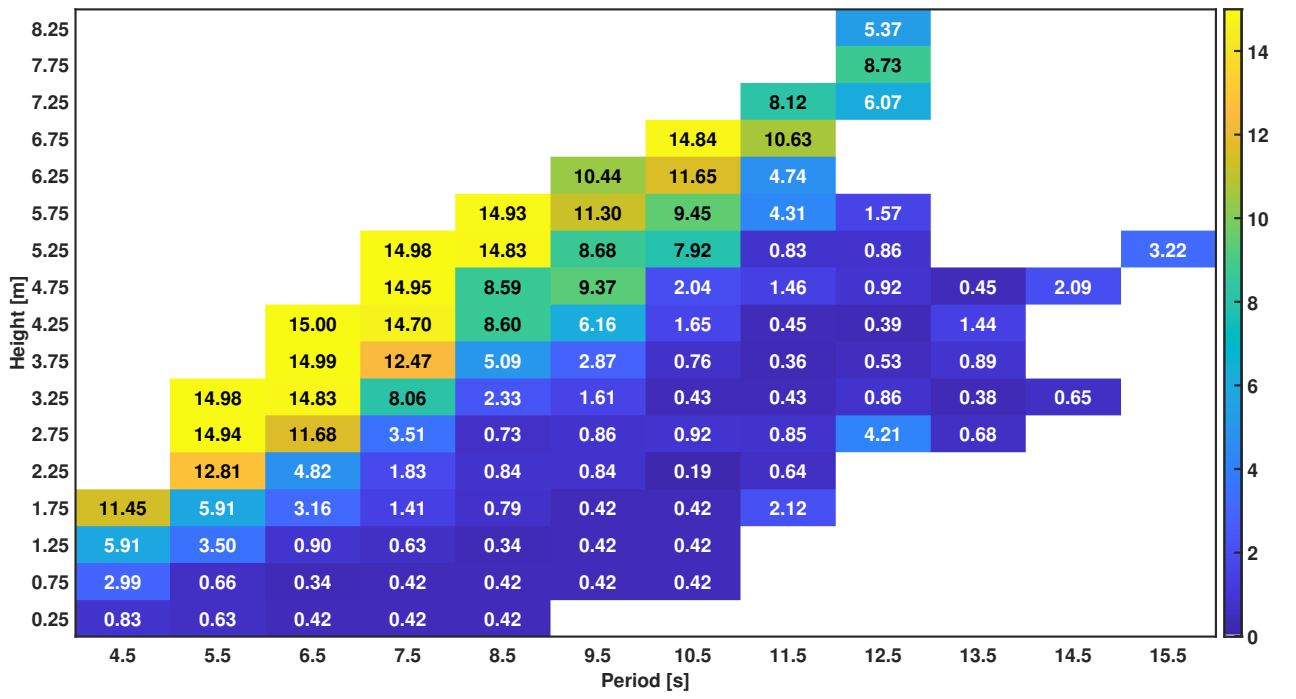

(b)

**Figure 42:** AEP of (a) WEC (b) WT in a hybrid system with 15 MW WT on a semi-submersible platform at the Humboldt site

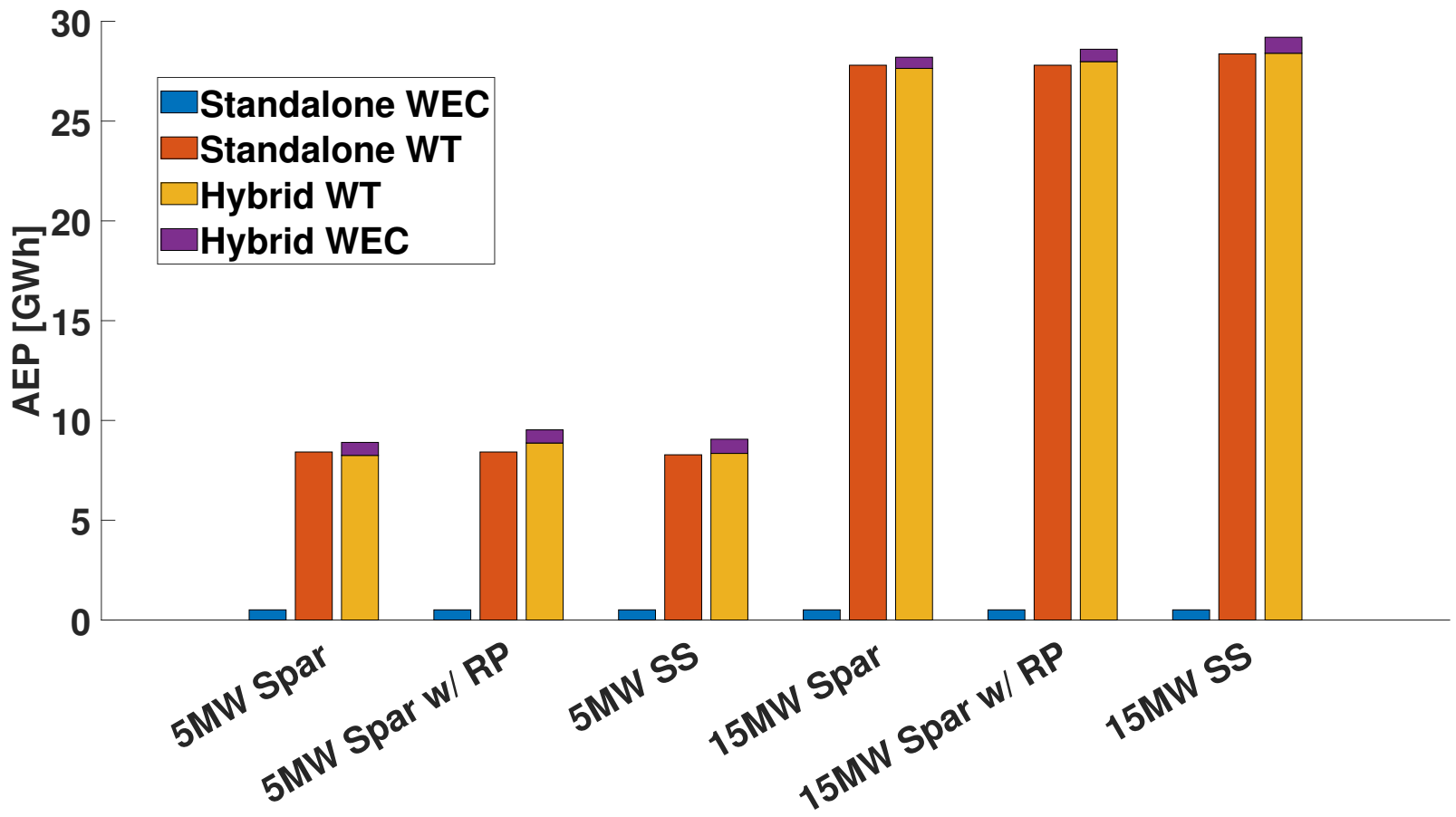

**Figure 43:** AEP of all configurations at the Humboldt Site

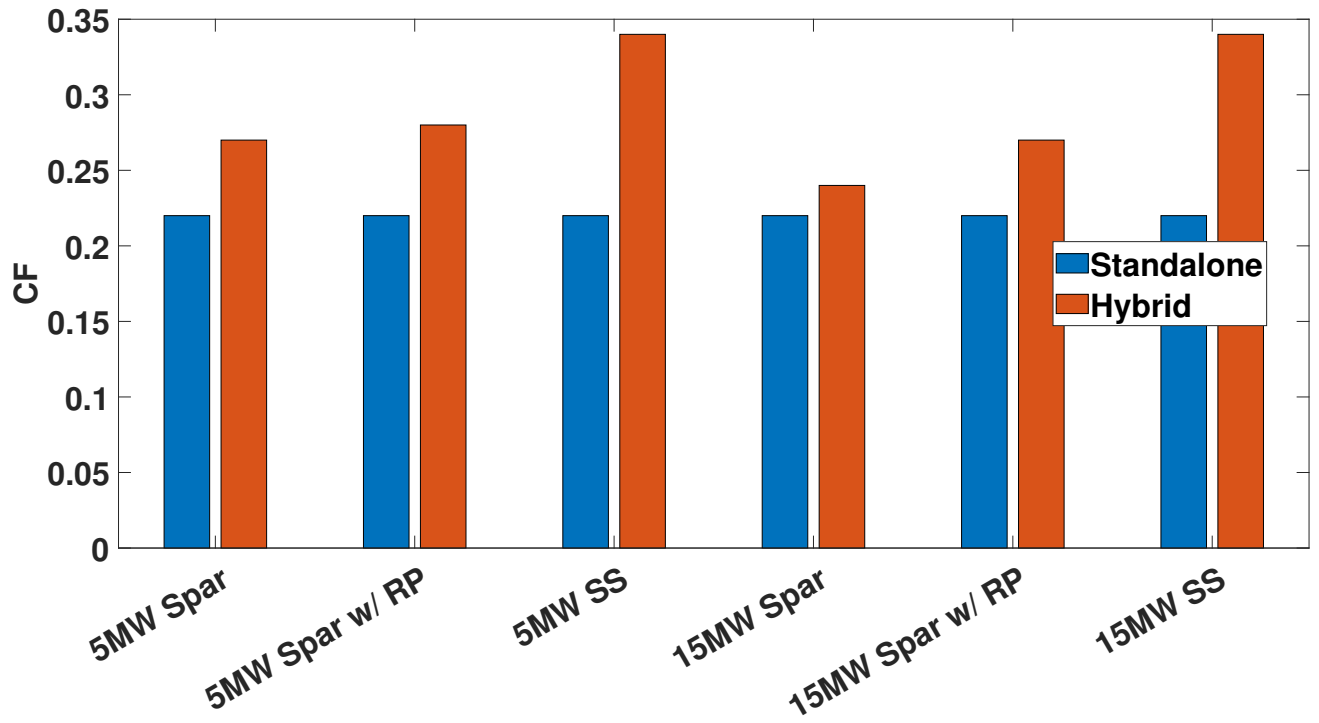

(a)

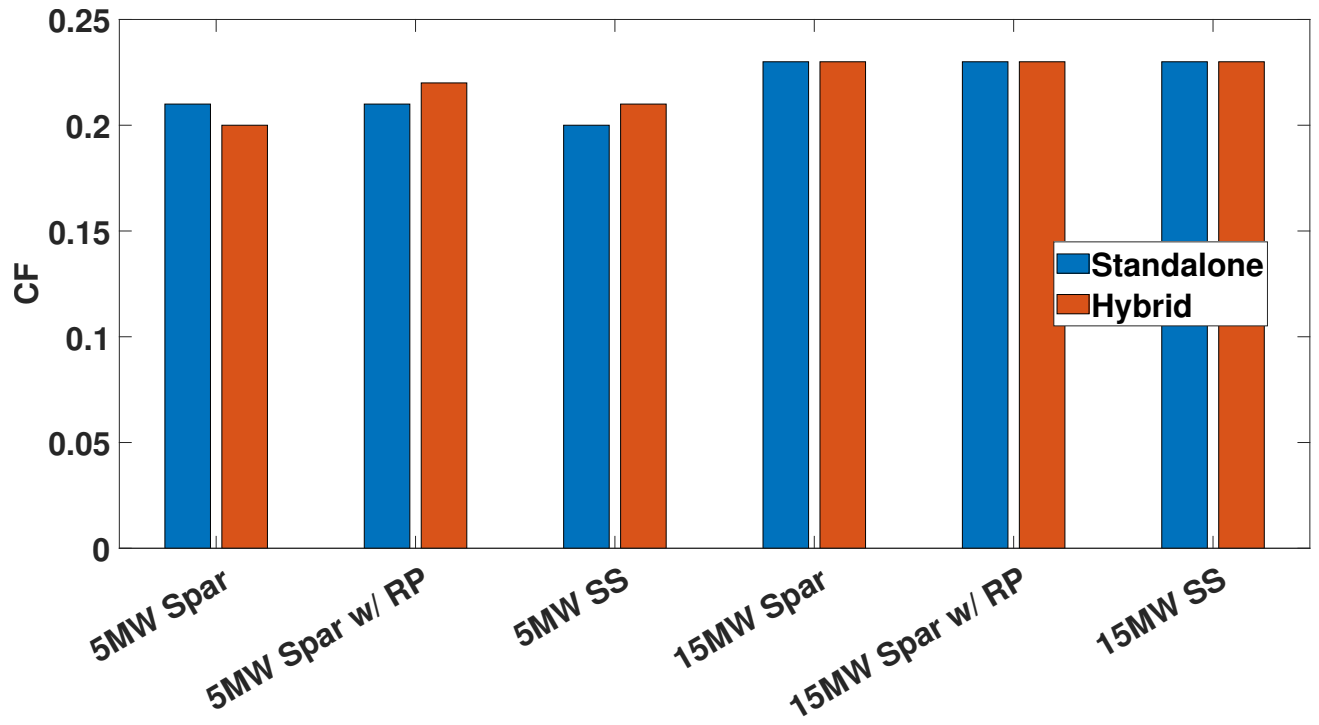

(b)

**Figure 44:** CF of (a) WEC (b) WT in all configurations at the Humboldt Site

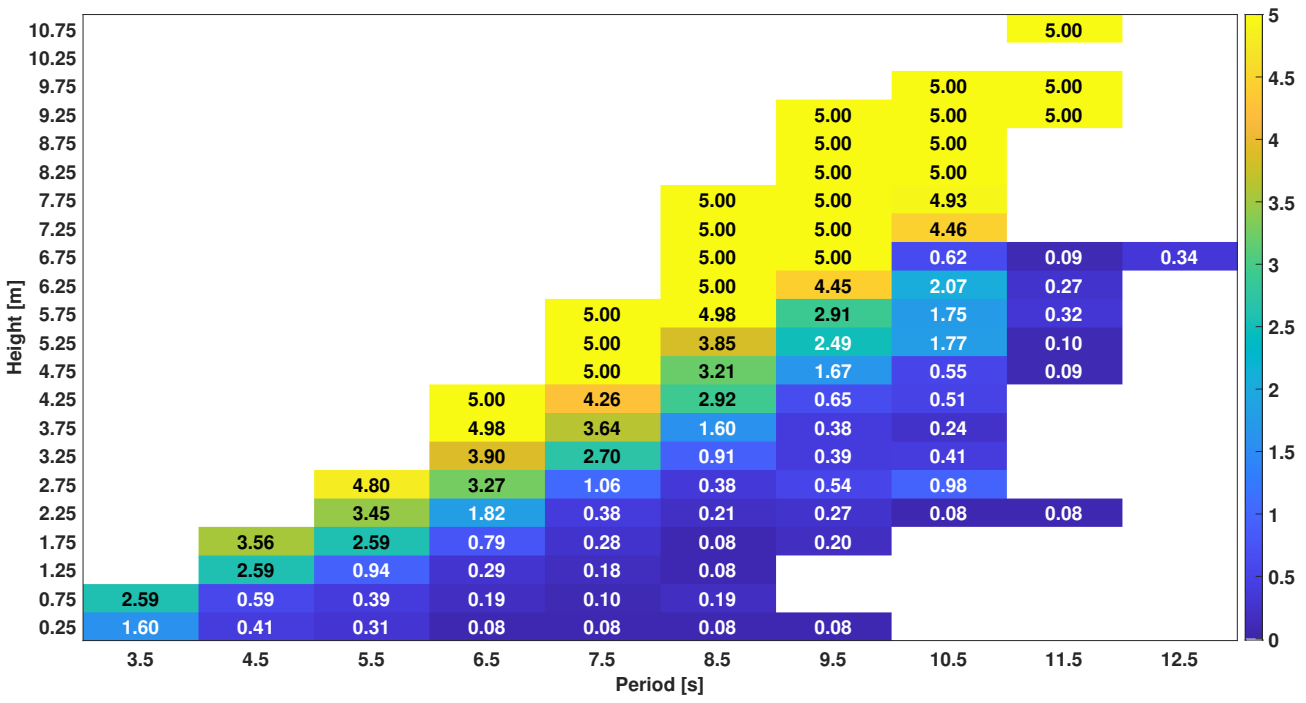

(a)

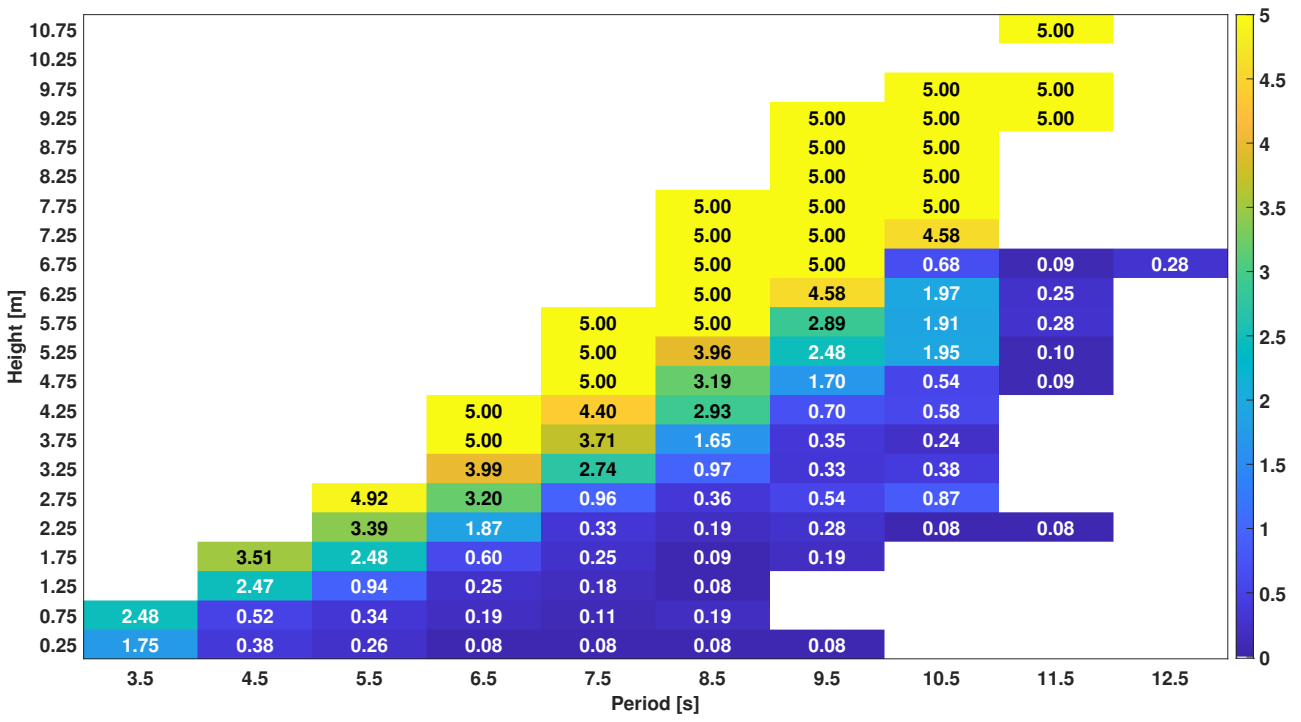

(b)

**Figure 45:** AEP of standalone 5 MW WT on a (a) spar platform and (b) semi-submersible platform at the Kodiak site

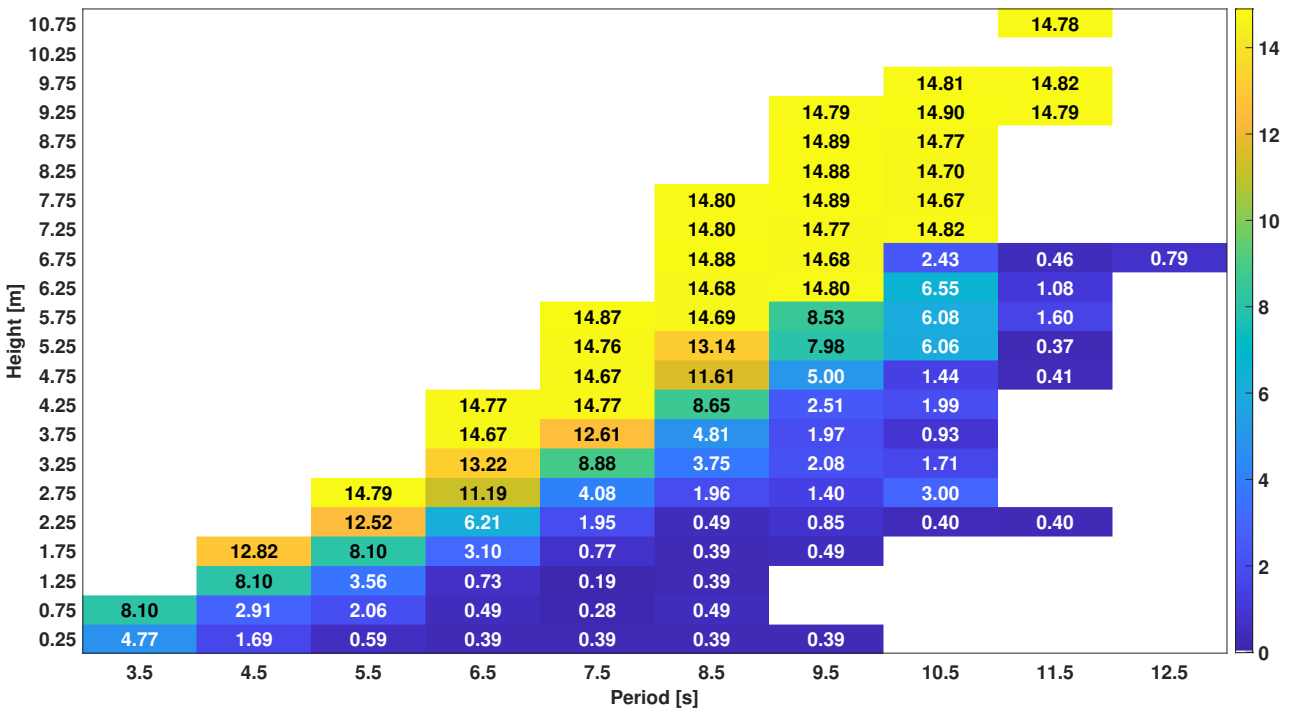

(a)

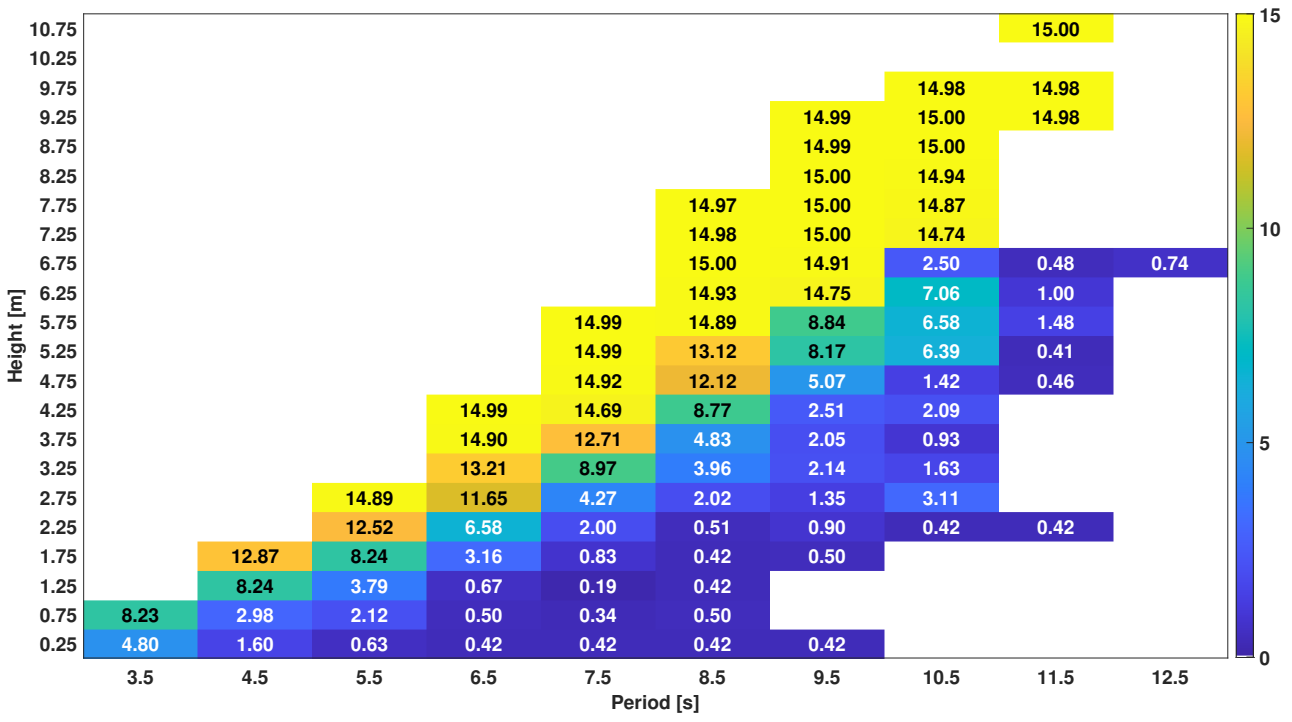

(b)

**Figure 46:** AEP of standalone 15 MW WT on a (a) spar platform and (b) semi-submersible platform at the Kodiak site

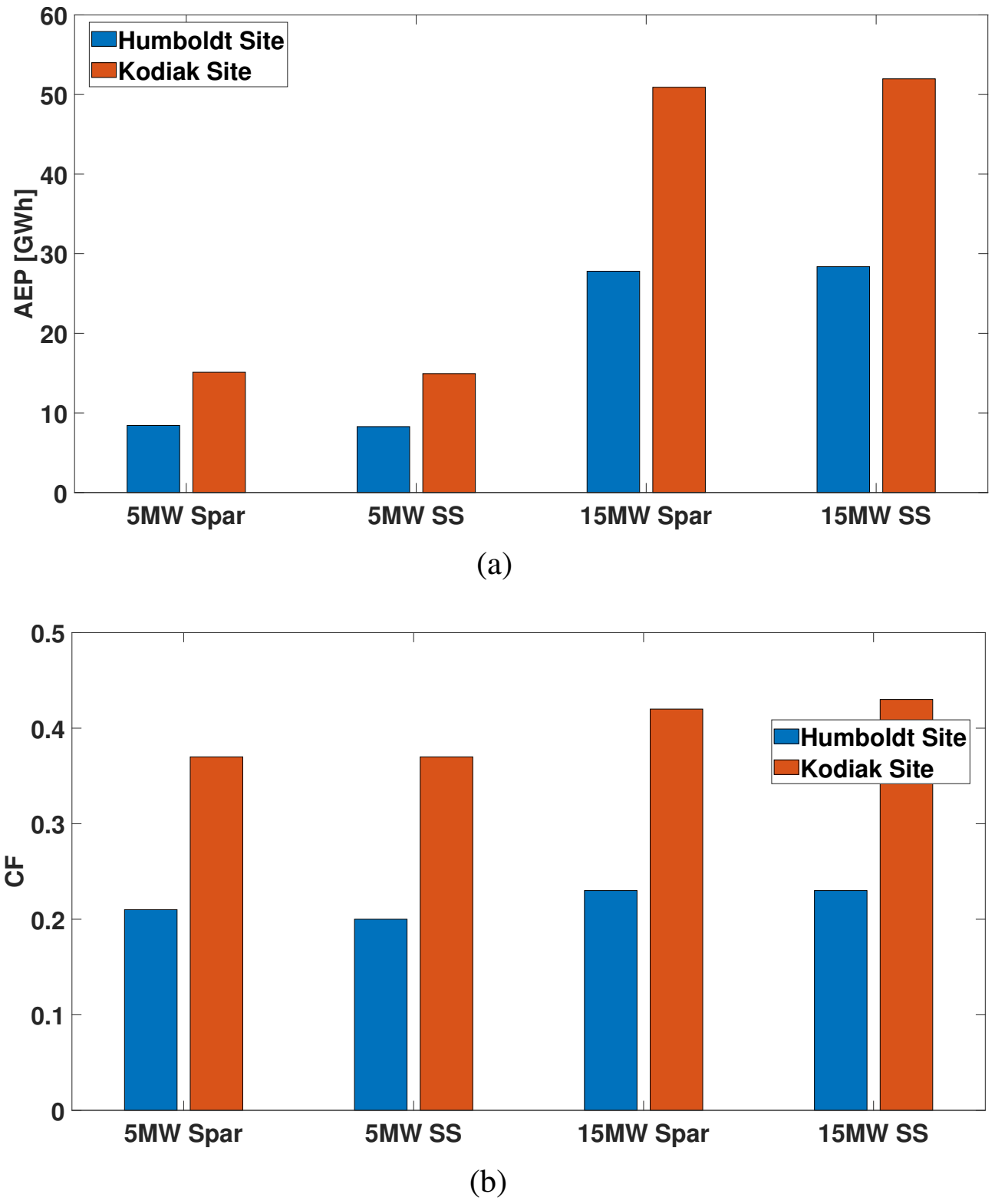

**Figure 47:** (a) AEP and (b) CF of standalone configurations at Kodiak Site

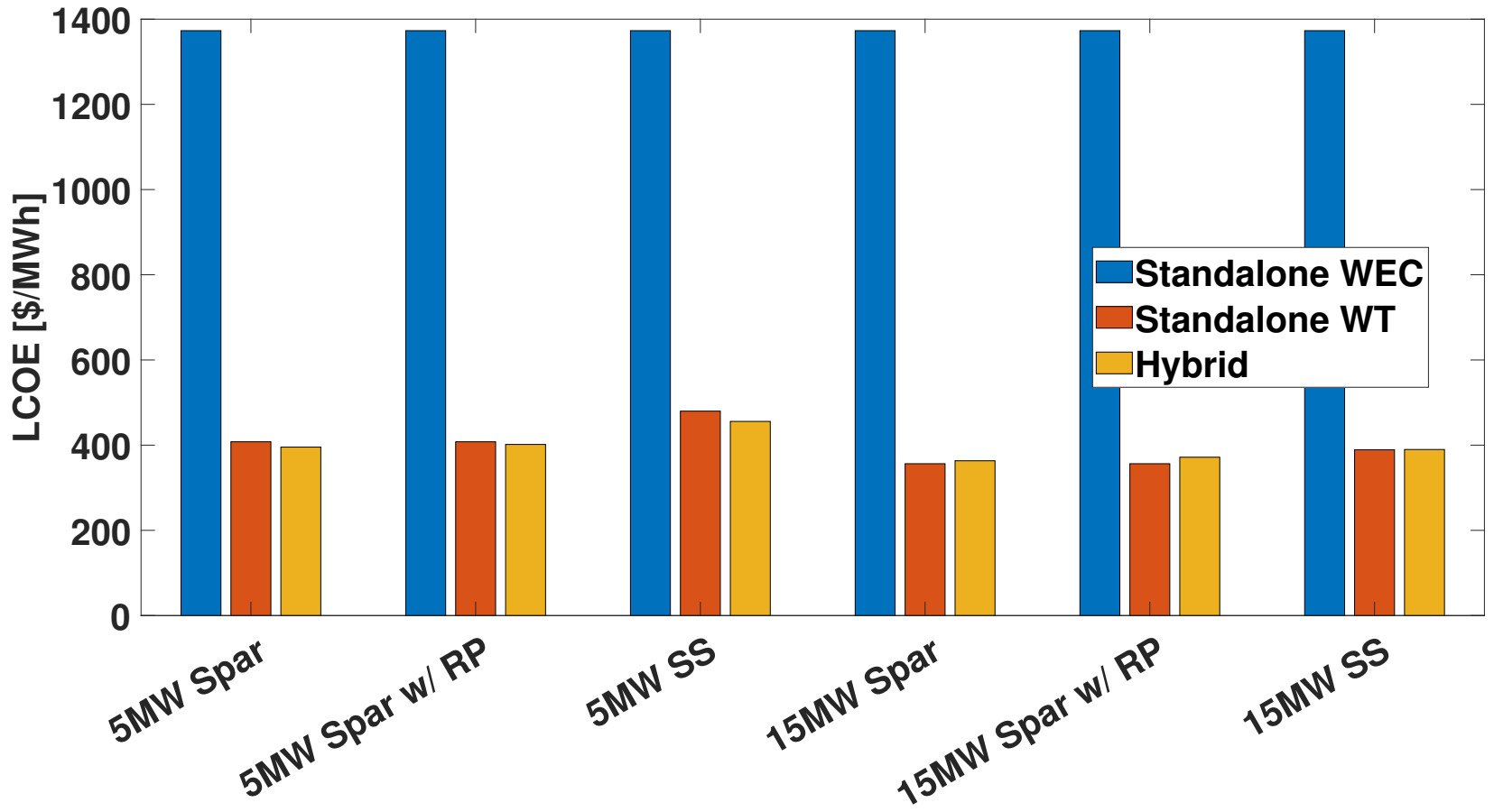

**Figure 48:** LCOE of a single unit within a 100-unit farm of different configurations

## 6.3. Economics

The economic performance of the hybrid and standalone systems is assessed using the LCOE metric. This analysis enables a comparison of cost-effectiveness between standalone WECs, standalone WTs, and hybrid configurations, highlighting the potential economic benefits of integrating wave and wind energy technologies. Importantly, these results represent one of the first techno-economic assessments for such hybrid offshore configurations, providing novel insights into the cost implications and feasibility of hybrid wind-wave energy systems.

### *6.3.1. LCOE*

The LCOE of each farm configuration was calculated based on Equation 3. The cost parameters for each configuration were drawn from publicly available data sources Neary et al. (2014); Beiter et al. (2016); Myhr et al. (2014); Martinez and Iglesias (2022) and adjusted to 2025 USD. The FCR was assumed to be 11% to match the rate in Neary et al. (2014); Beiter et al. (2016); Myhr et al. (2014); Martinez and Iglesias (2022). Figure 48 presents the LCOE for a single unit within the 100-unit farms.

A significant difference can be observed between the LCOE of standalone WEC systems and standalone WT systems. The WEC farm exhibits higher costs combined with lower energy production compared to the WT farm. Therefore, integrating WECs with WTs into a hybrid system can help reduce overall costs, particularly for the WEC component. By comparing the LCOE of the hybrid systems to that of the standalone systems, it is observed that, for the 5 MW WT configurations, the hybrid systems achieve a slight reduction in LCOE across all mooring foundations: approximately a 3% reduction for the spar platform and a 5% reduction for the semi-submersible platform. This improvement indicates that the increased energy production from the hybrid system effectively offsets the additional costs.

In contrast, for the 15 MW WT configurations, the hybrid system exhibits a slight LCOE increase: approximately 2% and 4% for the spar platform with and without the reaction plate, respectively, while the LCOE remains nearly unchanged for the semi-submersible platform. Although the hybrid system with the 15 MW WT increased the total AEP, it was insufficient to fully compensate for the additional costs.

These results suggest that the impact of hybridization on LCOE is sensitive to the scale of the WT. Specifically, smaller-scale turbines (5 MW) benefit more from hybridization in terms of LCOE reduction, while larger-scale turbines (15 MW) do not. This highlights the importance of carefully balancing system size and associated costs. It also suggests that further improvements could be achieved through concurrent co-design processes that optimize the sizing of both WTs and WECs based on location-specific and scale-specific conditions.

## 6.4. Synergy

The results in Figure 49 illustrate the symbiotic relationship using LCOE between the FWT and the WEC in a hybrid system considering a 100-unit hybrid system farm by integrating the RM3 with the 5 MW and 15 MW FWTs on a spar mooring platform with and without the reaction plate at the bottom and a semi-submersible mooring platform.

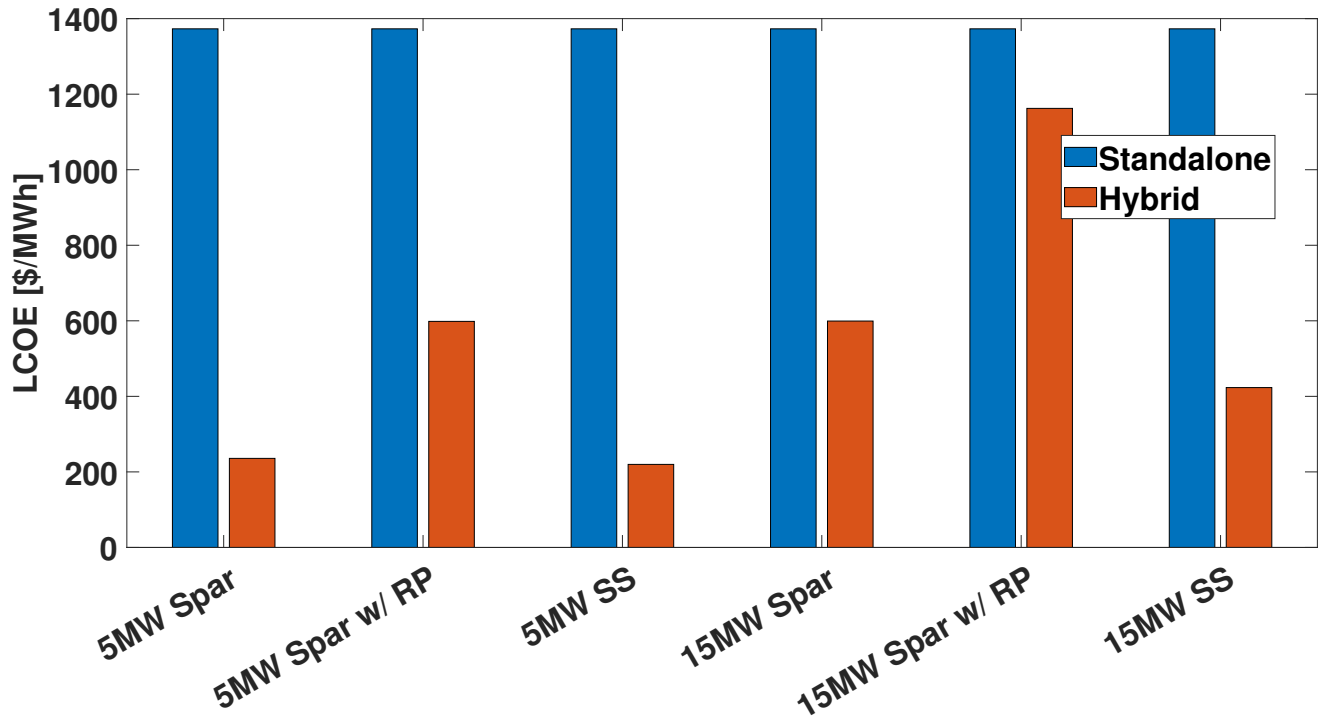

(a)

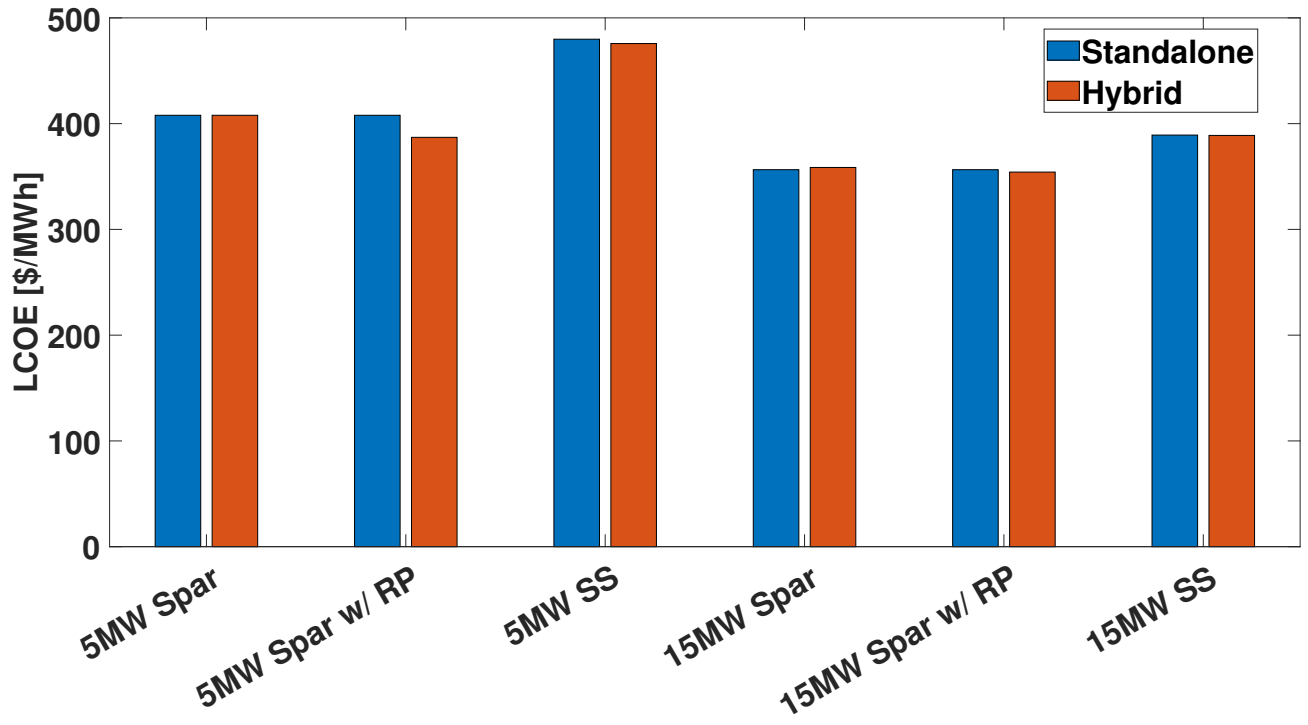

(b)

**Figure 49:** Comparison between the LCOE of (a) WEC and (b) WT in a 100 standalone RM3 farm and 100 standalone 5 MW and 15 MW WTs farms with a 100 hybrid system farm by integrating the RM3 with the 5 MW and 15 MW FWTs

These hybrid configurations are compared against standalone farms consisting of 100 units of RM3, 5 MW FWTs, and 15 MW FWTs.

Across all hybrid systems, the LCOE of the WEC decreased significantly from $1373/MWh to as low as $235.79/MWh, representing reductions ranging from 15% to 83%. This substantial improvement is expected, as key capital expenses (e.g., mooring foundations, electrical cables) and operational expenses are shared with the FWT, resulting in notable cost savings. For the 5 MW FWT, the LCOE decreased from $407.90/MWh to $387.04/MWh (approximately a 5% reduction) with the spar platform with the plate at the bottom, and from $479.84/MWh to $475.73/MWh (less than a 1% reduction) with the semi-submersible platform. However, no change was observed with the spar platform without the bottom plate. This behavior suggests a mutualistic symbiotic relationship between the WEC and the 5 MW FWT in most configurations, except for the spar platform without the plate, which exhibited a commensalistic relationship. Similarly, for the 15 MW FWT, the LCOE showed a slight reduction (less than 1%) in the hybrid configurations with the spar platform with a plate and the semi-submersible platform. In contrast, a slight LCOE increase (less than 1%) was observed for the hybrid system with the spar platform without a plate, indicating a parasitic relationship in that specific case. Overall, the synergy analysis reveals a relationship that predominantly resembles mutualism, where the WEC benefits significantly, while the FWT experiences slight improvements or negligible impacts. These results suggest that with further system-level optimization, full mutualism could potentially be achieved.

While LCOE is a valuable metric for assessing the economic feasibility of WEC integration, it does not capture additional values and benefits that the WEC may offer to the FWT. For instance, the integration could reduce fatigue loads on the FWT Kluger (2017) or decrease fatigue damage in mooring cables Clark and Paredes (2018). These effects could lower material costs and extend the operational lifespan of the FWT, further enhancing the overall techno-economic performance of hybrid systems. Such benefits, although not considered in this study, represent important areas for future research.

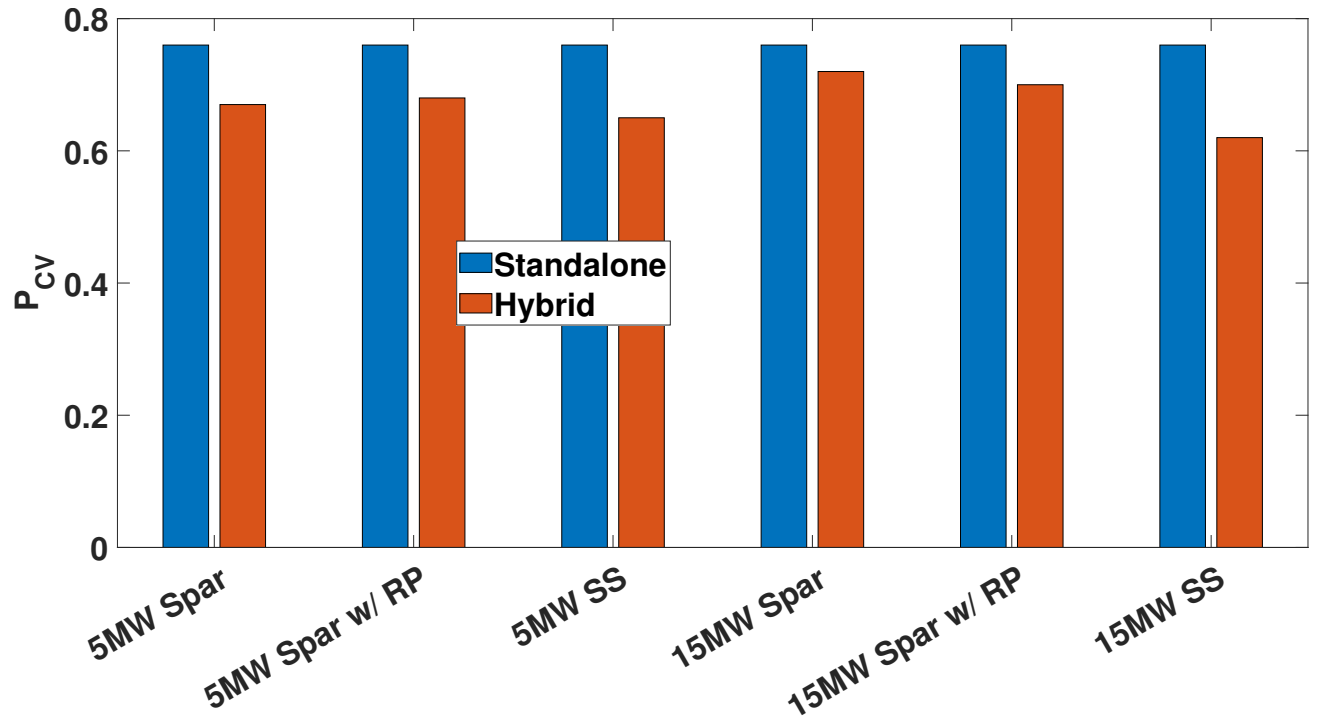

(a)

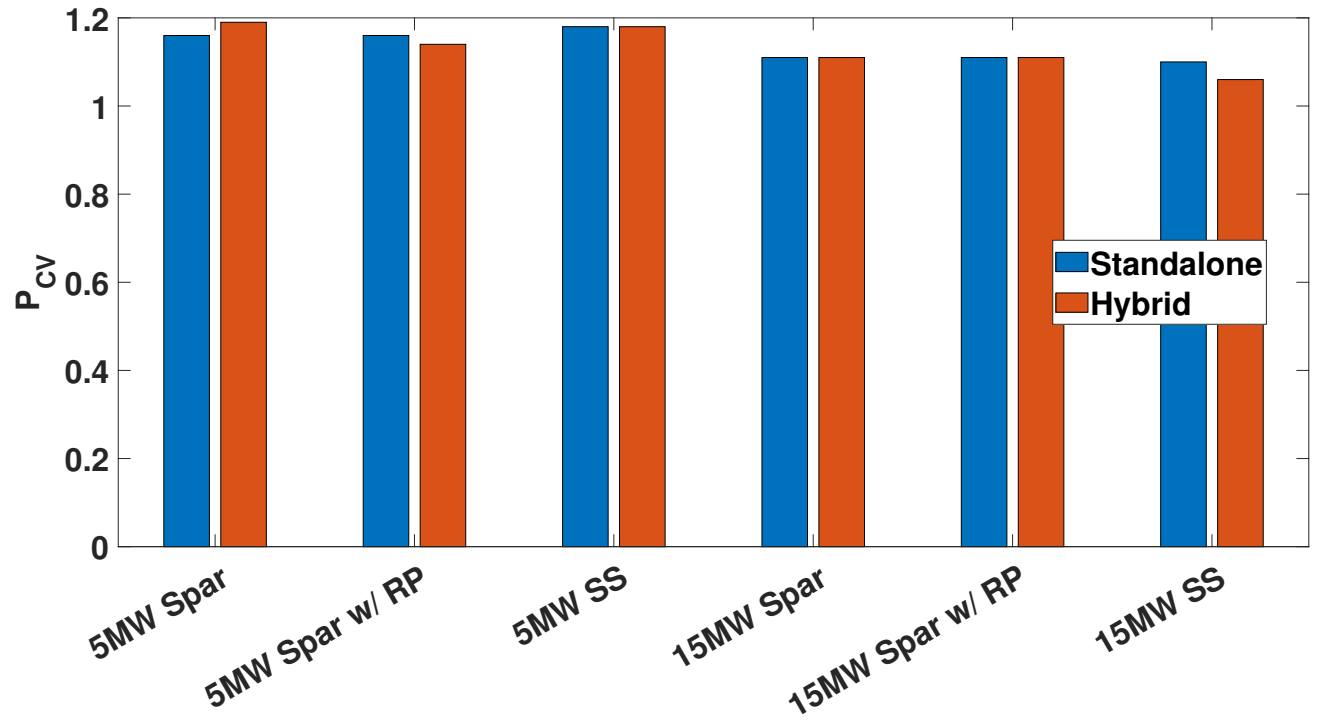

(b)

**Figure 50:** Comparison between the $P_{CV}$ of (a) WEC and (b) WT for RM3, standalone 5 MW, and 15 MW WTs systems with different configurations of hybrid systems by integrating the RM3 with the 5 MW and 15 MW FWTs

Additionally, Figure 50 illustrates the symbiotic relationship quantified through $P_{CV}$ between the FWT and the WEC in hybrid systems by integrating the RM3 with the 5 MW and 15 MW FWTs on a spar mooring platform with and without the reaction plate at the bottom and a semi-submersible mooring platform. These hybrid configurations are compared against standalone RM3, 5 MW FWTs, and 15 MW FWTs. The $P_{CV}$ of the WEC decreases across all different configurations, with reduction percentages ranging between 5% and 18%, demonstrating that the hybrid setup enhances the stability of wave energy output. For the WTs, the change in the $P_{CV}$ varies by configuration. For the 5 MW FWT, the $P_{CV}$ increases with the spar configuration without the reaction plate by about 3%, whereas it decreases with the reaction plate by about 2%., highlighting the improvement in performance due to the presence of the reaction plate. The $P_{CV}$ remains unchanged in the semi-submersible configuration. Consequently, the symbiotic relationship for the 5 MW depends on the configuration with no clear trend: mutualistic with the spar configuration with the reaction plate, commensalistic with the semi-submersible configuration, and parasitic with the spar configuration without the reaction plate. In contrast, the 15 MW FWT shows unchanged $P_{CV}$ for both spar configurations, with and without the reaction plate, and a 4% reduction with the semi-submersible configuration, corresponding to commensalistic and mutualistic interactions, respectively. Overall, the results indicate a predominantly commensalistic relationship, which can be further optimized toward mutualism. These findings demonstrate that hybridization consistently enhances WEC power stability while generally maintaining or slightly improving FWT stability, underscoring the potential of hybrid systems to provide more reliable and consistent energy generation.

## 7. Conclusions

A hybrid wind-wave system integrating the RM3 WEC and various FWTs was evaluated through numerical simulations using the coupled WEC-Sim, MOST, and MoorDyn models. The FWTs considered in this study include the NREL 5 MW WT on spar and semi-submersible platforms and the IEA 15 MW WT on similar platforms. Key findings include:

- **Reaction plate impact on AEP and LCOE**
  - The reaction plate plays a critical role in maximizing energy output.
  - Its removal led to substantial reductions in AEP across different water depths.
  - Extending the spar length while retaining the reaction plate improved energy production and reduced the LCOE by ~20%.
  - Cost savings from removing the reaction plate do not offset the resulting energy loss.
- **System dynamics**
  - The float exhibited increased motions in most directions, while the mooring platform experienced reduced surge and pitch responses, contributing to improved overall stability
  - Mooring forces showed mixed responses: reductions in surge and pitch forces were observed, while heave forces increased, indicating higher vertical loads on the system.
  - The addition of the reaction plate significantly reduces mooring loads and enhances stability for the 5 MW spar and semi-submersible platforms, particularly in surge, heave, and pitch under regular waves, with the most pronounced effects observed with the spar platform.
  - For the 15 MW systems, the reaction plate provides minimal improvement under both regular and irregular waves.
- **Power generation**
  - The power generation by the WT had a negligible change.
  - The power generation by the WEC was significantly increased, showing the hydrodynamic benefits of hybridization.
  - This dynamic interaction of the hybrid system has minimal impact on the aerodynamics while having a substantial impact on hydrodynamics.
  - Hybrid systems provide more stable power, with $P_{CV}$ about 50% lower than standalone systems.
- **Economic performance**
  - In a 100-unit hybrid farm, the WEC LCOE decreased by up to 83%, while the FWT experienced minimal changes in LCOE.
  - Hybrid systems with the NREL 5 MW WT showed overall LCOE reductions.
  - Hybrid systems with the IEA 15 MW WT showed a slight LCOE increase of ~4%.
  - $P_{CV}$ analysis confirmed enhanced stability of WEC power output (reductions of 5–18%), while FWTs showed configuration-dependent minor or negligible changes.
  - Hybridization demonstrates a predominantly *mutualistic relationship*: the WEC benefits significantly, while the FWT sees minor improvements or negligible impacts.
  - These results suggest that, with further system-level optimization, full mutualism could potentially be achieved.

While this study provides valuable insights into the performance and economic viability of hybrid FWT-WEC systems, it is important to acknowledge its limitations. The numerical simulations are based on linear wave theory and the WEC-Sim model, which relies on potential flow theory. These assumptions do not fully capture nonlinear hydrodynamic effects, such as viscous drag or wave-breaking phenomena, which are incorporated as external terms. Additionally, the control parameters of the WEC, represented by the damping coefficient of the PTO, remained constant and were not optimized. Although control co-design approaches could substantially enhance WEC performance, this aspect lies beyond the scope of the present work. Future work could employ more advanced modeling approaches to better account for hydrodynamic nonlinearities and refine the accuracy of the results and integrate control co-design strategies to further improve the power output by the WEC.

In conclusion, this study demonstrates the potential of hybrid offshore wind-wave systems to enhance energy efficiency and economic viability. By carefully considering critical structural components, such as the reaction plate, and optimizing system integration strategies, hybrid systems offer a promising pathway toward a more sustainable and cost-effective solution for offshore renewable energy generation.

## Acknowledgement

The authors gratefully acknowledge the financial support provided by the Cornell University Ezra Systems Scholars Program, which made this research possible through a fellowship award to Alaa Ahmed.

## A. Additional tables for detailed results

**Table A.1**

Difference in RMS of the float's heave response RMS between standalone and hybrid system for the six configurations under different regular wave conditions

| Configuration | Hybrid 5MW Spar | | | Hybrid 5MW Spar w/ RP | | | Hybrid 5MW SS | | | Hybrid 15MW Spar | | | Hybrid 15MW Spar w/ RP | | | Hybrid 15MW SS | | |
|---|---|---|---|---|---|---|---|---|---|---|---|---|---|---|---|---|---|---|
| Run number | Run 1 | Run 2 | Run 3 | Run 1 | Run 2 | Run 3 | Run 1 | Run 2 | Run 3 | Run 1 | Run 2 | Run 3 | Run 1 | Run 2 | Run 3 | Run 1 | Run 2 | Run 3 |
| Float (m) | 0.04 | 0.10 | 0.11 | 0.02 | -0.02 | -0.10 | 0.06 | 0.16 | 0.12 | -0.01 | -0.01 | -0.03 | 0.05 | -0.01 | -0.10 | 0.01 | 0.24 | 0.08 |

**Table A.2**

Difference in RMS of the float's heave response between standalone and hybrid system for the six configurations under different irregular wave conditions

| Configuration | Hybrid 5MW Spar | | | Hybrid 5MW Spar w/ RP | | | Hybrid 5MW SS | | | Hybrid 15MW Spar | | | Hybrid 15MW Spar w/ RP | | | Hybrid 15MW SS | | |
|---|---|---|---|---|---|---|---|---|---|---|---|---|---|---|---|---|---|---|
| Run number | Run 1 | Run 2 | Run 3 | Run 1 | Run 2 | Run 3 | Run 1 | Run 2 | Run 3 | Run 1 | Run 2 | Run 3 | Run 1 | Run 2 | Run 3 | Run 1 | Run 2 | Run 3 |
| Float (m) | 0 | 0.16 | -0.13 | 0 | -0.03 | -0.09 | 0.03 | 0.11 | 0.14 | -0.01 | 0 | 0 | -0.01 | -0.01 | -0.08 | 0.04 | 0.10 | 0.15 |

**Table A.3**

Difference in response RMS of the spar platform and mooring forces between standalone and hybrid system of 5MW WT and the float with spar platform with and without the reaction plate under different regular wave conditions with wind speed of 11.4 m/s

| Degree of Freedom | X | | | Y | | | Z | | | RX | | | RY | | | RZ | | |
|---|---|---|---|---|---|---|---|---|---|---|---|---|---|---|---|---|---|---|
| Run number | Run 1 | Run 2 | Run 3 | Run 1 | Run 2 | Run 3 | Run 1 | Run 2 | Run 3 | Run 1 | Run 2 | Run 3 | Run 1 | Run 2 | Run 3 | Run 1 | Run 2 | Run 3 |
| **Mooring platform w/o plate (m or deg)** | 0.15 | -0.23 | 0.19 | -0.01 | 0.01 | -0.05 | 0.03 | 0.18 | 0.50 | -0.02 | -0.02 | -0.37 | 0.07 | 0.03 | -0.42 | -0.02 | -0.04 | -0.41 |
| **Mooring platform w plate (m or deg)** | -0.67 | -0.42 | -0.24 | -0.17 | -0.13 | 0.02 | -0.03 | -0.06 | -0.18 | -0.09 | -0.04 | 0.21 | -0.19 | -0.13 | -0.11 | 0.03 | 0.15 | 0.76 |
| **Mooring Forces w/o plate (kNor kNm)** | 4.56 | -0.93 | 9.33 | -1.12 | 0.06 | -16.51 | 2.19 | 1.99 | 18.57 | -89.64 | -11.19 | -1578.57 | 434.75 | 154.78 | 460.12 | -20.49 | -10.90 | -340.59 |
| **Mooring Forces w plate (kNor kNm)** | -11.38 | -6.46 | -4.18 | -5.71 | -2.28 | 13.05 | -5.13 | -3.08 | -2.22 | -400.77 | -105.32 | 1105.19 | -748.53 | -347.92 | -250.76 | -73.00 | -37.01 | 69.28 |

**Table A.4**

Difference in response RMS of the spar platform and mooring forces between standalone and hybrid system of 5MW WT and the float with spar platform with and without the reaction plate under different irregular wave conditions with wind speed of 11.4 m/s

| Degree of Freedom | X | | | Y | | | Z | | | RX | | | RY | | | RZ | | |
|---|---|---|---|---|---|---|---|---|---|---|---|---|---|---|---|---|---|---|
| Run number | Run 1 | Run 2 | Run 3 | Run 1 | Run 2 | Run 3 | Run 1 | Run 2 | Run 3 | Run 1 | Run 2 | Run 3 | Run 1 | Run 2 | Run 3 | Run 1 | Run 2 | Run 3 |
| Mooring platform w/o plate (m or deg) | -3.54 | -3.76 | -3.60 | -0.85 | -0.81 | -0.78 | -0.06 | 0.00 | 0.07 | -0.22 | -0.21 | -0.21 | -0.46 | -0.54 | -0.79 | -0.12 | -0.17 | -0.24 |
| Mooring platform w plate (m or deg) | -0.13 | -0.36 | -0.20 | -0.07 | -0.05 | -0.05 | 0.00 | -0.03 | -0.07 | -0.04 | -0.04 | -0.04 | -0.17 | -0.16 | -0.13 | 0.05 | 0.04 | 0.10 |
| Mooring Forces w/o plate (kNor kNm) | -58.80 | -61.65 | -64.88 | -28.79 | -27.99 | -27.05 | -27.45 | -25.96 | -23.63 | -1889.06 | -1832.40 | -1783.88 | -3678.80 | -3890.44 | -4217.22 | -197.94 | -187.87 | -185.43 |
| Mooring Forces w plate (kNor kNm) | -0.45 | -4.05 | -1.59 | 1.10 | 1.35 | 0.84 | 0.42 | -1.42 | -0.55 | 123.47 | 142.81 | 101.26 | 68.98 | -151.60 | 16.22 | -40.14 | -35.09 | -35.81 |

**Table A.5**

Difference in response RMS of the semi-submersible platform and mooring forces between standalone and hybrid system of 5MW WT and the float with semi-submersible platform under different regular wave conditions with wind speed of 11.4 m/s

| Degree of Freedom | X | | | Y | | | Z | | | RX | | | RY | | | RZ | | |
|---|---|---|---|---|---|---|---|---|---|---|---|---|---|---|---|---|---|---|
| Run number | Run 1 | Run 2 | Run 3 | Run 1 | Run 2 | Run 3 | Run 1 | Run 2 | Run 3 | Run 1 | Run 2 | Run 3 | Run 1 | Run 2 | Run 3 | Run 1 | Run 2 | Run 3 |
| **Mooring platform (m or deg)** | -0.04 | 0.01 | 0.00 | 0.00 | 0.00 | -0.01 | 0.02 | -0.05 | 0.40 | -0.03 | -0.07 | -0.08 | -0.05 | -0.06 | -0.07 | 0.00 | -0.01 | -0.02 |
| **Mooring Forces (kNor kNm)** | -3.21 | 9.56 | 58.37 | -0.17 | -0.10 | -0.19 | 2.71 | 5.04 | 59.40 | -32.52 | -74.93 | -77.35 | -66.04 | -69.53 | 292.52 | -1.58 | -20.23 | -23.28 |

**Table A.6**

Difference in response RMS of the semi-submersible platform and mooring forces between standalone and hybrid system of 5MW WT and the float with semi-submersible platform under different irregular wave conditions with wind speed of 11.4 m/s

| Degree of Freedom | X | | | Y | | | Z | | | RX | | | RY | | | RZ | | |
|---|---|---|---|---|---|---|---|---|---|---|---|---|---|---|---|---|---|---|
| Run number | Run 1 | Run 2 | Run 3 | Run 1 | Run 2 | Run 3 | Run 1 | Run 2 | Run 3 | Run 1 | Run 2 | Run 3 | Run 1 | Run 2 | Run 3 | Run 1 | Run 2 | Run 3 |
| **Mooring platform (m or deg)** | 0.02 | 0.01 | 0.01 | -0.03 | -0.02 | -0.02 | -0.01 | 0.04 | -0.06 | -0.14 | -0.13 | -0.07 | -0.02 | -0.03 | -0.05 | 0.00 | 0.00 | 0.03 |
| **Mooring Forces (kN or kNm)** | 1.73 | 2.42 | 12.41 | -1.30 | -0.89 | -0.76 | 0.94 | 3.14 | 8.56 | -139.51 | -133.16 | -64.35 | -20.02 | -34.56 | -16.44 | 5.17 | 1.73 | 52.34 |

**Table A.7**

Difference in response RMS of the spar platform and mooring forces between standalone and hybrid system of 15MW WT and the float with spar platform with and without the reaction plate under different regular wave conditions with wind speed of 10.6 m/s

| Degree of Freedom | X | | | Y | | | Z | | | RX | | | RY | | | RZ | | |
|---|---|---|---|---|---|---|---|---|---|---|---|---|---|---|---|---|---|---|
| Run number | Run 1 | Run 2 | Run 3 | Run 1 | Run 2 | Run 3 | Run 1 | Run 2 | Run 3 | Run 1 | Run 2 | Run 3 | Run 1 | Run 2 | Run 3 | Run 1 | Run 2 | Run 3 |
| **Mooring platform w/o plate (m or deg)** | -0.18 | -0.21 | -0.31 | -0.01 | -0.14 | -0.06 | -0.18 | -0.13 | -0.19 | -0.02 | -0.05 | -0.16 | -0.07 | -0.07 | -0.25 | 0.00 | 0.02 | -0.01 |
| **Mooring platform w plate (m or deg)** | -0.35 | -0.29 | -0.16 | -0.22 | -0.05 | -0.09 | -0.17 | -0.16 | -0.35 | -0.04 | -0.04 | 0.02 | -0.11 | -0.06 | -0.03 | 0.00 | 0.03 | 0.07 |
| **Mooring Forces w/o plate (kNor kNm)** | -16.89 | -23.74 | -38.83 | 0.12 | -0.79 | -0.70 | -3.81 | -6.32 | -11.17 | -0.85 | -79.67 | -155.55 | -1068.43 | -1513.83 | -2575.62 | -25.14 | -76.34 | -274.65 |
| **Mooring Forces w plate (kNor kNm)** | -40.00 | -38.45 | -29.06 | -0.38 | 0.19 | 0.11 | -7.94 | -8.88 | -8.84 | -46.23 | -24.75 | 47.18 | -2598.89 | -2470.59 | -1842.56 | -53.89 | -68.05 | 28.00 |

**Table A.8**
Difference in response RMS of the spar platform and mooring forces between standalone and hybrid system of 15MW WT and the float with spar platform with and without the reaction plate under different irregular wave conditions with wind speed of 10.6 m/s

| Degree of Freedom | X | | | Y | | | Z | | | RX | | | RY | | | RZ | | |
|---|---|---|---|---|---|---|---|---|---|---|---|---|---|---|---|---|---|---|
| Run number | Run 1 | Run 2 | Run 3 | Run 1 | Run 2 | Run 3 | Run 1 | Run 2 | Run 3 | Run 1 | Run 2 | Run 3 | Run 1 | Run 2 | Run 3 | Run 1 | Run 2 | Run 3 |
| **Mooring platform w/o plate (m or deg)** | -0.78 | -0.80 | -0.69 | 0.53 | 0.47 | 0.51 | -0.32 | -0.31 | -0.31 | -0.16 | -0.16 | -0.08 | -0.12 | -0.13 | -0.19 | -0.01 | -0.01 | 0.00 |
| **Mooring platform w plate (m or deg)** | -0.32 | -0.35 | -0.26 | 0.02 | -0.03 | 0.02 | -0.29 | -0.29 | -0.38 | -0.09 | -0.07 | 0.02 | -0.15 | -0.13 | -0.09 | 0.00 | -0.01 | 0.02 |
| **Mooring Forces w/o plate (kNor kNm)** | -74.01 | -75.15 | -72.80 | 2.84 | 2.53 | 2.73 | -14.02 | -14.56 | -14.98 | 123.96 | 103.76 | 152.90 | -4620.75 | -4716.43 | -4621.09 | -256.72 | -263.12 | -130.70 |
| **Mooring Forces w plate (kNor kNm)** | -25.46 | -28.36 | -21.64 | 1.86 | 1.58 | 1.80 | -5.47 | -6.27 | -5.37 | 100.54 | 89.55 | 145.23 | -1305.82 | -1514.39 | -1099.30 | -146.40 | -125.65 | 30.68 |

**Table A.9**
Difference in response RMS of the semi-submersible platform and mooring forces between standalone and hybrid system of 15MW WT and the float with semi-submersible platform under different regular wave conditions with wind speed of 10.6 m/s

| Degree of Freedom | X | | | Y | | | Z | | | RX | | | RY | | | RZ | | |
|---|---|---|---|---|---|---|---|---|---|---|---|---|---|---|---|---|---|---|
| Run number | Run 1 | Run 2 | Run 3 | Run 1 | Run 2 | Run 3 | Run 1 | Run 2 | Run 3 | Run 1 | Run 2 | Run 3 | Run 1 | Run 2 | Run 3 | Run 1 | Run 2 | Run 3 |
| **Mooring platform (m or deg)** | -0.01 | 0.02 | 0.04 | 0.00 | 0.01 | 0.01 | 0.01 | 0.02 | -0.05 | 0.00 | 0.00 | 0.02 | 0.00 | -0.03 | -0.06 | 0.00 | 0.00 | 0.00 |
| **Mooring Forces (kNor kNm)** | -8.16 | -1.59 | 13.39 | -0.13 | 1.09 | 0.82 | -0.57 | 9.61 | -5.80 | -4.07 | -4.10 | 93.52 | -251.94 | -147.02 | -216.79 | 12.79 | -33.86 | 40.07 |

**Table A.10**
Difference in response RMS of the semi-submersible platform and mooring forces between standalone and hybrid system of 15MW WT and the float with semi-submersible platform under different irregular wave conditions with wind speed of 10.6 m/s

| Degree of Freedom | X | | | Y | | | Z | | | RX | | | RY | | | RZ | | |
|---|---|---|---|---|---|---|---|---|---|---|---|---|---|---|---|---|---|---|
| Run number | Run 1 | Run 2 | Run 3 | Run 1 | Run 2 | Run 3 | Run 1 | Run 2 | Run 3 | Run 1 | Run 2 | Run 3 | Run 1 | Run 2 | Run 3 | Run 1 | Run 2 | Run 3 |
| **Mooring platform (m or deg)** | -0.12 | -0.13 | -0.08 | 0.07 | 0.08 | 0.07 | -0.11 | -0.03 | -0.09 | 0.08 | 0.02 | 0.01 | -0.01 | -0.02 | -0.04 | 0.02 | 0.02 | 0.00 |
| **Mooring Forces (kNor kNm)** | -23.42 | -27.41 | -20.71 | 8.98 | 11.37 | 9.53 | -0.73 | 2.47 | -6.92 | 343.51 | 199.72 | 122.40 | -352.98 | -494.77 | -422.01 | 137.32 | 158.77 | 8.53 |

**Table A.11**
Difference in RMS of power between standalone and hybrid system for different configurations under different regular wave conditions

| Configuration | Hybrid 5MW Spar | | | Hybrid 5MW Spar w/ RP | | | Hybrid 5MW SS | | | Hybrid 15MW Spar | | | Hybrid 15MW Spar w/ RP | | | Hybrid 15MW SS | | |
|---|---|---|---|---|---|---|---|---|---|---|---|---|---|---|---|---|---|---|
| Run number | Run 1 | Run 2 | Run 3 | Run 1 | Run 2 | Run 3 | Run 1 | Run 2 | Run 3 | Run 1 | Run 2 | Run 3 | Run 1 | Run 2 | Run 3 | Run 1 | Run 2 | Run 3 |
| Power WT (MW) | 0.02 | 0.01 | 0 | 0.04 | 0.05 | 0.04 | 0.01 | 0.01 | 0.01 | 0.07 | 0.01 | -0.23 | 0.12 | 0.14 | 0.12 | 0 | -0.01 | -0.04 |
| Power WEC (kW) | 41.94 | 259.60 | 369.22 | 35.41 | 263.32 | 515.48 | 47.70 | 222.46 | 584.18 | 10.11 | 216.79 | 392.80 | 60.21 | 262.96 | 470.71 | 86.85 | 380.06 | 320.31 |

**Table A.12**
Difference in RMS of power between standalone and hybrid system for different configurations under different irregular wave conditions

| Configuration | Hybrid 5MW Spar | | | Hybrid 5MW Spar w/ RP | | | Hybrid 5MW SS | | | Hybrid 15MW Spar | | | Hybrid 15MW Spar w/ RP | | | Hybrid 15MW SS | | |
|---|---|---|---|---|---|---|---|---|---|---|---|---|---|---|---|---|---|---|
| Run number | Run 1 | Run 2 | Run 3 | Run 1 | Run 2 | Run 3 | Run 1 | Run 2 | Run 3 | Run 1 | Run 2 | Run 3 | Run 1 | Run 2 | Run 3 | Run 1 | Run 2 | Run 3 |
| Power WT (MW) | -3.52 | -3.51 | -3.52 | 0.06 | 0.07 | 0.05 | 0.01 | 0.01 | 0.01 | -0.06 | -0.08 | -0.09 | 0.09 | 0.09 | 0.10 | 0.03 | 0.01 | -0.01 |
| Power WEC (kW) | 8.59 | 36.57 | 253.69 | 3.83 | 70.42 | 275.59 | 6.19 | 115.83 | 295.85 | -5.69 | 59.31 | 240.54 | -1.17 | 73.23 | 272.49 | 35.60 | 138.52 | 362.85 |

**Table A.13**
Power coefficient of variation for standalone and hybrid systems of different configurations

| | Standalone WEC | Standalone 5MW WT | | Standalone 15MW WT | | Hybrid 5MW WT | | | Hybrid 15MW WT | | |
|---|---|---|---|---|---|---|---|---|---|---|---|
| | Spar platform | Spar platform | Semi-submersible platform | Spar platform | Semi-submersible platform | Spar platform | Spar platform with plate | Semi-submersible platform | Spar platform | Spar platform with plate | Semi-submersible platform |
| $P_{CV}$ | 0.76 | 1.16 | 1.18 | 1.11 | 1.10 | 0.68 | 0.67 | 0.65 | 0.7 | 0.68 | 0.60 |

**Table A.14**
AEP and CF of all configurations at the Humboldt Site

| System Type | Standalone WEC | Standalone 5 MW WT | | Standalone 15 MW WT | | Hybrid 5 MW WT | | | | | | Hybrid 15 MW WT | | | | | |
|---|---|---|---|---|---|---|---|---|---|---|---|---|---|---|---|---|---|
| | Spar platform | Spar platform | Semi-submersible platform | Spar platform | Semi-submersible platform | Spar platform | | Spar platform with plate | | Semi-submersible platform | | Spar platform | | Spar platform with plate | | Semi-submersible platform | |
| | | | | | | WEC | WT | WEC | WT | WEC | WT | WEC | WT | WEC | WT | WEC | WT |
| AEP (GWh) | 0.51 | 8.42 | 8.28 | 27.8 | 28.37 | 0.65 | 8.25 | 0.66 | 8.87 | 0.71 | 8.35 | 0.56 | 27.64 | 0.62 | 27.98 | 0.80 | 28.40 |
| CF | 0.22 | 0.21 | 0.20 | 0.23 | 0.23 | 0.27 | 0.20 | 0.28 | 0.22 | 0.34 | 0.21 | 0.24 | 0.23 | 0.27 | 0.23 | 0.34 | 0.23 |

**Table A.15**
AEP and CF of standalone configurations at Kodiak Site

| System Type | Standalone 5 MW WT | | Standalone 15 MW WT | |
|---|---|---|---|---|
| | Spar platform | Semi-submersible platform | Spar platform | Semi-submersible platform |
| AEP (GWh) | 15.11 | 14.94 | 50.91 | 51.97 |
| CF | 0.37 | 0.37 | 0.42 | 0.43 |

**Table A.16**
LCOE of a single unit within a 100-unit farm of different configurations

| System Type | Standalone WEC | Standalone 5 MW WT | | Standalone 15 MW WT | | Hybrid 5 MW WT | | | Hybrid 15 MW WT | | |
|---|---|---|---|---|---|---|---|---|---|---|---|
| | Spar platform | Spar platform | Semi-submersible platform | Spar platform | Semi-submersible platform | Spar platform | Spar platform with plate | Semi-submersible platform | Spar platform | Spar platform with plate | Semi-submersible platform |
| Total AEP (GWh) | 0.51 | 8.42 | 8.28 | 27.80 | 28.37 | 9.08 | 9.53 | 9.06 | 28.20 | 28.59 | 29.19 |
| CapEx (B$) | 0.53 | 2.29 | 2.78 | 6.51 | 7.54 | 2.43 | 2.65 | 2.92 | 6.82 | 7.17 | 7.85 |
| OpEx (M$) | 12.6 | 91.5 | 91.5 | 275 | 275 | 91.5 | 91.5 | 91.5 | 275 | 275 | 275 |
| LCOE ($/MWh) | 1373 | 407.9 | 479.84 | 356.44 | 389.18 | 395.39 | 401.67 | 455.73 | 363.34 | 371.68 | 389.79 |

**Table A.17**
Comparison between the LCOE for a 100 standalone RM3 farm and 100 standalone 5 MW and 15 MW WTs farms with a 100 hybrid system farm by integrating the RM3 with the 5 MW and 15 MW FWTs

| System Type | Standalone WEC | Standalone 5 MW WT | | Standalone 15 MW WT | | Hybrid 5 MW WT | | | | | | Hybrid 15 MW WT | | | | | |
|---|---|---|---|---|---|---|---|---|---|---|---|---|---|---|---|---|---|
| | Spar platform | Spar platform | Semi-submersible platform | Spar platform | Semi-submersible platform | Spar platform | | Spar platform with plate | | Semi-submersible platform | | Spar platform | | Spar platform with plate | | Semi-submersible platform | |
| | | | | | | WEC | WT | WEC | WT | WEC | WT | WEC | WT | WEC | WT | WEC | WT |
| LCOE ($/MWh) | 1373 | 407.90 | 479.84 | 356.44 | 389.18 | 235.79 | 407.90 | 598.48 | 387.04 | 219.87 | 475.73 | 599.36 | 358.55 | 1162.4 | 354.19 | 423.23 | 388.85 |

**Table A.18**
Comparison between the $P_{CV}$ for RM3, standalone 5 MW, and 15 MW WTs systems with different configurations of hybrid systems by integrating the RM3 with the 5 MW and 15 MW FWTs

| System Type | Standalone WEC | Standalone 5 MW WT | | Standalone 15 MW WT | | Hybrid 5 MW WT | | | | | | Hybrid 15 MW WT | | | | | |
|---|---|---|---|---|---|---|---|---|---|---|---|---|---|---|---|---|---|
| | Spar platform | Spar platform | Semi-submersible platform | Spar platform | Semi-submersible platform | Spar platform | | Spar platform with plate | | Semi-submersible platform | | Spar platform | | Spar platform with plate | | Semi-submersible platform | |
| | | | | | | WEC | WT | WEC | WT | WEC | WT | WEC | WT | WEC | WT | WEC | WT |
| $P_{CV}$ | 0.76 | 1.16 | 1.18 | 1.11 | 1.1 | 0.67 | 1.19 | 0.68 | 1.14 | 0.65 | 1.18 | 0.72 | 1.11 | 0.7 | 1.11 | 0.62 | 1.06 |

# References

Beiter, P., Musial, W., Smith, A., Kilcher, L., Damiani, R., Maness, M., Sirnivas, S., Stehly, T., Gevorgian, V., Mooney, M., et al., 2016. A spatial-economic cost-reduction pathway analysis for US offshore wind energy development from 2015–2030. Technical Report. National Renewable Energy Lab.(NREL), Golden, CO (United States).

Clark, C.E., Paredes, G.M., 2018. Effects of co-located floating wind-wave systems on fatigue damage of floating offshore wind turbine mooring cables, in: International Conference on Offshore Mechanics and Arctic Engineering, American Society of Mechanical Engineers. p. V001T01A024.

Cottura, L., Caradonna, R., Ghigo, A., Novo, R., Bracco, G., Mattiazzo, G., 2021. Dynamic modeling of an offshore floating wind turbine for application in the mediterranean sea. Energies 14, 248.

Dabar, O.A., Awaleh, M.O., Waberi, M.M., Ghiasirad, H., Adan, A.B.I., Ahmed, M.M., Nasser, M., Juangsa, F.B., Guirreh, I.A., Abdillahi, M.O., et al., 2024. Techno-economic and environmental assessment of green hydrogen and ammonia production from solar and wind energy in the republic of djibouti: A geospatial modeling approach. Energy Reports 12, 3671–3689.

Dallman, A.R., Neary, V.S., 2014. Characterization of US Wave Energy Converter (WEC) Test Sites: A Catalogue of Met-Ocean Data. Technical Report. Sandia National Lab.(SNL-NM), Albuquerque, NM (United States).

Feng, X., Lin, Y., Gu, Y., Li, D., Chen, B., Liu, H., Sun, Y., 2023. Preliminary stability design method and hybrid experimental validation of a floating platform for 10 mw wind turbine. Ocean Engineering 285, 115401.

Gaertner, E., Rinker, J., Sethuraman, L., Zahle, F., Anderson, B., Barter, G., Abbas, N., Meng, F., Bortolotti, P., Skrzypinski, W., Scott, G., Feil, R., Bredmose, H., Dykes, Shields, K.M., Allen, C., Viselli, A., 2020. Definition of the IEA 15-Megawatt Offshore Reference Wind. Technical Report. National Renewable Energy Lab.(NREL), Golden, CO (United States).

Gao, Q., Bechlenberg, A., Jayawardhana, B., Ertugrul, N., Vakis, A.I., Ding, B., 2024. Techno-economic assessment of offshore wind and hybrid wind–wave farms with energy storage systems. Renewable and Sustainable Energy Reviews 192, 114263.

Gao, Q., Khan, S.S., Sergiienko, N., Ertugrul, N., Hemer, M., Negnevitsky, M., Ding, B., 2022. Assessment of wind and wave power characteristic and potential for hybrid exploration in australia. Renewable and Sustainable Energy Reviews 168, 112747.

Garcia, M., Oliva, S., 2023. Technical, economic, and co2 emissions assessment of green hydrogen production from solar/wind energy: the case of chile. Energy 278, 127981.

Garcia-Sanz, M., 2024. Hybrid energy systems: Synergy margin and control co-design. Advanced Control for Applications: Engineering and Industrial Systems 6, e238.

Hall, M., 2020. Moordyn v2: New capabilities in mooring system components and load cases, in: International Conference on Offshore Mechanics and Arctic Engineering, American Society of Mechanical Engineers. p. V009T09A078.

He, G., Zhao, C., Liu, C., He, R., Luan, Z., 2024. Power absorption and dynamic response analysis of a hybrid system with a semi-submersible wind turbine and a salter's duck wave energy converter array. Energy 305, 132210.

Ibáñez-Rioja, A., Järvinen, L., Puranen, P., Kosonen, A., Ruuskanen, V., Hynynen, K., Ahola, J., Kauranen, P., 2023. Off-grid solar pv–wind power–battery–water electrolyzer plant: Simultaneous optimization of component capacities and system control. Applied energy 345, 121277.

Jonkman, B., et al., 2014. Turbsim user's guide v2. 00.00. Natl. Renew. Energy Lab .

Jonkman, J., 2010. Definition of the Floating System for Phase IV of OC3. Technical Report. National Renewable Energy Lab.(NREL), Golden, CO (United States).

Jonkman, J., Butterfield, S., Musial, W., Scott, G., 2009. Definition of a 5-MW reference wind turbine for offshore system development. Technical Report. National Renewable Energy Lab.(NREL), Golden, CO (United States).

Khan, T., Yu, M., Waseem, M., 2022. Review on recent optimization strategies for hybrid renewable energy system with hydrogen technologies: State of the art, trends and future directions. International Journal of Hydrogen Energy 47, 25155–25201.

Kluger, J.M., 2017. Synergistic design of a combined floating wind turbine-wave energy converter. Ph.D. thesis. Massachusetts Institute of Technology.

Kluger, J.M., Haji, M.N., Slocum, A.H., 2023. The power balancing benefits of wave energy converters in offshore wind-wave farms with energy storage. Applied Energy 331, 120389.

Kluger, J.M., Slocum, A.H., Sapsis, T.P., 2017. A first-order dynamics and cost comparison of wave energy converters combined with floating wind turbines, in: ISOPE International Ocean and Polar Engineering Conference, ISOPE. pp. ISOPE–I.

Li, M., Luo, H., Zhou, S., Kumar, G.M.S., Guo, X., Law, T.C., Cao, S., 2022. State-of-the-art review of the flexibility and feasibility of emerging offshore and coastal ocean energy technologies in east and southeast asia. Renewable and Sustainable Energy Reviews 162, 112404.

Lou, J., Robertson, B., 2024. Pacwave cotemporal resource assessment of waves, winds, and currents .

Martinez, A., Iglesias, G., 2022. Mapping of the levelised cost of energy for floating offshore wind in the european atlantic. Renewable and Sustainable Energy Reviews 154, 111889.

Meng, F., Sergiienko, N., Ding, B., Zhou, B., Da Silva, L.S.P., Cazzolato, B., Li, Y., 2023. Co-located offshore wind–wave energy systems: Can motion suppression and reliable power generation be achieved simultaneously? Applied Energy 331, 120373.

Moghimi, M., Derakhshan, S., Motawej, H., 2020. A mathematical model development for assessing the engineering and economic improvement of wave and wind hybrid energy system. Iranian Journal of Science and Technology, Transactions of Mechanical Engineering 44, 507–521.

Myhr, A., Bjerkseter, C., Ågotnes, A., Nygaard, T.A., 2014. Levelised cost of energy for offshore floating wind turbines in a life cycle perspective. Renewable energy 66, 714–728.

National Data Buoy Center, 2024a. Station 46001 - kodiak, ak URL: https://www.ndbc.noaa.gov/station_page.php?station=46001. accessed: April 27, 2025.

National Data Buoy Center, 2024b. Station 46022 - eureka, ca (166) URL: https://www.ndbc.noaa.gov/station_page.php?station=46022. accessed: April 27, 2025.

Neary, V.S., Lawson, M., Previsic, M., Copping, A., Hallett, K.C., Labonte, A., Rieks, J., Murray, D., 2014. Methodology for design and economic analysis of marine energy conversion (mec) technologies .

Neshat, M., Sergiienko, N.Y., Nezhad, M.M., da Silva, L.S., Amini, E., Marsooli, R., Garcia, D.A., Mirjalili, S., 2024. Enhancing the performance of hybrid wave-wind energy systems through a fast and adaptive chaotic multi-objective swarm optimisation method. Applied Energy 362, 122955.

Ringwood, J.V., 2025. Control co-design for wave energy systems. Applied Ocean Research 158, 104514.

Robertson, A., Jonkman, J., Masciola, M., Song, H., Goupee, A., Coulling, A., Luan, C., 2014. Definition of the semisubmersible floating system for phase II of OC4. Technical Report. National Renewable Energy Lab.(NREL), Golden, CO (United States).

Ruehl, K., Keester, A., Forbush, D., Grasberger, J., Husain, S., Leon, J., Ogden, D., Shabara, M., 2024. WEC-Sim v6.1. Technical Report. National Renewable Energy Lab.(NREL), Golden, CO (United States).

Schupp, M.F., Bocci, M., Depellegrin, D., Kafas, A., Kyriazi, Z., Lukic, I., Schultz-Zehden, A., Krause, G., Onyango, V., Buck, B.H., 2019. Toward a common understanding of ocean multi-use. Frontiers in Marine Science 6, 165.

Sergiienko, N.Y., Xue, L., da Silva, L.S., Ding, B., Cazzolato, B.S., 2025. Statistical analysis of floating hybrid wind–wave energy systems. Applied Energy 401, 126669.

Ströfer, C.A.M., Gaebele, D.T., Coe, R.G., Bacelli, G., 2023. Control co-design of power take-off systems for wave energy converters using wecopttool. IEEE Transactions on Sustainable Energy 14, 2157–2167.

Tang, O., Rehme, J., Cerin, P., 2022. Levelized cost of hydrogen for refueling stations with solar pv and wind in sweden: On-grid or off-grid? Energy 241, 122906.

Wei, Z., Cao, F., Cao, C., Han, Z., Shi, H., Ji, T., 2025. Experimental study on the effects of an array of concentric wave energy converters on the dynamic of semi-submersible floating wind turbine. Renewable Energy , 122495.

Xu, S., Ji, B., Xu, F., Chen, C., 2024. Dynamic response and power performance of a combined semi-submersible floating wind turbine and point absorber wave energy converter array. Renewable Energy 237, 121903.

Zhao, H., Stansby, P., Liao, Z., Li, G., 2024. Multi-objective optimal control of a hybrid offshore wind turbine platform integrated with multi-float wave energy converters. Energy 312, 133547.

Zhu, K., Shi, H., Zheng, S., Michele, S., Cao, F., 2023. Hydrodynamic analysis of hybrid system with wind turbine and wave energy converter. Applied Energy 350, 121745.

Zhu, X., Chen, J., Wang, L., Wang, P., Fang, Z., Du, X., Gao, F., 2024. Optimum power control and coordinate sizing for the stand-alone wind-energy storage integrated hydrogen production system. International Journal of Hydrogen Energy 88, 803–814.

Zou, S., Robertson, B., Roach, A., Mundon, T., Rosenberg, B., Penalba, M., 2024. Wave energy converter arrays: A methodology to assess performance considering the disturbed wave field. Renewable Energy 229, 120719.